\documentclass[11pt]{article}
\usepackage{amsmath,amscd,amssymb}
\usepackage[english]{babel}
\usepackage[applemac]{inputenc}
\usepackage{graphicx}
\newtheorem{Lem}{Lemma \thesection.}
\newtheorem{Th}[Lem]{Theorem \thesection.}
\newtheorem{Cor}[Lem]{Corollary \thesection.}
\newtheorem{Def}[Lem]{Definition \thesection.}
\newtheorem{Prop  et Def}[Lem]{Proposition et Definition \thesection.}
\newtheorem{Ex}[Lem]{Examples \thesection.}
\newtheorem{Prop}[Lem]{Proposition \thesection.}
\newtheorem{Rem}[Lem]{Remark \thesection.}
\newtheorem{Lem and Def}[Lem]{Lemma and Definition \thesection.}
\newtheorem{Not}[Lem]{Notations \thesection.}
\def\cal{\mathcal}
\def\bb{\mathbb} 
\def\frak{\mathfrak }
\def\a{\alpha }
\def\b{\beta }

\def\d{\delta }
\def\D{\Delta }

\def\G{\Gamma }
\def\k{\kappa }
\def\l{\lambda }
\def\m{\mu }
\def\o{\omega }
\def\O{\Omega }

\def\r{\rho }
\def\s{\sigma }

\def\t{\theta }

\def\f{\varphi }
\def\dim{\rm dim\; }

\def\dps{\displaystyle }

\def\Card{{\rm Card}\, }
\def\deg{{\rm deg\ }}
\def\det{{\rm det}\ }

\newcommand{\tlowername}[2]%
{$\stackrel{\makebox[1pt]{#1}}%
{\begin{picture}(0,0)%
\put(0,0){\makebox(0,6)[t]{\makebox[1pt]{$#2$}}}%
\end{picture}}$}%

\newcommand{\AR}[1]%
{\begin{picture}(#1,0)%
\put(0,0){\vector(1,0){#1}}%
\end{picture}}%

\newcommand{\DOTAR}[1]%
{\NUMBEROFDOTS=#1%
\divide\NUMBEROFDOTS by 3%
\begin{picture}(#1,0)%
\multiput(0,0)(3,0){\NUMBEROFDOTS}{\circle*{1}}%
\put(#1,0){\vector(1,0){0}}%
\end{picture}}%

\newcommand{\MONO}[1]%
{\begin{picture}(#1,0)%
\put(0,0){\vector(1,0){#1}}%
\put(2,-2){\line(0,1){4}}%
\end{picture}}%

\newcommand{\EPI}[1]%
{\begin{picture}(#1,0)(-#1,0)%
\put(-#1,0){\vector(1,0){#1}}%
\put(-6,-2){\line(0,1){4}}%
\end{picture}}%

\newcommand{\BIMO}[1]%
{\begin{picture}(#1,0)(-#1,0)%
\put(-#1,0){\vector(1,0){#1}}%
\put(-6,-2){\line(0,1){4}}%
\put(-#1,-2){\hspace{2pt}\line(0,1){4}}%
\end{picture}}%

\newcommand{\BIAR}[1]%
{\begin{picture}(#1,4)%
\put(0,0){\vector(1,0){#1}}%
\put(0,4){\vector(1,0){#1}}%
\end{picture}}%

\newcommand{\EQL}[1]%
{\begin{picture}(#1,0)%
\put(0,1){\line(1,0){#1}}%
\put(0,-1){\line(1,0){#1}}%
\end{picture}}%

\newcommand{\ADJAR}[1]%
{\begin{picture}(#1,4)%
\put(0,0){\vector(1,0){#1}}%
\put(#1,4){\vector(-1,0){#1}}%
\end{picture}}%

\newcommand{\BKAR}[1]%
{\begin{picture}(#1,0)%
\put(#1,0){\vector(-1,0){#1}}%
\end{picture}}%

\newcommand{\BKDOTAR}[1]%
{\NUMBEROFDOTS=#1%
\divide\NUMBEROFDOTS by 3%
\begin{picture}(#1,0)%
\multiput(#1,0)(-3,0){\NUMBEROFDOTS}{\circle*{1}}%
\put(0,0){\vector(-1,0){0}}%
\end{picture}}%

\newcommand{\BKMONO}[1]%
{\begin{picture}(#1,0)(-#1,0)%
\put(0,0){\vector(-1,0){#1}}%
\put(-2,-2){\line(0,1){4}}%
\end{picture}}%

\newcommand{\BKEPI}[1]%
{\begin{picture}(#1,0)%
\put(#1,0){\vector(-1,0){#1}}%
\put(6,-2){\line(0,1){4}}%
\end{picture}}%

\newcommand{\BKBIMO}[1]%
{\begin{picture}(#1,0)%
\put(#1,0){\vector(-1,0){#1}}%
\put(6,-2){\line(0,1){4}}%
\put(#1,-2){\hspace{-2pt}\line(0,1){4}}%
\end{picture}}%

\newcommand{\BKBIAR}[1]%
{\begin{picture}(#1,4)%
\put(#1,0){\vector(-1,0){#1}}%
\put(#1,4){\vector(-1,0){#1}}%
\end{picture}}%

\newcommand{\BKADJAR}[1]%
{\begin{picture}(#1,4)%
\put(0,4){\vector(1,0){#1}}%
\put(#1,0){\vector(-1,0){#1}}%
\end{picture}}%

\newcommand{\lowername}[2]%
{$\stackrel{\makebox[1pt]{#1}}%
{\begin{picture}(0,0)%
\truex{600}%
\put(0,0){\makebox(0,\value{x})[t]{\makebox[1pt]{$#2$}}}%
\end{picture}}$}%

\newcommand{\hcase}[2]%
{\makebox[0pt]%
{\raisebox{-1pt}[0pt][0pt]{#1{#2}}}}%

\newcommand{\Hcase}[3]%
{\makebox[0pt]
{\raisebox{-1pt}[0pt][0pt]%
{$\stackrel{\makebox[0pt]{$\textstyle{#2}$}}{#1{#3}}$}}}%

\newcommand{\hcasE}[3]%
{\makebox[0pt]%
{\raisebox{-9pt}[0pt][0pt]%
{\lowername{#1{#3}}{#2}}}}%

\newcommand{\hbicase}[2]%
{\makebox[0pt]%
{\raisebox{-2.5pt}[0pt][0pt]{#1{#2}}}}%

\newcommand{\Hbicase}[4]%
{\makebox[0pt]
{\raisebox{-10.5pt}[0pt][0pt]%
{$\stackrel{\makebox[0pt]{$\textstyle{#2}$}}%
{\mbox{\lowername{#1{#4}}{#3}}}$}}}%

\newcommand{\EAR}[1]%
{\begin{picture}(#1,0)%
\put(0,0){\vector(1,0){#1}}%
\end{picture}}%

\newcommand{\EDOTAR}[1]%
{\truex{100}\truey{300}%
\NUMBEROFDOTS=#1%
\divide\NUMBEROFDOTS by \value{y}%
\begin{picture}(#1,0)%
\multiput(0,0)(\value{y},0){\NUMBEROFDOTS}%
{\circle*{\value{x}}}%
\put(#1,0){\vector(1,0){0}}%
\end{picture}}%

\newcommand{\EMONO}[1]%
{\begin{picture}(#1,0)%
\put(0,0){\vector(1,0){#1}}%
\truex{300}\truey{600}%
\put(\value{x},-\value{x}){\line(0,1){\value{y}}}%
\end{picture}}%

\newcommand{\EEPI}[1]%
{\begin{picture}(#1,0)(-#1,0)%
\put(-#1,0){\vector(1,0){#1}}%
\truex{300}\truey{600}\truez{800}%
\put(-\value{z},-\value{x}){\line(0,1){\value{y}}}%
\end{picture}}%

\newcommand{\EBIMO}[1]%
{\begin{picture}(#1,0)(-#1,0)%
\put(-#1,0){\vector(1,0){#1}}%
\truex{300}\truey{600}\truez{800}%
\put(-\value{z},-\value{x}){\line(0,1){\value{y}}}%
\put(-#1,-\value{x}){\hspace{3pt}\line(0,1){\value{y}}}%
\end{picture}}%

\newcommand{\EBIAR}[1]%
{\truex{400}%
\begin{picture}(#1,\value{x})%
\put(0,0){\vector(1,0){#1}}%
\put(0,\value{x}){\vector(1,0){#1}}%
\end{picture}}%

\newcommand{\EEQL}[1]%
{\begin{picture}(#1,0)%
\truex{200}%
\put(0,\value{x}){\line(1,0){#1}}%
\put(0,0){\line(1,0){#1}}%
\end{picture}}%

\newcommand{\EADJAR}[1]%
{\truex{400}%
\begin{picture}(#1,\value{x})%
\put(0,0){\vector(1,0){#1}}%
\put(#1,\value{x}){\vector(-1,0){#1}}%
\end{picture}}%

\newcommand{\ear}%
{\hspace{\SOURCE\unitlength}%
\hcase{\EAR}{\ARROWLENGTH}}%

\newcommand{\Ear}[1]%
{\hspace{\SOURCE\unitlength}%
\Hcase{\EAR}{#1}{\ARROWLENGTH}}%

\newcommand{\eaR}[1]%
{\hspace{\SOURCE\unitlength}%
\hcasE{\EAR}{#1}{\ARROWLENGTH}}%

\newcommand{\edotar}%
{\hspace{\SOURCE\unitlength}%
\hcase{\EDOTAR}{\ARROWLENGTH}}%

\newcommand{\Edotar}[1]%
{\hspace{\SOURCE\unitlength}%
\Hcase{\EDOTAR}{#1}{\ARROWLENGTH}}%

\newcommand{\edotaR}[1]%
{\hspace{\SOURCE\unitlength}%
\hcasE{\EDOTAR}{#1}{\ARROWLENGTH}}%

\newcommand{\emono}%
{\hspace{\SOURCE\unitlength}%
\hcase{\EMONO}{\ARROWLENGTH}}%

\newcommand{\Emono}[1]%
{\hspace{\SOURCE\unitlength}%
\Hcase{\EMONO}{#1}{\ARROWLENGTH}}%

\newcommand{\emonO}[1]%
{\hspace{\SOURCE\unitlength}%
\hcasE{\EMONO}{#1}{\ARROWLENGTH}}%

\newcommand{\eepi}%
{\hspace{\SOURCE\unitlength}%
\hcase{\EEPI}{\ARROWLENGTH}}%

\newcommand{\Eepi}[1]%
{\hspace{\SOURCE\unitlength}%
\Hcase{\EEPI}{#1}{\ARROWLENGTH}}%

\newcommand{\eepI}[1]%
{\hspace{\SOURCE\unitlength}%
\hcasE{\EEPI}{#1}{\ARROWLENGTH}}%

\newcommand{\ebimo}%
{\hspace{\SOURCE\unitlength}%
\hcase{\EBIMO}{\ARROWLENGTH}}%

\newcommand{\Ebimo}[1]%
{\hspace{\SOURCE\unitlength}%
\Hcase{\EBIMO}{#1}{\ARROWLENGTH}}%

\newcommand{\ebimO}[1]%
{\hspace{\SOURCE\unitlength}%
\hcasE{\EBIMO}{#1}{\ARROWLENGTH}}%

\newcommand{\eiso}%
{\hspace{\SOURCE\unitlength}%
\Hcase{\EAR}{\cong}{\ARROWLENGTH}}%

\newcommand{\Eiso}[1]%
{\hspace{\SOURCE\unitlength}%
\Hcase{\EAR}{\cong#1}{\ARROWLENGTH}}%

\newcommand{\eisO}[1]%
{\hspace{\SOURCE\unitlength}%
\hcasE{\EAR}{\cong#1}{\ARROWLENGTH}}%

\newcommand{\ebiar}%
{\hspace{\SOURCE\unitlength}%
\hbicase{\EBIAR}{\ARROWLENGTH}}%

\newcommand{\Ebiar}[2]%
{\hspace{\SOURCE\unitlength}%
\Hbicase{\EBIAR}{#1}{#2}{\ARROWLENGTH}}%

\newcommand{\eeql}%
{\hspace{\SOURCE\unitlength}%
\hbicase{\EEQL}{\ARROWLENGTH}}%

\newcommand{\eadjar}%
{\hspace{\SOURCE\unitlength}%
\hbicase{\EADJAR}{\ARROWLENGTH}}%

\newcommand{\Eadjar}[2]%
{\hspace{\SOURCE\unitlength}%
\Hbicase{\EADJAR}{#1}{#2}{\ARROWLENGTH}}%

\newcommand{\WAR}[1]%
{\begin{picture}(#1,0)%
\put(#1,0){\vector(-1,0){#1}}%
\end{picture}}%

\newcommand{\WDOTAR}[1]%
{\truex{100}\truey{300}%
\NUMBEROFDOTS=#1%
\divide\NUMBEROFDOTS by \value{y}%
\begin{picture}(#1,0)%
\multiput(#1,0)(-\value{y},0){\NUMBEROFDOTS}%
{\circle*{\value{x}}}%
\put(0,0){\vector(-1,0){0}}%
\end{picture}}%

\newcommand{\WMONO}[1]%
{\begin{picture}(#1,0)(-#1,0)%
\put(0,0){\vector(-1,0){#1}}%
\truex{300}\truey{600}%
\put(-\value{x},-\value{x}){\line(0,1){\value{y}}}%
\end{picture}}%

\newcommand{\WEPI}[1]%
{\begin{picture}(#1,0)%
\put(#1,0){\vector(-1,0){#1}}%
\truex{300}\truey{600}\truez{800}%
\put(\value{z},-\value{x}){\line(0,1){\value{y}}}%
\end{picture}}%

\newcommand{\WBIMO}[1]%
{\begin{picture}(#1,0)%
\put(#1,0){\vector(-1,0){#1}}%
\truex{300}\truey{600}\truez{800}%
\put(\value{z},-\value{x}){\line(0,1){\value{y}}}%
\put(#1,-\value{x}){\hspace{-3pt}\line(0,1){\value{y}}}%
\end{picture}}%

\newcommand{\WBIAR}[1]%
{\truex{400}%
\begin{picture}(#1,\value{x})%
\put(#1,0){\vector(-1,0){#1}}%
\put(#1,\value{x}){\vector(-1,0){#1}}%
\end{picture}}%

\newcommand{\WADJAR}[1]%
{\truex{400}%
\begin{picture}(#1,\value{x})%
\put(0,\value{x}){\vector(1,0){#1}}%
\put(#1,0){\vector(-1,0){#1}}%
\end{picture}}%

\newcommand{\war}%
{\hspace{\SOURCE\unitlength}%
\hcase{\WAR}{\ARROWLENGTH}}%

\newcommand{\War}[1]%
{\hspace{\SOURCE\unitlength}%
\Hcase{\WAR}{#1}{\ARROWLENGTH}}%

\newcommand{\waR}[1]%
{\hspace{\SOURCE\unitlength}%
\hcasE{\WAR}{#1}{\ARROWLENGTH}}%

\newcommand{\wdotar}%
{\hspace{\SOURCE\unitlength}%
\hcase{\WDOTAR}{\ARROWLENGTH}}%

\newcommand{\Wdotar}[1]%
{\hspace{\SOURCE\unitlength}%
\Hcase{\WDOTAR}{#1}{\ARROWLENGTH}}%

\newcommand{\wdotaR}[1]%
{\hspace{\SOURCE\unitlength}%
\hcasE{\WDOTAR}{#1}{\ARROWLENGTH}}%

\newcommand{\wmono}%
{\hspace{\SOURCE\unitlength}%
\hcase{\WMONO}{\ARROWLENGTH}}%

\newcommand{\Wmono}[1]%
{\hspace{\SOURCE\unitlength}%
\Hcase{\WMONO}{#1}{\ARROWLENGTH}}%

\newcommand{\wmonO}[1]%
{\hspace{\SOURCE\unitlength}%
\hcasE{\WMONO}{#1}{\ARROWLENGTH}}%

\newcommand{\wepi}%
{\hspace{\SOURCE\unitlength}%
\hcase{\WEPI}{\ARROWLENGTH}}%

\newcommand{\Wepi}[1]%
{\hspace{\SOURCE\unitlength}%
\Hcase{\WEPI}{#1}{\ARROWLENGTH}}%

\newcommand{\wepI}[1]%
{\hspace{\SOURCE\unitlength}%
\hcasE{\WEPI}{#1}{\ARROWLENGTH}}%

\newcommand{\wbimo}%
{\hspace{\SOURCE\unitlength}%
\hcase{\WBIMO}{\ARROWLENGTH}}%

\newcommand{\Wbimo}[1]%
{\hspace{\SOURCE\unitlength}%
\Hcase{\WBIMO}{#1}{\ARROWLENGTH}}%

\newcommand{\wbimO}[1]%
{\hspace{\SOURCE\unitlength}%
\hcasE{\WBIMO}{#1}{\ARROWLENGTH}}%

\newcommand{\wiso}%
{\hspace{\SOURCE\unitlength}%
\Hcase{\WAR}{\cong}{\ARROWLENGTH}}%

\newcommand{\Wiso}[1]%
{\hspace{\SOURCE\unitlength}%
\Hcase{\WAR}{#1}{\ARROWLENGTH}}%

\newcommand{\wisO}[1]%
{\hspace{\SOURCE\unitlength}%
\hcasE{\WAR}{#1}{\ARROWLENGTH}}%

\newcommand{\wbiar}%
{\hspace{\SOURCE\unitlength}%
\hbicase{\WBIAR}{\ARROWLENGTH}}%

\newcommand{\Wbiar}[2]%
{\hspace{\SOURCE\unitlength}%
\Hbicase{\WBIAR}{#1}{#2}{\ARROWLENGTH}}%

\newcommand{\weql}%
{\hspace{\SOURCE\unitlength}%
\hbicase{\EEQL}{\ARROWLENGTH}}%

\newcommand{\wadjar}%
{\hspace{\SOURCE\unitlength}%
\hbicase{\WADJAR}{\ARROWLENGTH}}%

\newcommand{\Wadjar}[2]%
{\hspace{\SOURCE\unitlength}%
\Hbicase{\WADJAR}{#1}{#2}{\ARROWLENGTH}}%

\newcommand{\Vbicase}[4]{\makebox[0pt]%
{\makebox[0pt][r]{\raisebox{0pt}[0pt][0pt]{$#2$\hspace{4pt}}}#1{#4}%
\makebox[0pt][l]{\raisebox{0pt}[0pt][0pt]{\hspace{5pt}$#3$}}}}%

\newcommand{\SAR}[1]%
{\begin{picture}(0,0)%
\put(0,0){\makebox(0,0)%
{\begin{picture}(0,#1)%
\put(0,#1){\vector(0,-1){#1}}%
\end{picture}}}\end{picture}}%

\newcommand{\SDOTAR}[1]%
{\truex{100}\truey{300}%
\NUMBEROFDOTS=#1%
\divide\NUMBEROFDOTS by \value{y}%
\begin{picture}(0,0)%
\put(0,0){\makebox(0,0)%
{\begin{picture}(0,#1)%
\multiput(0,#1)(0,-\value{y}){\NUMBEROFDOTS}%
{\circle*{\value{x}}}%
\put(0,0){\vector(0,-1){0}}%
\end{picture}}}\end{picture}}%

\newcommand{\SMONO}[1]%
{\begin{picture}(0,0)%
\put(0,0){\makebox(0,0)%
{\begin{picture}(0,#1)%
\put(0,#1){\vector(0,-1){#1}}%
\truex{300}\truey{600}%
\put(0,#1){\begin{picture}(0,0)%
\put(-\value{x},-\value{x}){\line(1,0){\value{y}}}\end{picture}}%
\end{picture}}}\end{picture}}%

\newcommand{\SEPI}[1]%
{\begin{picture}(0,0)%
\put(0,0){\makebox(0,0)%
{\begin{picture}(0,#1)%
\put(0,#1){\vector(0,-1){#1}}%
\truex{300}\truey{600}\truez{800}%
\put(-\value{x},\value{z}){\line(1,0){\value{y}}}%
\end{picture}}}\end{picture}}%

\newcommand{\SBIMO}[1]%
{\begin{picture}(0,0)%
\put(0,0){\makebox(0,0)%
{\begin{picture}(0,#1)%
\put(0,#1){\vector(0,-1){#1}}%
\truex{300}\truey{600}\truez{800}%
\put(0,#1){\begin{picture}(0,0)%
\put(-\value{x},-\value{x}){\line(1,0){\value{y}}}\end{picture}}%
\put(-\value{x},\value{z}){\line(1,0){\value{y}}}%
\end{picture}}}\end{picture}}%

\newcommand{\SBIAR}[1]%
{\begin{picture}(0,0)%
\truex{200}%
\put(0,0){\makebox(0,0)%
{\begin{picture}(0,#1)\put(-\value{x},#1){\vector(0,-1){#1}}%
\put(\value{x},#1){\vector(0,-1){#1}}%
\end{picture}}}\end{picture}}%

\newcommand{\SEQL}[1]%
{\begin{picture}(0,0)%
\truex{100}%
\put(0,0){\makebox(0,0)%
{\begin{picture}(0,#1)\put(-\value{x},#1){\line(0,-1){#1}}%
\put(\value{x},#1){\line(0,-1){#1}}%
\end{picture}}}\end{picture}}%

\newcommand{\Sisov}[2]%
{\Vbicase{\SAR}{#1\hspace{-2pt}}{\hspace{-2pt}\cong}{#200}}%

\newcommand{\NAR}[1]%
{\begin{picture}(0,0)%
\put(0,0){\makebox(0,0)%
{\begin{picture}(0,#1)\put(0,0){\vector(0,1){#1}}%
\end{picture}}}\end{picture}}%

\newcommand{\NDOTAR}[1]%
{\truex{100}\truey{300}%
\NUMBEROFDOTS=#1%
\divide\NUMBEROFDOTS by \value{y}%
\begin{picture}(0,0)%
\put(0,0){\makebox(0,0)%
{\begin{picture}(0,#1)%
\multiput(0,0)(0,\value{y}){\NUMBEROFDOTS}%
{\circle*{\value{x}}}%
\put(0,#1){\vector(0,1){0}}%
\end{picture}}}\end{picture}}%

\newcommand{\NMONO}[1]%
{\begin{picture}(0,0)%
\put(0,0){\makebox(0,0)%
{\begin{picture}(0,#1)%
\put(0,0){\vector(0,1){#1}}%
\truex{300}\truey{600}%
\put(-\value{x},\value{x}){\line(1,0){\value{y}}}%
\end{picture}}}%
\end{picture}}%

\newcommand{\NEPI}[1]%
{\begin{picture}(0,0)%
\put(0,0){\makebox(0,0)%
{\begin{picture}(0,#1)%
\put(0,0){\vector(0,1){#1}}%
\truex{300}\truey{600}\truez{800}%
\put(0,#1){\begin{picture}(0,0)%
\put(-\value{x},-\value{z}){\line(1,0){\value{y}}}\end{picture}}%
\end{picture}}}\end{picture}}%

\newcommand{\NBIMO}[1]%
{\begin{picture}(0,0)%
\put(0,0){\makebox(0,0)%
{\begin{picture}(0,#1)%
\put(0,0){\vector(0,1){#1}}%
\truex{300}\truey{600}\truez{800}%
\put(-\value{x},\value{x}){\line(1,0){\value{y}}}%
\put(0,#1){\begin{picture}(0,0)%
\put(-\value{x},-\value{z}){\line(1,0){\value{y}}}\end{picture}}%
\end{picture}}}\end{picture}}%

\newcommand{\NBIAR}[1]%
{\begin{picture}(0,0)%
\truex{200}%
\put(0,0){\makebox(0,0)%
{\begin{picture}(0,#1)\put(-\value{x},0){\vector(0,1){#1}}%
\put(\value{x},0){\vector(0,1){#1}}%
\end{picture}}}\end{picture}}%

\newcommand{\Nisov}[2]%
{\Vbicase{\NAR}{#1\hspace{-2pt}}{\hspace{-2pt}\cong}{#200}}%

\newcommand{\NEDOTAR}%
{\truex{100}\truey{212}%
\NUMBEROFDOTS=5800%
\divide\NUMBEROFDOTS by \value{y}%
\begin{picture}(0,0)%
\multiput(-2900,-2900)(\value{y},\value{y}){\NUMBEROFDOTS}%
{\circle*{\value{x}}}%
\put(2900,2900){\vector(1,1){0}}%
\end{picture}}%

\newcommand{\SWDOTAR}%
{\truex{100}\truey{212}%
\NUMBEROFDOTS=5800%
\divide\NUMBEROFDOTS by \value{y}%
\begin{picture}(0,0)%
\multiput(2900,2900)(-\value{y},-\value{y}){\NUMBEROFDOTS}%
{\circle*{\value{x}}}%
\put(-2900,-2900){\vector(-1,-1){0}}%
\end{picture}}%

\newcommand{\SEDOTAR}%
{\truex{100}\truey{212}%
\NUMBEROFDOTS=5800%
\divide\NUMBEROFDOTS by \value{y}%
\begin{picture}(0,0)%
\multiput(-2900,2900)(\value{y},-\value{y}){\NUMBEROFDOTS}%
{\circle*{\value{x}}}%
\put(2900,-2900){\vector(1,-1){0}}%
\end{picture}}%

\newcommand{\NWDOTAR}%
{\truex{100}\truey{212}%
\NUMBEROFDOTS=5800%
\divide\NUMBEROFDOTS by \value{y}%
\begin{picture}(0,0)%
\multiput(2900,-2900)(-\value{y},\value{y}){\NUMBEROFDOTS}%
{\circle*{\value{x}}}%
\put(-2900,2900){\vector(-1,1){0}}%
\end{picture}}%

\newcommand{\ENEAR}[2]%
{\makebox[0pt]{\begin{picture}(0,0)%
\put(0,-150){\makebox(0,0){\begin{picture}(0,0)%
\put(-6600,-3300){\vector(2,1){13200}}%
\truex{200}\truey{800}\truez{600}%
\put(-\value{x},\value{x}){\makebox(0,\value{z})[r]{${#1}$}}%
\put(\value{x},-\value{y}){\makebox(0,\value{z})[l]{${#2}$}}%
\end{picture}}}\end{picture}}}%

\newcommand{\ESEAR}[2]%
{\makebox[0pt]{\begin{picture}(0,0)%
\put(0,-150){\makebox(0,0){\begin{picture}(0,0)%
\put(-6600,3300){\vector(2,-1){13200}}%
\truex{200}\truey{800}\truez{600}%
\put(\value{x},\value{x}){\makebox(0,\value{z})[l]{${#1}$}}%
\put(-\value{x},-\value{y}){\makebox(0,\value{z})[r]{${#2}$}}%
\end{picture}}}\end{picture}}}%

\newcommand{\WNWAR}[2]%
{\makebox[0pt]{\begin{picture}(0,0)%
\put(0,-150){\makebox(0,0){\begin{picture}(0,0)%
\put(6600,-3300){\vector(-2,1){13200}}%
\truex{200}\truey{800}\truez{600}%
\put(\value{x},\value{x}){\makebox(0,\value{z})[l]{${#1}$}}%
\put(-\value{x},-\value{y}){\makebox(0,\value{z})[r]{${#2}$}}%
\end{picture}}}\end{picture}}}%

\newcommand{\WSWAR}[2]%
{\makebox[0pt]{\begin{picture}(0,0)%
\put(0,-150){\makebox(0,0){\begin{picture}(0,0)%
\put(6600,3300){\vector(-2,-1){13200}}%
\truex{200}\truey{800}\truez{600}%
\put(-\value{x},\value{x}){\makebox(0,\value{z})[r]{${#1}$}}%
\put(\value{x},-\value{y}){\makebox(0,\value{z})[l]{${#2}$}}%
\end{picture}}}\end{picture}}}%

\newcommand{\NNEAR}[2]%
{\raisebox{-1pt}[0pt][0pt]{\begin{picture}(0,0)%
\put(0,0){\makebox(0,0){\begin{picture}(0,0)%
\put(-3300,-6600){\vector(1,2){6600}}%
\truex{100}\truez{600}%
\put(-\value{x},\value{x}){\makebox(0,\value{z})[r]{${#1}$}}%
\put(\value{x},-\value{z}){\makebox(0,\value{z})[l]{${#2}$}}%
\end{picture}}}\end{picture}}}%

\newcommand{\SSWAR}[2]%
{\raisebox{-1pt}[0pt][0pt]{\begin{picture}(0,0)%
\put(0,0){\makebox(0,0){\begin{picture}(0,0)%
\put(3300,6600){\vector(-1,-2){6600}}%
\truex{100}\truez{600}%
\put(-\value{x},\value{x}){\makebox(0,\value{z})[r]{${#1}$}}%
\put(\value{x},-\value{z}){\makebox(0,\value{z})[l]{${#2}$}}%
\end{picture}}}\end{picture}}}%

\newcommand{\SSEAR}[2]%
{\raisebox{-1pt}[0pt][0pt]{\begin{picture}(0,0)%
\put(0,0){\makebox(0,0){\begin{picture}(0,0)%
\put(-3300,6600){\vector(1,-2){6600}}%
\truex{200}\truez{600}%
\put(\value{x},\value{x}){\makebox(0,\value{z})[l]{${#1}$}}%
\put(-\value{x},-\value{z}){\makebox(0,\value{z})[r]{${#2}$}}%
\end{picture}}}\end{picture}}}%

\newcommand{\NNWAR}[2]%
{\raisebox{-1pt}[0pt][0pt]{\begin{picture}(0,0)%
\put(0,0){\makebox(0,0){\begin{picture}(0,0)%
\put(3300,-6600){\vector(-1,2){6600}}%
\truex{200}\truez{600}%
\put(\value{x},\value{x}){\makebox(0,\value{z})[l]{${#1}$}}%
\put(-\value{x},-\value{z}){\makebox(0,\value{z})[r]{${#2}$}}%
\end{picture}}}\end{picture}}}%

\newcommand{\Necurve}[2]%
{\begin{picture}(0,0)%
\truex{1300}\truey{2000}\truez{200}%
\put(0,\value{x}){\oval(#200,\value{y})[t]}%
\put(0,\value{x}){\makebox(0,0){\begin{picture}(#200,0)%
\put(#200,0){\vector(0,-1){\value{z}}}%
\put(0,0){\line(0,-1){\value{z}}}\end{picture}}}%
\truex{2500}%
\put(0,\value{x}){\makebox(0,0)[b]{${#1}$}}%
\end{picture}}%

\newcommand{\Nwcurve}[2]%
{\begin{picture}(0,0)%
\truex{1300}\truey{2000}\truez{200}%
\put(0,\value{x}){\oval(#200,\value{y})[t]}%
\put(0,\value{x}){\makebox(0,0){\begin{picture}(#200,0)%
\put(#200,0){\line(0,-1){\value{z}}}%
\put(0,0){\vector(0,-1){\value{z}}}\end{picture}}}%
\truex{2500}%
\put(0,\value{x}){\makebox(0,0)[b]{${#1}$}}%
\end{picture}}%

\newcommand{\Securve}[2]%
{\begin{picture}(0,0)%
\truex{1300}\truey{2000}\truez{200}%
\put(0,-\value{x}){\oval(#200,\value{y})[b]}%
\put(0,-\value{x}){\makebox(0,0){\begin{picture}(#200,0)%
\put(#200,0){\vector(0,1){\value{z}}}%
\put(0,0){\line(0,1){\value{z}}}\end{picture}}}%
\truex{2500}%
\put(0,-\value{x}){\makebox(0,0)[t]{${#1}$}}%
\end{picture}}%

\newcommand{\Swcurve}[2]%
{\begin{picture}(0,0)%
\truex{1300}\truey{2000}\truez{200}%
\put(0,-\value{x}){\oval(#200,\value{y})[b]}%
\put(0,-\value{x}){\makebox(0,0){\begin{picture}(#200,0)%
\put(#200,0){\line(0,1){\value{z}}}%
\put(0,0){\vector(0,1){\value{z}}}\end{picture}}}%
\truex{2500}%
\put(0,-\value{x}){\makebox(0,0)[t]{${#1}$}}%
\end{picture}}%

\newcommand{\Escurve}[2]%
{\begin{picture}(0,0)%
\truex{1400}\truey{2000}\truez{200}%
\put(\value{x},0){\oval(\value{y},#200)[r]}%
\put(\value{x},0){\makebox(0,0){\begin{picture}(0,#200)%
\put(0,0){\vector(-1,0){\value{z}}}%
\put(0,#200){\line(-1,0){\value{z}}}\end{picture}}}%
\truex{2500}%
\put(\value{x},0){\makebox(0,0)[l]{${#1}$}}%
\end{picture}}%

\newcommand{\Encurve}[2]%
{\begin{picture}(0,0)%
\truex{1400}\truey{2000}\truez{200}%
\put(\value{x},0){\oval(\value{y},#200)[r]}%
\put(\value{x},0){\makebox(0,0){\begin{picture}(0,#200)%
\put(0,0){\line(-1,0){\value{z}}}%
\put(0,#200){\vector(-1,0){\value{z}}}\end{picture}}}%
\truex{2500}%
\put(\value{x},0){\makebox(0,0)[l]{${#1}$}}%
\end{picture}}%

\newcommand{\Wscurve}[2]%
{\begin{picture}(0,0)%
\truex{1300}\truey{2000}\truez{200}%
\put(-\value{x},0){\oval(\value{y},#200)[l]}%
\put(-\value{x},0){\makebox(0,0){\begin{picture}(0,#200)%
\put(0,0){\vector(1,0){\value{z}}}%
\put(0,#200){\line(1,0){\value{z}}}\end{picture}}}%
\truex{2400}%
\put(-\value{x},0){\makebox(0,0)[r]{${#1}$}}%
\end{picture}}%

\newcommand{\Wncurve}[2]%
{\begin{picture}(0,0)%
\truex{1300}\truey{2000}\truez{200}%
\put(-\value{x},0){\oval(\value{y},#200)[l]}%
\put(-\value{x},0){\makebox(0,0){\begin{picture}(0,#200)%
\put(0,0){\line(1,0){\value{z}}}%
\put(0,#200){\vector(1,0){\value{z}}}\end{picture}}}%
\truex{2400}%
\put(-\value{x},0){\makebox(0,0)[r]{${#1}$}}%
\end{picture}}%

\newcount\SCALE%

\newcount\NUMBER%

\newcount\LINE%

\newcount\COLUMN%

\newcount\WIDTH%

\newcount\SOURCE%

\newcount\ARROW%

\newcount\TARGET%

\newcount\ARROWLENGTH%

\newcount\NUMBEROFDOTS%

\newcounter{x}%

\newcounter{y}%

\newcounter{z}%

\newcounter{horizontal}%

\newcounter{vertical}%

\newskip\itemlength%

\newskip\firstitem%

\newskip\seconditem%

\newcommand{\printarrow}{}%

\newcommand{\truex}[1]{%
\NUMBER=#1%
\multiply\NUMBER by 100%
\divide\NUMBER by \SCALE%
\setcounter{x}{\NUMBER}}%

\newcommand{\truey}[1]{%
\NUMBER=#1%
\multiply\NUMBER by 100%
\divide\NUMBER by \SCALE%
\setcounter{y}{\NUMBER}}%

\newcommand{\truez}[1]{%
\NUMBER=#1%
\multiply\NUMBER by 100%
\divide\NUMBER by \SCALE%
\setcounter{z}{\NUMBER}}%

\newcommand{\changecounters}[1]{%
\SOURCE=\ARROW%
\ARROW=\TARGET%
\settowidth{\itemlength}{#1}%
\ifdim \itemlength > 2800\unitlength%
\addtolength{\itemlength}{-2800\unitlength}%
\TARGET=\itemlength%
\divide\TARGET by 1310%
\multiply\TARGET by 100%
\divide\TARGET by \SCALE%
\else%
\TARGET=0%
\fi%
\ARROWLENGTH=5000%
\advance\ARROWLENGTH by -\SOURCE%
\advance\ARROWLENGTH by -\TARGET%
\advance\SOURCE by -\TARGET}%

\newcommand{\initialize}[1]{%
\LINE=0%
\COLUMN=0%
\WIDTH=0%
\ARROW=0%
\TARGET=0%
\changecounters{#1}%
\renewcommand{\printarrow}{#1}%
\begin{center}%
\vspace{10pt}%
\begin{picture}(0,0)}%

\newcommand{\n}[1]{%
\changecounters{\mbox{$#1$}}%
\put(\COLUMN,\LINE){\makebox(0,0){\printarrow}}%
\thinlines%
\renewcommand{\printarrow}{\mbox{$#1$}}%
\advance\COLUMN by 4000}%

\newcommand{\nn}[1]{%
\put(\COLUMN,\LINE){\makebox(0,0){\printarrow}}%
\thinlines%
\ifnum \WIDTH < \COLUMN%
\WIDTH=\COLUMN%
\else%
\fi%
\advance\LINE by -4000%
\COLUMN=0%
\ARROW=0%
\TARGET=0%
\changecounters{\mbox{$#1$}}%
\renewcommand{\printarrow}{\mbox{$#1$}}}%

\newcommand{\conclude}{%
\put(\COLUMN,\LINE){\makebox(0,0){\printarrow}}%
\thinlines%
\ifnum \WIDTH < \COLUMN%
\WIDTH=\COLUMN%
\else%
\fi%
\setcounter{horizontal}{\WIDTH}%
\setcounter{vertical}{-\LINE}%
\end{picture}}%

\newcommand{\diag}{%
\conclude%
\raisebox{0pt}[0pt][\value{vertical}\unitlength]{}%
\hspace*{\value{horizontal}\unitlength}%
\vspace{10pt}%
\end{center}%
\setlength{\unitlength}{1pt}}%

\newcommand{\diagv}[3]{%
\conclude%
\NUMBER=#1%
\rule{0pt}{\NUMBER pt}%
\hspace*{-#2pt}%
\raisebox{0pt}[0pt][\value{vertical}\unitlength]{}%
\hspace*{\value{horizontal}\unitlength}
\NUMBER=#3%
\advance\NUMBER by 10%
\vspace*{\NUMBER pt}%
\end{center}%
\setlength{\unitlength}{1pt}}%

\newcommand{\N}[1]%
{\raisebox{0pt}[7pt][0pt]{$#1$}}%

\newcommand{\crosslength}[2]{%
\settowidth{\firstitem}{#1}%
\settowidth{\seconditem}{#2}%
\ifdim\firstitem < \seconditem%
\itemlength=\seconditem%
\else%
\itemlength=\firstitem%
\fi%
\divide\itemlength by 2%
\hspace{\itemlength}}%

\begin{document}

\title {Quadratic forms and singularities of genus one or two}
\author {Georges Dloussky\footnote{Georges Dloussky, dloussky@cmi.univ-mrs.fr, , version December 2007. 
}
}
\date{ }

\maketitle
\abstract{We study singularities obtained by the contraction of the maximal divisor in compact (non-k\H ahlerian) surfaces which contain global spherical shells. These singularities are of genus 1 or 2, may be $\bb Q$-Gorenstein, numerically Gorenstein or Gorenstein. A family of polynomials depending on the configuration of the curves computes the discriminants of the quadratic forms of these singularities. We introduce a multiplicative branch topological invariant which determines the twisting of a non-vanishing holomorphic 1-form on the complement of the singular point.} 
\tableofcontents

\section{Introduction}
We are interested in a large class of singularities which generalize cusps, obtained by the contraction of all the rational curves in compact surfaces $S$ wich contain global spherical shells. Particular cases are Inoue-Hirzebruch surfaces with  two ``dual'' cycles of rational curves. The duality can be explained by the construction of these surfaces by sequences of blowing-ups \cite{D2}. Several authors have studied cusps \cite{H,KA, N3, N4, LW}. In general, the maximal divisor is composed of a cycle with branches. These (non-k\"ahlerian) surfaces contain exactly $n=b_2(S)$
 rational
curves. The intersection matrices $M(S)$ have been completely
classified \cite{N2, D1};  they are negative definite in all cases except when the maximal divisor is a cycle $D$ of $n$ rational curves
such that
$D^{2}=0$. In this article we study normal singularities obtained  by contraction of the exceptional divisor and the link between the 
intersection matrix and global topological or analytical properties of the surface $S$. They are elliptic or of genus two in which case they are Gorenstein. Using the existence of global section on $S$ of $-mK_S\otimes L$ for a suitable integer $m\ge 1$ and a flat line bundle $L\in H^1(S,\bb C^\star)$, we show that these singularities are $\bb Q$-Gorenstein (resp. numerically Gorenstein) if and only if the global property $H^0(S,-mK_S)\neq 0$ (resp. $H^0(S, -K_S\otimes L)\neq 0$) holds. The main part of this article is devoted to the study of the discriminant of the quadratic form associated to the singularity. In \cite{D1} the quadratic form has been decomposed into a sum of squares. The intersection matrix is completely determined by the sequence $\s$ of (opposite) self-intersections of the rational curves  when taken in the canonical order, i.e. the order in which the curves are obtained in a repeted sequence of blowing-ups. Let $(Y,y)=(Y_\s,y)$  be the associated singularity obtained by the contraction of the rational curves. We introduce  a family of polynomials $P_\s=P_{A(\s)}$ which have integer values on integers, depending on the configuration of the dual graph of the singularity, such that the discriminant is the square of this polynomial. When we fix the sequence $\s$ we obtain an integer $\D_\s$ which is a multiplicative topological invariant  i.e. satisfies $\D_{\s\s'}=\D_\s\D_{\s'}$. We show that $\D_\s$ is equal to the product of the determinants of the branches. We apply this result to determine the twisting integer of holomorphic 1-forms in a neighbourhood of the singularity. We develop here rather the algebraic point of view, however these singularities have deep relations with properties of compact complex surfaces $S$  containing global spherical shells,  the classification of singular contracting germs of mappings and dynamical systems: for instance,  the integer  $\D_\s$ is equal to the integer $k=k(S)$ wich appears in the normal form of contracting germs $F(z_1,z_2)=(\l z_1z_2^s+P(z_2),z_2^k)$ which define $S$ \cite{D2germe, DO, DO2, F}.\\
I thank Karl Oeljeklaus for fruitful discussions on that subject.
\section{Preliminaries}
\subsection{Basic results on singularities}
 Let $D_{0},\dots ,D_{n-1}$ be compact curves on a (not necessarilly 
compact) complex surface  $X$, and $D=D_0+\cdots+D_{n-1}$ the associated reduced divisor. We assume  that  $D$ is exceptional i.e. the intersection matrix M of $D$ is 
negative definite. We denote by $\cal O_X$ the structural sheaf of $X$, $K_X=\det T^\star X$ the canonical bundle and by $\O^2_X$ its sheaf of sections. It is well known by Grauert's theorem that there exists a proper 
mapping \hbox {$\Pi
:X\rightarrow Y$} such that each connected component of $|D|=\cup_i D_i$ is contracted onto a 
point 
$y$
which is a normal singularity of $Y$. For $|D|$ connected, denote  by 
$$r:H^{0}(X, \Omega
_{X}^{2})\rightarrow   H^{0}(Y\backslash \lbrace y\rbrace, 
\Omega 
_{Y\backslash\lbrace
y\rbrace}^{2}) $$
 the canonical morphism induced by $\Pi $. We define the 
{\bf geometric genus}
of the singularity $(Y, y)$ by 
$$p_{g} = p_{g}(Y,y) = h^{0}(Y, R^{1}\Pi _{*}{\cal 
O}_{X}).$$ 
When $Y$ is Stein, we have  $p_{g}  = \dim  H^{0}(Y\backslash\lbrace y\rbrace,
\Omega _{Y\backslash\lbrace y\rbrace}^{2}) / rH^{0}(X, \Omega 
_{X}^{2})$.\\

 A singularity $( Y, y)$ is called {\bf rational} (resp. {\bf elliptic})
 if $p_{g}( Y, y)=0$  (resp.
$p_{g}(Y, y)=1$).
Therefore a singularity is rational if for every 2-form $\o$ on $Y\backslash\lbrace y\rbrace$, the 2-form $\Pi^\star \o$ extends to a 2-form on $X$.\\
\begin{Prop} \label{SS} Let $\Pi:X\to Y$ be the proper morphism obtained by the contraction of an exceptional divisor:\\
1) $p_{g}=h^{0}(Y, R^{1}\Pi _{*}{\cal 
O}_{X})$ is independant of the choice of 
the desingularization.
\begin{itemize}
\item If $X$ is compact then $p_g=\chi(O_Y)-\chi(O_X)$
 \item If $X$ is spc and $Y$ is Stein then $p_g=h^1(X,O_X)$
\end{itemize}
2) The following sequence $$ 0\rightarrow  H^{1}( Y, {\cal 
O}_{Y})
\rightarrow  H^{1}( X, {\cal O}_{X})
\rightarrow H^{0}Y, R^{1}\Pi _\star{\cal O}_{X})\rightarrow  
H^{2}(Y,
{\cal O}_{Y})\rightarrow  H^{2}( X, {\cal O}_{X})$$
is exact.
\end{Prop}

We give now a criterion of rationality \cite{PIN}, p. 152:
\begin{Prop} \label{SRat}
Let $\Pi :X\rightarrow Y$ be the minimal resolution of the 
singularity 
$(Y, y)$ and denote
by $D_{i}$ the irreducible components of the exceptional divisor $D$. 
If 
$(Y, y)$ is rational,
then:\\
i) the curves $D_{i}$ are regular and rational\\
ii) for $i\neq j$, $D_i\cap D_j=\emptyset$ or $D_{i}$ meets $D_{j}$ tranversally. If $D_{i}$, 
$D_{j}$, $D_{k}$ are
distinct irreducible components, $D_{i}\cap D_{j}\cap D_{k}$ is 
empty\\
iii)  D contains no cycle.
\end{Prop}
\begin{Def} A singularity $(Y,y)$ is called {\bf Gorenstein}  if the dualizing sheaf $\omega_{Y}$
is trivial, i.e. there exists a small 
neighbourhood $ U$ of $y$ and a non-vanishing holomorphic 2-form on
$U\setminus\lbrace y\rbrace$.
\end{Def}
 Since there is only a finite 
number of linearly
independant 2-forms in the complement of the exceptional divisor $D$ 
modulo 
$H^{0}(X,\Omega^{2})$, a 2-form extends meromorphically
across $D$. Therefore we have (see \cite{SA})

\begin{Lem}\label{Gor} Let $Y$ be a Gorenstein normal surface 
and 
$\Pi : X\rightarrow Y$ be the
minimal desingularization. Then there is a unique effective divisor 
$\Delta$ 
on $X$  supported
in $D=\Pi^{-1}(Sing(Y))$ such that  $$\omega_{X}\simeq
\Pi^\star\omega_{Y}\mathop{\otimes}_{{\cal O}_{X}}{\cal 
O}_{X}(-\Delta)$$ 
\end{Lem}

\subsection{Lattices}
 Here are recalled some well known facts
about lattices (see \cite{SE}).
 We call lattice, denoted by $\bigl
(L, < \ .\ ,\ .\ >
\bigr )$, a free {\bf Z}-module $L$, endowed with an integral non
degenerate symetric bilinear
form
$$\begin{array}{lccc}
< \ .\ ,\ .\ > \ :&L\times L&\longrightarrow &{\bf Z}\\
&(x, y)&\longmapsto &< \ x\ ,\ y\ > .
\end{array}$$
If $B=\{e_{1},\ldots , e_{n}\}$ is a basis of $L$, the determinant of the
matrix
$$\bigl (< \ e_{i}\ ,\ e_{j}\ > \bigr )_{1\leq i, j\leq n}, $$ is
independent of the choice of the
basis; this integer, denoted by $d(L)$ is called the discriminant of the
lattice. A lattice is unimodular
if $d(L)=\pm 1$. Let  $L^{\vee}:= Hom_{\bb Z}(L, {\bb Z})$ be the dual of
$L$.
The mapping
$$\begin{array}{cccc}
\phi : &L&\longrightarrow &L^{\vee}\\
&x&\longmapsto &< \ .\ ,\ x\ > 
\end{array}$$
identifies $L$
with a sublattice of $L^{\vee}$ of same rank, since
$d(L)\not= 0$.
Moreover, if
$L_{\bb Q}:=L\otimes_{\bb Z}{\bb Q}$,
 it is possible to identify $L^{\vee}$
with the sub-{$\bb Z$}-module 
$$ \bigl\lbrace x\in  L_{\bb Q}\ \mid \ \forall y\in L,\ < \ x\ ,\ y\ > \in {\bb Z}
\bigr \rbrace$$
of  $L_{\bb Q}$. So, we may write $L\subset L^{\vee}\subset L_{\bb Q}$,
where $L$ and
$L^{\vee}$ have same rank. \hfill\break
\begin{Lem}\label{detcarre}
 1) The index of $L$ in $L^{\vee}$ is $\mid d(L)\mid$.\hfill\break
2) If $M$ is a submodule of $L$ of the same rank, then the index of $M$ in
$L$ satifies
$$[L\ :\ M]^{2}\ =\ d(M)\ d(L)^{-1}.$$
In particular $d(M)$ and $d(L)$ have same sign.
\end{Lem}

\subsection{Surfaces with global spherical shells}
We recall  some properties of these surfaces which have been
first introduced by Ma. Kato \cite{K} and we refer to \cite{D1} for details.
\begin{Def} Let $S$ be a compact complex surface. We say that $S$ contains a global spherical shell, if there is a biholomorphic map $\f:U\to S$ from a neighbourhood $U\subset \bb C^2\setminus\{0\}$ of the sphere $S^3$ into $S$ such that $S\setminus \f(S^3)$ is connected.
\end{Def}
Hopf surfaces are the simplest examples of surfaces with GSS (see \cite{D1}), however they contain no rational curves and elliptic curves have self-intersection equal to $0$, hence no singularity can be obtained.\\

Let $S$ be a minimal surface containing a GSS with $n=b_2(S)$. It is known that $S$ contains $n$ rational curves and to each curve it is possible to associate a contracting germ of mapping $F=\Pi\s=\Pi_0\cdots\Pi_{n-1}\s:(\bb C^2,0) \to (\bb C^2,0)$ where $\Pi=\Pi_0\cdots\Pi_{n-1}:B^\Pi\to B$ is a sequence of $n$ blowing-ups. If we want to obtain a minimal surface, the sequence of blowing-ups has to be done in the following way:
\begin{itemize}
\item $\Pi_0 $ blows up the origin of the two dimensional unit ball $B$, 
\item $\Pi_1$ blows up a point $O_0\in C_0=\Pi_0^{-1}(0)$,\ldots 
\item $\Pi_{i+1}$ blows up a point $O_{i}\in C_{i}=\Pi_{i}^{-1}(O_{i-1})$, for $i=0,\ldots,n-2$, and 
\item $\s:\bar B\to B^\Pi$ sends isomorphically a neighbourhood of $\bar B$ onto a small ball in $B^\Pi$ in such a way that $\s(0)\in C_{n-1}$. 
\end{itemize}

 It is easy to see that
the homological groups
satisfy $$H_{1}(S,{\bb Z})\simeq {\bb Z},\quad H_{2}(S,{\bb Z})\simeq {\bb
Z}^n$$ In particular,
$b_{2}(S)=n$.\\
The universal covering space $(\tilde S,  \omega , S)$ of S contains only
rational curves
$(C_{i})_{i\in {\bb Z}}$ with a canonical order relation, ``the order of
creation'' (\cite{D1}, p 29).
Following \cite{D1}, we can associate to S the following invariants:
\begin{itemize}
\item  The family
of opposite
self-intersection of curves of the universal covering space of S, denoted by
$$a(S) := (a_{i})_{i\in {\bb Z}} = (- C_{i}^{2})_{i\in {\bb Z}}$$
 this
family is periodic of period
$n$,
\item 
$$\displaystyle\sigma _{n}(S) := \sum_{i=j}^{j+n-1}a_{i} =
-\sum_{i=0}^{n-1}D_{i}^{2}+
2\ \sharp \lbrace rational\, curves\, with\, nodes\rbrace$$  where j is any
index, and the $D_{i}$ are
the rational curves of  S. It can be seen that $2n\leq \sigma _{n}(S)\leq 3n$ (\cite{D1}, p 43).
\item  The intersection matrix of the $n$ rational curves of $S$,
$$M(S):=(D_{i}. D_{j}).$$
{\bf Important Remark}: The essential fact useful to understand the dual graph or the intersection matrix is that
\begin{itemize}
\item if $a_i=-D_i^2=2$ then $D_i$ meets $D_{i+1}$,
\item if $a_i=-D_i^2=3$ then $D_i$ meets $D_{i+2}$,\ldots, 
\item if $a_i=-D_i^2=k+2$ then $D_i$ meets $D_{i+k+1}$,
\end{itemize}
 the indices being in $\bb Z/n\bb Z$, in particular $D_i$ may meet itself: we obtain a rational curve with double point.
\item $n$ classes of contracting holomorphic germs of mappings $F=\Pi \sigma :({\bb
C}^2, 0)\rightarrow
({\bb C}^2, 0)$ (\cite{D1}, p 32).
\end{itemize}
\begin{Prop}\label{matriceintdetcarre}
 Let $S$ be a surface containing a GSS with $b_{2}(S)=n$, $D_{0}, \ldots ,
D_{n-1}$ the n rational curves and $M(S)$ the intersection matrix. \hfill\break
1) If $\sigma_{n}(S)=2n$, then $det M(S)=0$.\hfill\break
2) If $\sigma_{n}(S)>2n$, then $\displaystyle \sum_{0\leq i\leq n-1}{\bf
Z}D_{i}$ is a complete
sublattice of $H_{2}(S, {\bf Z})$ and its index satisfies
$$\bigl [H_{2}(S, {\bf Z})\ :\ \sum_{0\leq i\leq n-1}{\bf Z}D_{i}\bigr
]^{2}\ =\ det M(S).$$
In particular, $det M(S)$ is the square of an integer $\geq 1$.
\end{Prop}
{\bf Proof: } If $\sigma_{n}(S)=2n$, $S$ is an Inoue surface; if
$\sigma_{n}(S)>2n$, $det
M(S)\not= 0$ so the sublattice is complete and the result is a mere consequence
of lemma \ref{detcarre}. \hfill $\Box$\\

In order to give a precise description of the intersection matrix we need the following definitions:
\begin{Def} Let $1\leq
p\leq n$.
A p-uple $\sigma =(a_{i},\dots , a_{i+p-1})$  of $a(S)$ is called 
\begin{itemize}
\item a
{\bf singular p-sequence} of $a(S)$ if
$$\displaystyle\sigma =(\underbrace{p+2, 2,\dots , 2}_{p}).$$
 It will be
denoted by $s_{p}$.
\item   a  {\bf regular
p-sequence} of $a(S)$ if
$$\displaystyle\sigma =(\underbrace{2, 2,\dots , 2}_{p})$$
 and $\sigma$ has
no common element
with a singular sequence. Such a p-uple will be denoted by $r_{p}$.
\end{itemize}
\end{Def}
  For example $s_1=(3), s_2=(4, 2), s_3=(5,
2, 2) ,\dots $ are singular sequences, $r_3=(2, 2, 2)$ is a regular sequence. It is easy to see that
 if we want to have, for example, a curve with self-intersection -4,
necessarily, the 
curve which follows in the sequence of repeted blowing-ups must have self-intersection -2, so it is easy to see (\cite{D1}, p39), that 
$a(S)$ admits a unique
partition by $N$ singular sequences and $\r\leq N$ regular sequences of maximal length.
More precisely, since
$a(S)$ is periodic it is possible to find a n-uple $\sigma$ such that 
$$\sigma
=\sigma_{p_{0}}\cdots \sigma_{p_{\r+N-1}},$$ 
where $\sigma_{p_{i}}$ is a
regular or a singular
$p_{i}$-sequence with $$\sum_{i=0}^{N+\r-1}p_{i} = n$$ and if
$\sigma_{p_{i}}$ is  regular it is between (mod. $N+\r$) two singular sequences.\hfill\break
{\bf Notation}: We shall write $$a(S) = (\overline{\s})=(\overline {\sigma_{p_{0}}\cdots
\sigma_{p_{N+\r-1}}}).$$
The sequence $\s$ is overlined to indicate  that the sequence $\s$ is infinitely repeted to obtain the sequence $a(S)=(a_i)_{i\in\bb Z}$. The sequence $a(S)$ may be defined by another period. For example
$$a(S)=(\overline {\sigma_{p_{1}}\cdots
\sigma_{p_{N+\r-1}}\s_{p_0}}).$$
If $\sigma _{n}(S)=2n$, $a(S)=({\overline r_{n}})$; if $\sigma
_{n}(S)=3n$,
$a(S)$ is only composed  of singular sequences and S is called a
Inoue-Hirzebruch surface.
Moreover if $a(S)$ is composed by the repetition of an even (resp. odd)
number of sequences
$\sigma_{p_{i}}$,  we shall say that S is an even (resp. odd)
Inoue-Hirzebruch surface. An
even (resp. odd) Inoue-Hirzebruch surface has exactly 2  cycles (resp. 1 cycle) of rational curves. Another
used terminology is respectively hyperbolic Inoue surface and half Inoue
surface.\\
 We recall that for any  $VII_{0}$-class surface without non
constant meromorphic
functions, the numerical characters of S are [10, I p755, II p683]
$$h^{0, 1}=1, \,  h^{1, 0}=h^{2, 0}=h^{0, 2}= 0,\, -c_{1}^{2}=c_{2}=b_{2}(S), \,
b_{2}^{+}=0, \,  b_{2}^{-}=b_{2}(S)$$

We shall need in the sequel the explicit description of the dual graph which is composed of a cycle and  branches. A branch $A_s$ determines and is determined by a piece $\G_s$ of the cycle $\G$.

\begin{Th}[\cite{D1} thm 2.39]\label{Arbresetcycle} Let $S$ be a minimal surface containing a GSS, $n=b_2(S)$, $D_0,\ldots,D_{n-1}$ its $n$ rational curves and $D=D_0+\cdots+D_{n-1}$.\\
1) If $\s_n(S)=2n$, then $D$ is a cycle and $D_i^2=-2$ for $i=0,\ldots,n-1$.\\
2) If $2n<\s_n(S)<3n$, then there are $\r(S)\ge 1$ branches and 
$$D=\sum_{s=1}^{\r(S)} (A_s+\G_s)$$
where
\begin{description}
\item{i)} $A_s$ is a branch for $s=1,\ldots,\r(S)$,
\item{ii)} $\G=\sum_{s=1}^{\r(S)}\G_s$ is a cycle,
\item{iii)} $A_s$ and $\G_s$ are defined in the following way: For each sequence of integers
$$(a_{t+1},\ldots,a_{t+l+k_1+\cdots+k_p+2})=(r_l s_{k_1}\cdots s_{k_p} 2 a_{t+l+k_1+\cdots+k_p+2})$$
contained in $a(S)=(\overline{\s_0\cdots \s_{N+\r-1}})$, where 
\begin{itemize}
\item $l\geq 1$ and $r_l$ is a regular $l$-sequence,
\item $p\ge 1$,   $i=1,\ldots,p$, $k_i\ge 1$ and $s_{k_i}$, is a singular $k_i$-sequence,
\end{itemize}
 we have the following decomposition into branch $A_s$ and corresponding piece of cycle $\G_s$ (where $p=p_s$ to simplify notations):

$$\hspace{-7mm}\left\{\begin{array}{lcl}
Selfint(A_s)&=&(\underbrace{2,\ldots,2}_{k_1-1},\  k_2+2,\  \underbrace{ 2,\ldots,2}_{k_3-1}, \ldots,\  k_{p-1}+2,\   \underbrace{2,\ldots,2}_{k_p-1},\  2) \\
&&\\
&&\hfill  If\  p\equiv 1 ({\rm mod}\ 2)\\
&&\\
Selfint(\G_s)&=&( \underbrace{2,\ldots,2}_{l-1},\ k_1+2,\  \underbrace{2,\ldots,2}_{k_2-1},\ldots,\ k_{p-2}+2,\  \underbrace{2,\ldots,2}_{k_{p-1}-1},\ k_p+2)
\end{array} \right.$$

$$\hspace{-7mm}\left\{\begin{array}{lcl}
Selfint(A_s)&=&(\underbrace{2,\ldots,2}_{k_1-1},\  k_2+2,\  \underbrace{ 2,\ldots,2}_{k_3-1}, \ldots,\  k_{p-2}+2,\   \underbrace{2,\ldots,2}_{k_{p-1}-1},\  k_p+2) \\
&&\\
&&\hfill  If\  p\equiv 0 ({\rm mod}\ 2)\\
&&\\
Selfint(\G_s)&=&( \underbrace{2,\ldots,2}_{l-1},\ k_1+2,\  \underbrace{2,\ldots,2}_{k_2-1},\ldots,\ k_{p-1}+2,\  \underbrace{2,\ldots,2}_{k_{p}-1},\ 2)
\end{array} \right.$$

\item{iv)} The top of the branch $A_s$ is its first vertex (or curve); the root of $A_s$  is the first vertex (or curve) of $\G_t$ where $t=s+1$ (mod $\r(S)$).

\end{description}

\noindent 3) If $\s_n(S)=3n$, $D$ has no branch and
\begin{description}
\item{i)} If $a(S)=(\overline{s_{k_1}\cdots s_{k_{2p}}})$ then
$$D=\G+\G'$$
where $\G$ and $\G'$ are two cycles
$$\hspace{-7mm}\left\{\begin{array}{lcl}
Selfint(\G)&=&(k_1+2, \underbrace{2,\ldots,2}_{k_2-1},\ k_3+2, \underbrace{2,\ldots,2}_{k_4-1},\ldots,\ k_{2p-1}+2, \underbrace{2,\ldots,2}_{k_{2p}-1})\\
&&\\
Selfint(\G')&=&(\underbrace{2,\ldots,2}_{k_1-1},\ k_2+2, \underbrace{2,\ldots,2}_{k_3-1}, k_4+2, \ldots, \underbrace{2,\ldots,2}_{k_{2p-1}-1}, k_{2p}+2)
\end{array} \right.$$
\item{ii)} If $a(S)=(\overline{s_{k_1}\cdots s_{k_{2p+1}}})$ then $D$ is contains only one cycle and
$$\begin{array}{lcl}
Selfint(D) &=&
(k_1+2, \underbrace{2,\ldots,2}_{k_2-1},\ k_3+2, \underbrace{2,\ldots,2}_{k_4-1},\ldots,\ k_{2p+1}+2, \\
&&\\
&&\underbrace{2,\ldots,2}_{k_{1}-1}, k_2+2,  \underbrace{2,\ldots,2}_{k_{3}-1},\ldots, k_{2p}+2, \underbrace{2,\ldots,2}_{k_{2p+1}-1})
\end{array}$$

\end{description}
\end{Th}

\subsection{Intersection matrix of the exceptional divisor}
Let $\s=\s_0\cdots\s_{N+\r-1}$ where $\s_i=r_{p_i}=(2,2,\ldots,2)$ is a regular sequence of length $p_i$ or $\s_i=s_{p_i}=(p_i+2,2,\ldots,2)$ is a singular sequence of length $p_i$, $i=0,\ldots,N+\r-1$. We suppose that 
\begin{itemize}
\item there are $N$ singular sequences and $\r\leq N$ regular sequences
\item if $\s_i$ is regular, then $\s_{i-1}$ and $\s_{i+1}$ are singular, indices being in $\bb Z/(N+\r)\bb Z$.
\end{itemize}
Let $n=\sum_{i=0}^{N+\r-1}p_i$ be the number of integers in the sequence $\s$.
\begin{Ex} For $0\le N\le 3$ we have the following possible sequences:\begin{itemize}
\item If $N=0$, $\s=r_n$,
\item If $N=1$, $\s=s_n$ or $\s=s_{p}r_m$,  $p+m=n$,
\item If $N=2$, $\s=s_{p_0}s_{p_1}$, $\s=s_{p_0}s_{p_1}r_{m_0}$, $\s=s_{p_0}r_{m_0}s_{p_1}$, $\s=s_{p_0}r_{m_0}s_{p_1}r_{m_1}$,
\item  If $N=3$, $\s=s_{p_0}s_{p_1}s_{p_2}$\\
 $\s=s_{p_0}r_{m_0}s_{p_1}s_{p_2}$, $\s=s_{p_0}s_{p_1}r_{m_0}s_{p_2}$, $\s=s_{p_0}s_{p_1}s_{p_2}r_{m_0}$, \\
 $\s=s_{p_0}r_{m_0}s_{p_1}r_{m_1}s_{p_2}$, $\s=s_{p_0}s_{p_1}r_{m_0}s_{p_2}r_{m_1}$, $\s=s_{p_0}r_{m_0}s_{p_1}s_{p_2}r_{m_1}$,\\
  $\s=s_{p_0}r_{m_0}s_{p_1}r_{m_1}s_{p_2}r_{m_2}$.
  \end{itemize}
  \end{Ex}
To a sequence $\s$ we associate a symmetric matrix of type $(n,n)$, $M(\s)=(m_{ij})$  ``written on a torus'', with indices in $\bb Z/n\bb Z$ defined in the following way: if $\s=\s_0\cdots\s_{N+\r-1}=(a_0,\ldots,a_{n-1})$
\begin{description}
\item{i)} $m_{ii}=\left\{\begin{array}{lcl}a_i&if&a_i\neq n+1\\n-1&if&a_i=n+1\end{array}\right.$
\item{ii)} For $0\leq i<j\leq n-1$,
$$m_{ij}=m_{ji}=\left\{\begin{array}{cl} 
-2 & if \  j=i+m_{ii}-1 \quad{\rm and}\quad i=j+m_{jj}-1 \quad{\rm mod}\ n\\
-1 & if\  j=i+m_{ii}-1 \quad{\rm or \ else}\quad i=j+m_{jj}-1 \quad{\rm mod}\ n\\
0 & {\rm in\  all\  other\  cases}
  \end{array}\right.$$
\end{description}
\begin{Th}[D1,N1] 1) Let $S$ be a minimal complex compact surface containing a GSS with $n=b_2(S)>0$, then $S$ contains $n$ rational curves $D_0,\ldots,D_{n-1}$ and there exists $\s$ such that the intersection matrix $M(S)$ of the rational curves in $S$ satisfy
$$M(S)=-M(\s).$$
Moreover the curve $D_i$  is non-singular if and only if $a_i\neq n+1$.\\
Conversely, for any $\s$ there exists a surface $S$ containing a GSS such that $M(S)=-M(\s)$.\\
2) For any $\s\neq r_n$, $M(\s)$ is positive definite.
\end{Th}
\begin{Ex}{\rm 1) For $\s=r_n$, $M(\s)$ is not positive definite. The dual graph of the curves has $n$ vertices
\begin{center}
\includegraphics[width=4cm]{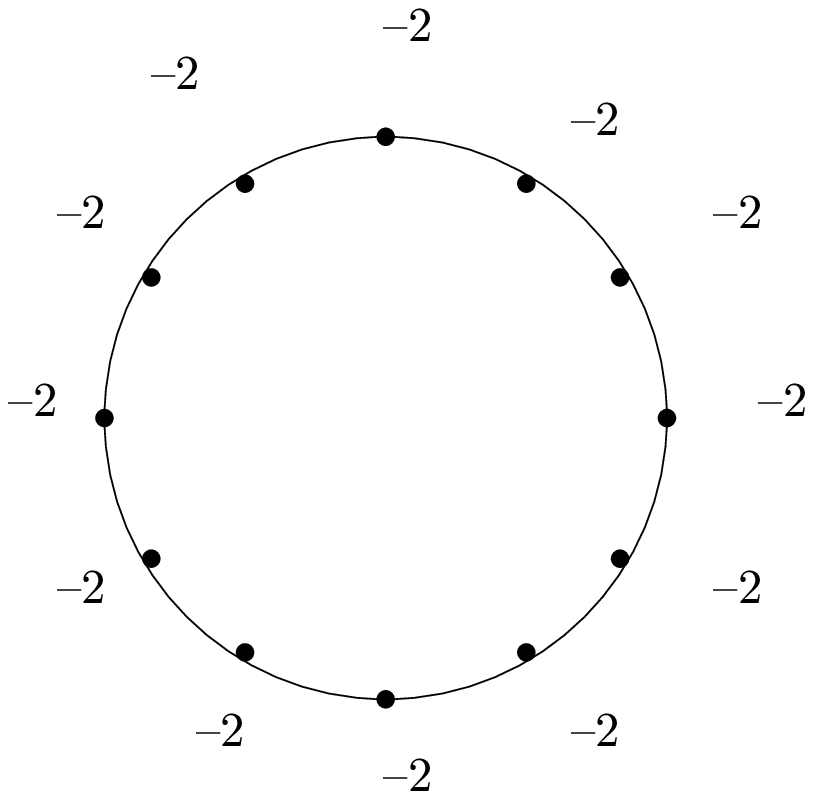}
\end{center}
This configuration of curves appear on Enoki surfaces \cite{E,N1,D1}.\\
2) If $\s=s_{p_0}\cdots s_{p_{N-1}}$ we obtain respectively one or two cycles if $N$ is odd (resp. even). The singularities are cusps and surfaces are odd (resp. even) Inoue-Hirzebruch surfaces \cite{I,N1,D1}. When there are two cycles, one of the two cycles determines the other. For example, if $\s=s_{p_0}s_{p_1}s_{p_2}s_{p_3}$, we obtain a cycle with $p_1+p_3$ curves and another with $p_0+p_2$ curves.
\begin{center}
\includegraphics[width=10cm]{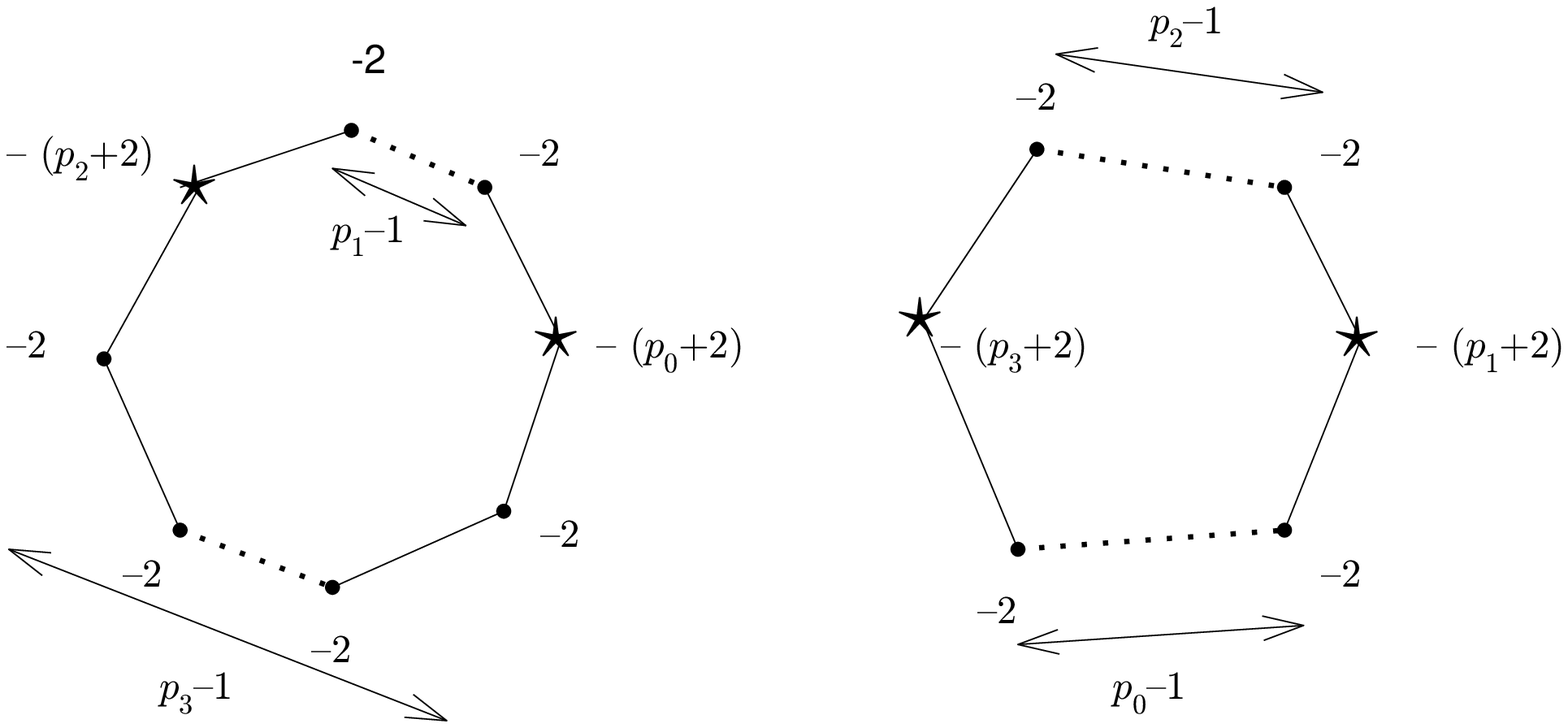}
\end{center}
3) The intermediate case \cite{N1,D1,DO}. There are branches and the number of branches is equal to the number of regular sequences in $\s$. For example, if $\s=r_{p_0}s_{p_1}$ the dual graph is
\begin{center}
\includegraphics[width=12cm]{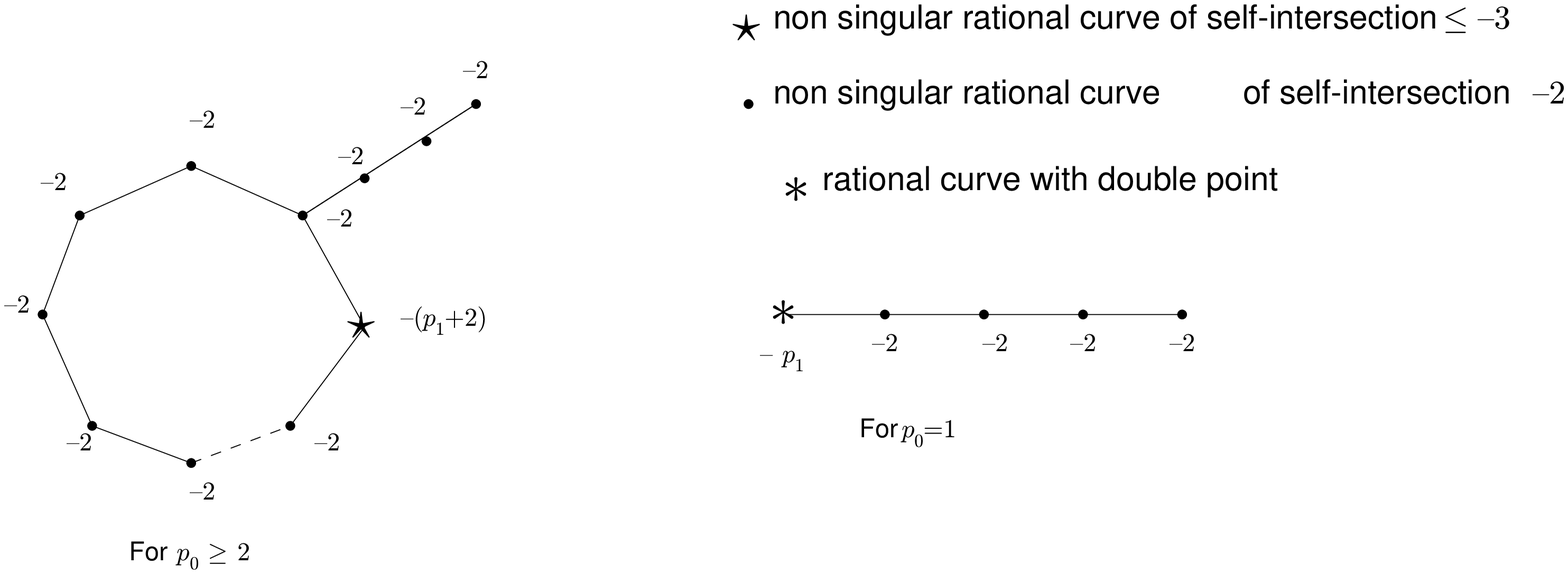}
\end{center}}
\end{Ex}

\section{Normal singularities associated to surfaces with GSS}
\subsection{Genus of the singularities}
If $S$ is a Inoue-Hirzebruch surface we obtain  by contraction
of a cycle, a singularity
called a cusp. They appear also in the compactification of Hilbert modular surfaces \cite{H}. We are interested here in the general
situation of any surface containing a GSS.\\

\begin{Prop}\label{cycle} Let $S$ be a compact complex surface of class VII$_0$ without non 
constant meromorphic. It is supposed that $n:=b_2(S)>0$, the maximal divisor $D$ is not trivial and the 
intersection matrix $M(S)$ is negative definite. Denote by $\Pi:S\to \bar S$ the contraction of the curves onto
isolated singular points. Then the following properties are equivalent:\\
i) $D$ contains a cycle of rational curves;\\
ii) $H^1(\bar S, {\mathcal O}_{\bar S})=0.$
\end{Prop}

Proof: {\bf i) $\Rightarrow$ ii)}  By Proposition \ref{SS}, the sequence
$$ 0\rightarrow  H^{1}(\bar S, {\cal 
O}_{\bar S})
\rightarrow  H^{1}( S, {\cal O}_{S})
\rightarrow H^{0}(\bar S, R^{1}\Pi _{\star}{\cal O}_{S})\rightarrow  
H^{2}(\bar S,
{\cal O}_{\bar S})\rightarrow 0.\leqno{(\ast)}$$
is exact. If $D$ contains a cycle then $h^{0}(\bar S, R^{1}\Pi _{*}{\cal O}_{S})\geq 1$. We suppose that 
$ h^{1}(\bar S, {\cal O}_{\bar S})=1$ and we shall derive a contradiction. With these assumptions,
$h^0(\bar S,\omega_{\bar S})=h^2(\bar S,{\mathcal O}_{\bar S})=1$ and $h^{0}(\bar S, R^{1}\Pi _{\star}{\cal O}_{S})= 1$
since $\bar S$ has no non constant meromorphic functions. Denote by  $x_i$, $i=0,\ldots, p$  the singular points of $\bar S$, 
$\Gamma_i=\Pi^{-1}(x_i)$ and 
$p_g(S,x_i)$ the geometric genus of $(S,x_i)$. Then $\sum p_g(S,x_i)=h^{0}(\bar S, R^{1}\Pi _{\star}{\cal O}_{S})= 1$, therefore
there are  rational singular points and one  elliptic  singular point. Moreover these singularities are Gorenstein because $h^0(\bar S,\omega_{\bar S})=1$ and a non-trivial section cannot vanish because there is no more curves. Hence there are rational double points with trivial
canonical divisor and one minimal elliptic singularity, $(S,x_0)$ with canonical divisor $\Gamma_0$ which is a cusp. Since there is a global meromorphic 2-form
on $S$, $n=-K_S^2=-\Gamma_0^2$. By \cite{N1}, $S$ is an odd Inoue-Hirzebruch surface (i.e. with one cycle); but such a surface has no canonical
divisor (see for example \cite{D2})\ldots a contradiction.\\
{\bf ii) $\Rightarrow$ i)}  By the exact sequence ($\ast$), $ h^{0}(\bar S, R^{1}\Pi _{*}{\cal O}_{S})\leq 2$ 
without any assumption and $1\leq h^{0}(\bar S, R^{1}\Pi _{*}{\cal O}_{S})$ by ii). Therefore there is a singular point, say $(S,x_0)$ such 
that $p_g(S,x_0)\geq 1$. If $\Gamma_0$ would be simply connected, then taking a 3-cover space $S'$ of $S$ we would obtain 3 copies of $\Gamma_0$ hence $ h^{0}(\bar S', R^{1}\Pi _{*}{\cal O}_{S'})\geq 3$ which is impossible since $S'$ remains in the VII$_0$-class, has no non-constant meromorphic functions and has to satisfy
$ h^{0}(\bar S', R^{1}\Pi _{*}{\cal O}_{S'})\leq 2$. $\Box$

\begin{Lem} \label{SEx}Let $S$ be a  surface with a GSS and such that 
$b_2(S)>0$.
Let $D$ be
the maximal divisor of $S$ and $\Pi:S\to \bar S$ be the 
contraction of $D$.  Then the sequence
$$ 0\rightarrow  H^{1}( S, {\cal O}_{S})
\rightarrow H^{0}(\bar S, R^{1}\Pi _{*}{\cal O}_{S})\rightarrow  
H^{2}(\bar S,
{\cal O}_{\bar S})\rightarrow  0$$
is exact and we have
$$1\leq h^{0}(\bar S, R^{1}\Pi _{*}{\cal O}_{S})=h^{0}(\bar S,
{\omega}_{\bar S}) +1 \leq 2 \leqno (\dag).$$
\end{Lem}
Proof:  By Proposition \ref{cycle} we have the desired exact sequence. 
Since $S$
has no non constant meromorphic functions,  the dimension of  
$H^{0}(\bar S, {\omega}_{\bar S})$
is 0 or 1.$\Box$\\
The proof of the following theorem follows the arguments of 
\cite{MER} Corollaire.
\begin{Th} Let $S$ be
a  surface with a GSS such that $2n<\sigma_{n}(S)\leq 3n$. Let $C$ be a connected component of the maximal 
divisor $D$ and let 
\hbox {$\Pi :S\rightarrow\bar S$} be the
contraction of  $C$, $\{x\}=\Pi(C)$. Then:\\
\noindent 1)  $p_{g}(\bar S, x) = 1 \, or \, 2$.\\
\noindent 2) If $2n<\sigma_{n}(S)<3n$ then $|D|$ is connected and the following 
conditions are equivalent:\par
i) $p_{g}(\bar S, x) = 2$\par
ii) the dualizing sheaf of $\bar S$ is trivial i.e. $\omega_{\bar 
S}\simeq
{\cal O}_{\bar S}$\par 
iii) the anticanonical bundle $-K$ is defined by a positive divisor 
$\Gamma$ i.e.
$\omega_{S}\simeq {\cal O}_{S}(-\Gamma )$ where $\Gamma > 0$.\par
iv) $(\bar S, p)$ is a Gorenstein singularity.\\
\noindent 3)  If S is an even Inoue-Hirzebruch surface, each cycle 
gives a minimally elliptic singularity and
the dualizing sheaf of $\bar S$ is trivial. In particular 
singularities are Gorenstein.\hfill\break
\noindent 4) If S is an odd Inoue-Hirzebruch surface the cycle gives 
a minimally elliptic singularity but 
the dualizing sheaf of $\bar S$ is not trivial. The singularity is 
still Gorenstein.\hfill\break
\end{Th}
Proof: 1) A connected component contains a cycle and we apply Lemma \thesection.\ref{SEx}. \\
2)   $i)\Longleftrightarrow ii)$: Notice that a global section of 
$\omega _{\bar S}$ cannot vanish 
since there is no curve.   Therefore by $(\dag)$ $p_{g}(\bar S, p) = 
2$ 
if and only if $\omega _{\bar S}$ is trivial.\\
$ii)\Rightarrow iii)$ By Lemma \ref{Gor}.\\
$iii)\Rightarrow ii)$ Let $\bar U=\bar S\setminus\lbrace x\rbrace$, 
$U=\Pi^{-1}(\bar U)$ and
$i:\bar U \lhook\joinrel\rightarrow \bar S$ the inclusion. We have since 
$\bar S$ is normal
$$\omega_{\bar S}=i_{*}\omega_{\bar U}\simeq 
i_{*}\Pi_{*}\omega_{U}\simeq
 i_{*}\Pi_{*}{\cal O}_{U}\simeq i_{*}{\cal O}_{\bar U}\simeq {\cal 
O}_{\bar S}$$
 Trivially $ii)\Rightarrow iv)$, we shall  proove $iv)\Rightarrow 
i)$. In fact, suppose 
that $p_{g}(\bar S,x)=1$, then by  \cite{L2} theorem 3.10, the 
singularity would be 
minimally elliptic, but it is impossible since in the case 
$2n<\sigma_{n}(S)<3n$ the maximal
 divisor contains a cycle with at least one branch \cite{D1} p113.\\
4) Suppose that $S$ is an even Inoue-Hirzebruch surface 
then   the sheaf $R^{1}\Pi _{*}{\cal O}_{S}$ is supported by two 
points. By $(\dag)$  and
Proposition 1.\ref{SRat},  $h^{0}(\bar S, R^{1}\Pi _{*}{\cal 
O}_{S})=2$,  and both singularities are minimally elliptic (see 
\cite{L2} p 1266). \\
5) It is 
well known  (\cite{I} or \cite{D2} Prop.2.14) that the canonical   
line bundle $K$ 
of an odd Inoue-Hirzebruch  surface is not
given by a divisor. The surface $S$ admits 
a double covering  by an even 
Inoue-Hirzebruch surface. By 4) the singularity is minimally elliptic and Gorenstein.\hfill $\Box$ 
\begin{Rem}{\rm Conditions i) and iv) are local conditions, though ii) 
and 
iii) are global ones.}
\end{Rem}

\subsection{$\bb Q$-Gorenstein and numerically Gorenstein singularities}
\begin{Def} Let $D$ be a connected exceptional divisor in $X$ and $\Pi:X\to \bar X$ the contraction onto $x=\Pi(D)\in\bar X$. Then $(\bar X,x)$ is a {\bf numerically Gorenstein} (resp. {\bf $\bb Q$-Gorenstein}) singularity if the positive numerically anticanonical $\bb Q$-divisor $D_{-K}$ is a divisor (resp. there exists an integer $m$ and a spc neighbourhood $U$ of $D$ such that the $m$-anticanonical bundle $K^{-m}$ has a section on $U$).
\end{Def}

If $S$ contains a GSS, then the fundamental group satisfies $\pi_1(S)=\bb Z$. Any topologically trivial line bundle is in $H^1(S,\bb C^\star)\simeq \bb C^\star$ and given by a representation of $\pi_1(S)$ in $\bb C^\star$. Therefore we shall denote  topologically trivial line bundles by $L^\a$ for $\a\in\bb C^\star$. 

\begin{Prop} Let $S$ be a compact complex surface containing a GSS of intermediate type, i.e $2n<\s_n(S)<3n$, $\Pi:S\to \bar S$ the contraction of the maximal divisor and $x=\Pi(D)$ the singular point of $\bar S$. Then 
\begin{description}
\item{i)} $(\bar S,x)$ is numerically Gorenstein  if and only if there exists a unique $\kappa\in\bb C^\star$  such that 
$$H^0(S,K^{-1}\otimes L^\kappa)\neq 0,$$
\item{ii)}   $(\bar S,x)$ is $\bb Q$-Gorenstein  if and only if there exists an integer $m\ge 1$ such that $$ H^0(S,K^{-m})\neq 0.$$
\end{description}
\end{Prop}
{\bf Proof}: i) The sufficient condition is  evident and the necessary condition derives from \cite{DO} thm 4.5.\\
ii) The sufficient condition is evident. Conversely, suppose that there exists an open neighbourhood $U$ of $D$ with $0\neq \t\in H^0(U, K_U^{-m})$,  non vanishing outside the exceptional divisor. Since the curves are a basis of $H^2(S,\bb Q)$, $K_S^{-m}$ is numerically equivalent to a positive divisor. The exponential exact sequence for surfaces of class VII$_0$ (\cite{KO} I,p766 and I (14) p756), yields the exact sequence
$$1\to H^1(S,\bb C^\star)\to H^1(S,\cal O_S^\star)\stackrel{c_1}{\to}H^2(S,\bb Z)\to 0$$
where $\bb C^\star\simeq H^1(S,\bb C^\star)$. Therefore there exists a unique $\k\in\bb C^\star$ such that 
$$H^0(S,K_S^{-m}\otimes L^{\k})\neq 0.$$
 Let $0\neq \o\in H^0(S,K_S^{-m}\otimes L^{\k})$. Flat line bundles are defined by a representation of $\pi_1(S)$ in $\bb C^\star$. Since in the intermediate case the cycle $\G$ of rational curves fulfils $H_1(\G,\bb Z)=H_1(S,\bb Z)$, the restriction $H^1(S,\bb C^\star)\to H^1(U,\bb C^\star)$ is an isomorphism. Then $\t/\o\in H^0(U,L^{1/\k})$ may vanish or may have a pole only on the exceptional divisor. Since the intersection matrix is negative definite, $\t/\o$ cannot vanish and $L^{1/\k}_{\mid U}$ is holomorphically trivial, hence $\k=1$.\hfill$\Box$

\begin{Ex}
In the exemple \cite{DO} 4.9, there is a family of surfaces with two rational curves, one rational curve with double point $D_0$ and a non-singular rational curve $D_1$, $D_0^2=-1$, $D_1^2=-2$ and $D_0D_1=1$. The obtained singularity is Gorenstein elliptic for $\a=\pm i$ and deforms into non-Gorenstein numerically Gorenstein elliptic singularities for other values of the parameter $\a$.
\end{Ex}

\section{Discriminant of the singularities}
\subsection{A family $\frak P$ of polynomials}
  For an integer $N\geq 1$, we denote ${\bb
Z}/N{\bb Z} =\{
\dot 0, \dot 1,\dots  , \buildrel .\over{ {N-1}}\}$. Let $$A=
\{\dot a_{1},\dots , \dot a_{p}\}
\subset  {\bb Z}/N{\bb Z}$$ a subset with p elements, $0\leq p\leq N$. We
may suppose that we have
$$0\leq a_{1}<a_{2}<\dots <a_{p}\leq N-1$$ which allows to define a
partition ${\cal
A}=(A_{i})_{1\leq i\leq p}$ of ${\bb Z}/N{\bb Z}$, where 

$$A_1:=\{\dot k\in {\bb Z}/N{\bb Z}\mid 0\le k\le a_1 \ or\ a_p<k\le N-1\} $$
$$A_i:=\{\dot k\in {\bb Z}/N{\bb Z}\mid a_{i-1}<k\le a_i\} \ for
\ 2\le i\le p.$$
When $A=\emptyset$, $\cal A$ is the trivial partition and $A_1={\bb
Z}/N{\bb Z}$.\hfill\break
\begin{center}
\includegraphics[width=5cm]{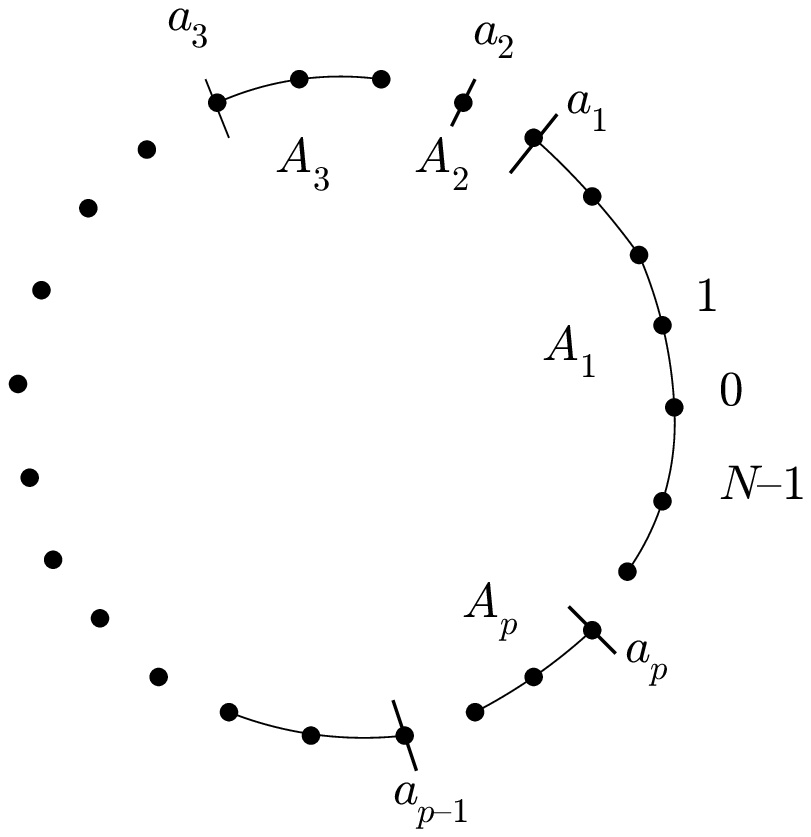}
\end{center}
\begin{Def} Let $N\geq 1$, $A\subset {\bb Z}/N{\bb Z}$ and
$B\subsetneq {\bb
Z}/N{\bb Z}$.\hfill\break
1) We shall say that $B$ is a {\bf generating allowed subset} relatively to $A$
if $B$ satisfies one of
the following conditions:\par
i)\ \ $B = \{\dot a\}$ \ with   $\dot a\in A$.\par
ii)\ $B = \{\dot k, \buildrel .\over{ {k+1}}\}$ and there exists
$1\leq i \leq p$ such that
$B\subset A_i$.\\
2) We shall say that $B$ is an  {\bf allowed subset} relatively to $A$ if $B$
admits a (possibly empty)
partition into generating allowed subsets.\\
The set of all allowed subsets will be denoted by ${\cal P}_A$.
\end{Def}

\begin{Def} \label{DefPol} For every $N\geq 0$, let $\frak P _N$ be the family of
polynomials defined in
the following way: $\frak P _{0} = \{0\}$.\\
If $N\geq 1$, $\frak P _{N} \subset {\bb Z}\lbrack X_0, \ldots ,
X_{N-1}\rbrack$ is the set of
polynomials $$P_{A}(X_0, \ldots , X_{N-1}) = \sum_{B\in {\cal P}_A}\quad
\prod_{i\notin
B}X_i \quad for \quad A\subset {\bf Z}/N{\bf Z}$$
We shall denote $$\frak P = \bigcup_{N\geq 0} \frak P_N$$
the union of all these polynomials.
\end{Def}

\begin{Ex}\label{Exemplesdepol} For $N=1$, there is only one polynomial $\frak P_{1 }= \{X\}$.\\
For $N=2$, 
$$\begin{array}{lcl}
\frak P_{2} &=& \bigl \{ P_{\emptyset}(X_0,X_1)=X_{0}X_{1},\
P_{\{0\}}(X_0,X_1)=X_{0}X_{1}+X_{1}, \\  
&&\\
&&\quad P_{\{1\}}(X_0,X_1)=X_{0}X_{1}+X_{0}, P_{\{0,
1\}}(X_0,X_1)=X_{0}X_{1}+X_{0} +X_{1}\bigr\}
\end{array}$$
For $N=3$, $\frak P_3$  contains the following polynomials
$$\left\{ \begin{array}{l}
P_{\emptyset}(X_{0}, X_{1}, X_{2}) =X_{0}X_{1}X_{2}+X_{0}+X_{1}+X_{2},\\
\\
  P_{\{0\}}(X_{0}, X_{1}, X_{2}) = X_{0}X_{1}X_{2}+X_{1}X_{2}+X_{0}+X_{2},\\
\\
P_{\{0, 1\}}(X_{0}, X_{1}, X_{2}) =
X_{0}X_{1}X_{2}+X_{1}X_{2}+X_{0}X_{2}+X_{1}+X_{2}, \\
\\
P_{\{0, 1, 2\}}(X_{0}, X_{1}, X_{2}) =
X_{0}X_{1}X_{2}+X_{1}X_{2}+X_{0}X_{2}+X_{0}X_{1}+X_{0}+X_{1}+X_{2}
\end{array}\right.$$ 
 and those obtained by circular
permutation of the variables.
\end{Ex}
Next proposition \ref{Premprop} gives the first properties of polynomials of
$\frak P$, lemma \ref{Annpartperm} shows
that  by vanishing of variables corresponding to an allowed subset, we
shall still obtain polynomials
of $\frak P$, proposition \ref{premier} shows that these polynomials are irreducible, finally proposition \ref{caracterisationpol} gives a characterization of the family $\frak P$.

\begin{Prop} \label{Premprop} 1) If $N\not= N'$, then $\frak P_{N}\bigcap \frak P_{N'} = \emptyset$\hfill\break
2) For $N\geq 2$,  the mapping
$$\begin{array}{ccc}
\frak P ({\bb Z}/N{\bb Z})&\longrightarrow&\frak P_N\\
A&\mapsto&P_A
\end{array}$$
is a bijection from the set $\frak P(\bb Z/N\bb Z)$ of subsets of ${\bb Z}/N{\bb Z}$ on $\frak P_N$. In
particular, if $N\geq
2$, $\frak P_N$ has $2^N$ elements.\\
3) If $A\subset {\bb Z}/N{\bb Z}$, then:
\begin{description}
\item{i)} $\deg P_A = N$ and $\displaystyle \prod_{i=0}^{N-1}X_i$ is the only
monomial of $P_A$ of
degree $N$.
\item{ii)} For $N\geq 2$, homogeneous part of $P_A$ of degree $N-1$ has $Card A$
monomials and these 
are $$\prod_{i\not= a}X_{i} \qquad for \ every \ a\in A$$
In particular homogeneous part of $P_A$ of degree $N-1$ determines $A$ and $P_A$
uniquely.
\item{iii)} $P_{A}(0) = 0$.
\end{description}
4) If $P(X_{0}, \ldots, X_{N-1})\in \frak P_{N}$ and $\a$ is a circular
permutation of $\{0,\ldots, N-1\}$ then $P(X_{\a (0)}, \ldots, X_{\a (N-1)})\in \frak P_{N}$.
\end{Prop}
{\bf Proof :} 1) derives from 3) i); 2) from 3) ii). Besides, the only
monomial of degree $N$ is
obtained for $B=\emptyset \in {\cal P}_A$, monomials of degree $N-1$ are
obtained for one
element subsets $\{a\}\in {\cal P}_A$. The integer  $N$ being fixed, these
monomials determine
$A$ and $P_A$. Finally, an allowed subset is by definition different from ${\bb Z}/N{\bb
Z}$, so we have the assertion 3) iii). Assertion 4) is evident.\hfill $\Box$

\begin{Lem and Def}\label{Partiespermisesfixes}
Let $A\subset {\bb Z}/N{\bb Z}$, ${\cal A}=(A_{i})_{1\leq i\leq p}$ the partition of  ${\bb Z}/N{\bb Z}$
defined by $A$ and let $B\in {\cal P}_{A}$.\\
1) Consider subsets of $B$ of the type $I=\{\buildrel .\over{
{j+1}}, \ldots , \buildrel
.\over{ {j+k}}\}$ such that:
\begin{description}
\item{i)} $\buildrel.\over{ {j+k}} \in A$,
\item{ii)} $I\subset B$ is maximal for inclusion,
\end{description}
Then $I$ is an allowed subset relatively to $A$ which will be called an
{\bf allowed subset fixed to $A$}.
The element $\dot j$ will be called the {\bf spring} of $I$.\\
\begin{center}
\includegraphics[width=7cm]{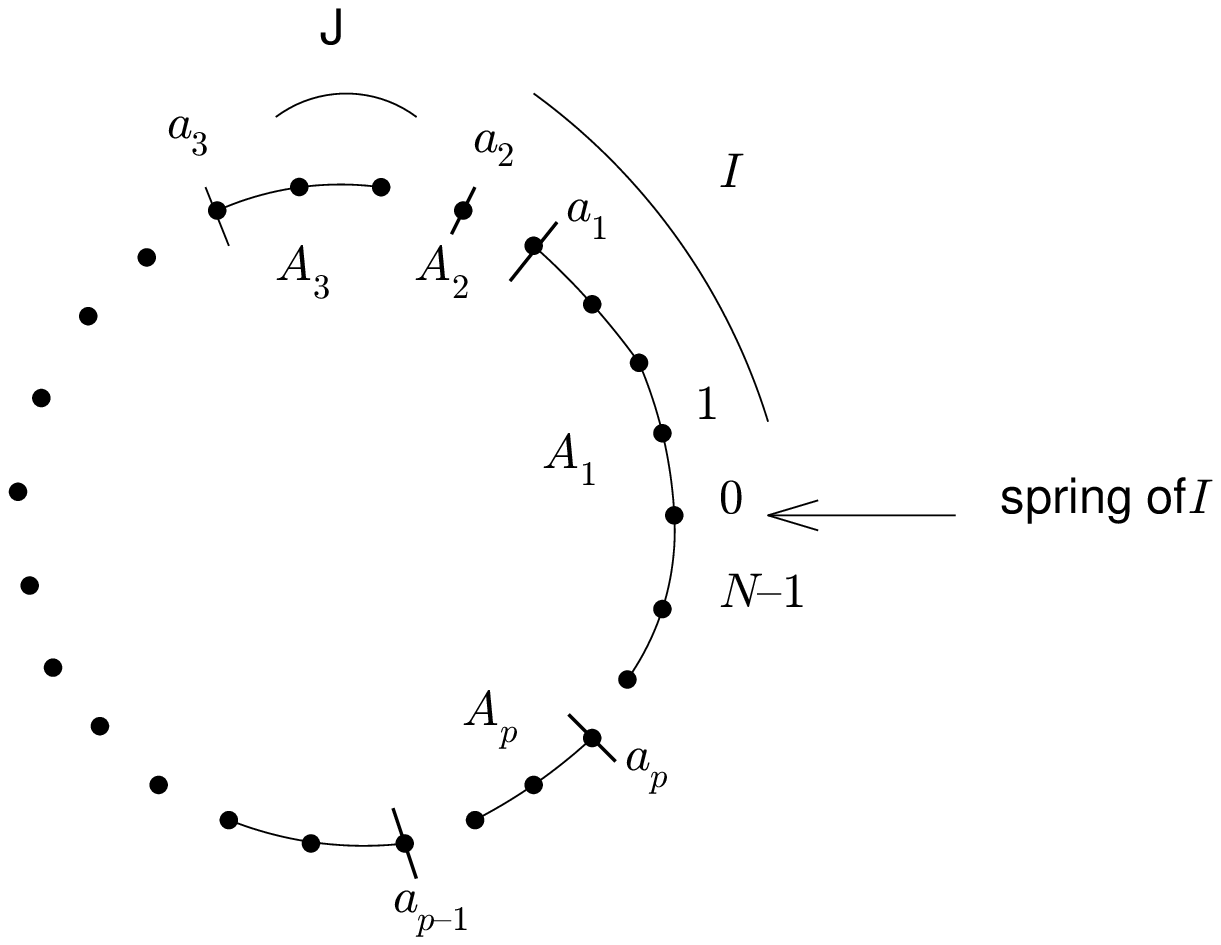}
\end{center}

2) Let $S_{B}$ be the set of springs of allowed subsets 
fixed to $A$, then we
have $S_{B}\cap B=\emptyset$.\hfill\break
3) Consider subsets of $B$ of the type $J=\{\buildrel .\over{
{j+1}}, \ldots , \buildrel
.\over{{j+2k}}\}$ such that:
\begin{description}
\item{i)} there exists $i$, $1\leq i\leq p$ such that $J\subset A_i$,
\item{ii)} $J\subset B$ is maximal for inclusion,
\item{iii)} For every allowed subset $I$, fixed to $A$, we have $J\cap
I=\emptyset$
\end{description}
then $J$ is an allowed subset relatively to $A$, which be called a
{\bf wandering allowed
subset}.\\
4) $B$ admits a unique partition by fixed allowed subsets and wandering
allowed subsets. This
partition will be called the {\bf canonical partition} of $B$.
\end{Lem and Def}
{\bf Proof: } clear.\hfill $\Box$ 

\begin{Rem}{\rm If $X\subset {\bb Z}/N{\bb Z}$ is not empty and $N'=\Card
 X$, canonical
action of ${\bb Z}/N{\bb Z}$ on itself induces an action of ${\bb Z}/N'{\bb
Z}$ on $X$, denoted
by $ +'$, defined in the following way: If $\dot x\in X$, let $j\geq
1$ be the least integer
such that $\buildrel .\over{{x+j}} \in X$; we set $\dot x 
+' \dot 1 =
\buildrel .\over{{x+j}}$.}
\end{Rem}

\begin{Lem}\label{Annulationdevariables}
Let $A\subset {\bb Z}/N{\bb Z}$,  and $B\in {\cal P}_{A}$. Denote by $B'$
the complement
of $B$ in  ${\bb Z}/N{\bb Z}$,  $N'=\Card B'$ and let $\varphi
:B'\rightarrow {\bb Z}/N'{\bb Z}$ a
bijection compatible with action of ${\bb Z}/N'{\bb Z}$ on $B'$ and on ${\bb
Z}/N'{\bb Z}$. If
$$A'=\varphi ( A\cap B')\bigcup  \varphi (S_{B})$$
where $S_B$ is the set of springs of $B$, then:\\
the mapping
$$\begin{array}{rccc}
\bar\varphi : &\{C\in  {\cal P}_{A}\mid  B\subset C \} & \longrightarrow
& \frak P  ({\bb Z}/N'{\bb
Z})\\
&&&\\
&C&\longmapsto &\varphi (C\cap B')
\end{array}$$
is a bijection from  $\{C\in  {\cal P}_{A} \mid  B\subset C \}$ on
${\cal P}_{A'}$.
\end{Lem}
{\bf Proof}: 1) $\bar\f$ is clearly injective.\\
2) Let $C\in \cal P_A$ such that $B\subset C$; to show that $\bar\f(C)\in\cal P_{A'}$, it is sufficient to show that if $I\subset C$ (resp. $J\subset J$) is an allowed subset fixed to $A$ (resp. a wandering allowed subset) belonging to the canonical partition of $C$, then $\f(I\cap B')\in \cal P_{A'}$ (resp. $\f(J\cap B')\in \cal P_{A'}$). On this purpose, we notice that if the last element of $I$ belongs to $A\cap B$, then  $\f(I\cap B')$ is an allowed subset with  last element in $\f(S_B)$; if the last element of $I$ is in $A\cap B'$, $\f(I\cap B')$ is an allowed subset with  the same last element in $\f(A\cap B')$. Therefore in both cases $\f(I\cap B')$ is an allowed subset fixed to $A'$. Besides, $J\cap A=\emptyset$ and $J\cap S_B=\emptyset$, hence $\f(J\cap B')$ is contained in an interval of the partition of $\bb Z/N'\bb Z$ associated to A'; $J$ has an even number of elements and $J\cap B'$ also. Finally,  $J\cap B'$ is a wandering allowed subset.\\
3) Let $C'\in \cal P_{A'}$ and $C=\f^{-1}(C')\cup B$. Then $C\in \cal P_A$, therefore $\bar\f$ is surjective.\hfill $\Box$

\begin{Lem} \label{Annpartperm}
 Let $P_{A} \in {\frak P}_{N}$,  $B\subset {\bb Z}/N{\bb Z}$
an allowed subset relatively to $A$, $B'$ the complement of $B$ in
${\bb Z}/N{\bb Z}$  and
$N'=\Card B'$. Then, identifying ${\bb Z}\lbrack X_{i}, i\in B' \rbrack$ with
${\bb Z}\lbrack X_{0}, \ldots , X_{N'-1}\rbrack$, there exists  $A'\subset
{\bb Z}/N'{\bb Z}$ such
that $$P_{A}(X_{i}=0,\ i\in B)\ =\ P_{A'}.$$
\end{Lem}
{\bf Proof :} In $P_{A}(X_{i}=0,\ i\in B)$ remain only monomials
$\displaystyle \prod_{i\notin
C}X_{i}$ of $P_A$ such that $B\subset C$; we then conclude by lemma
\ref{Annulationdevariables}.\hfill$\Box$

\begin{Prop}\label{premier}
 1) If $A=\emptyset$ and $N$ is even (resp. odd), $P_A$ has only
monomials of even (resp. odd) degrees.\\
2) If $N\geq 3$ and $P_{A} \in {\frak P}_{N}$, $P_A$ is irreducible in
${\bb Z}\lbrack X_{0}, \ldots , X_{N-1}\rbrack$.
\end{Prop}
{\bf Proof }: 1) If $B\in\cal P_A$ then $\Card B= 0$ mod 2.\\
2){\bf  First case}: $A=\emptyset$. The polynomial $P=P_A$ is invariant by circular permutation of the variables. Suppose that $P=P_1P_2$, with $P_j\in\bb Q[X_0,\ldots,X_{N-1}]$, $j=1,2$ and $P_1$ irreducible, $P_2\not\in\bb Q$. Fix a variable, say $X_i$, then $\deg_{X_i} P=1$; therefore the degree of one polynomial is zero and the degree of the other is one. Hence $P_1$ and $P_2$ have different variables. Denote by $I_j$, $j=1,2$ the subsets of indices $i$ such that $P_j$ depends on $X_i$. By  \ref{Premprop} 1) and 3),
$$P_j(X_i,i\in I_j)= \l_j \prod_{i\in I_j} X_i  \quad {\rm mod}\ (X_i,i\in I_j)^{\Card I_j -2} , \quad \l_1\l_2=1,$$
and $P_j$ contains only monomials the degree of which is of the same parity as $\Card I_j$, because $P_1$ and $P_2$ depend on different variables.\\
We show now that $P_1$ cannot depend on two consecutive variables: in fact, we could choose $X_i$ and $X_{i+1}$ in such a way that $P_1$ should not depend on $X_{i+2}$. However $P$ is stable by circular permutation, then
$$P(X)=P_1(X_i,i\in I_1) P_2(X_i,i\in I_2) = P_1(X_{i+1},i\in I_1) P_2(X_{i+1},i\in I_2)$$
where $P_1(X_{i+1},i\in I_1)$ is irreducible but cannot divide neither $P_1(X_i,i\in I_1)$ neither $P_2(X_i,i\in I_2)$, which is impossible since $\bb Q[X_0,\ldots,X_{N-1}]$ is factorial.\\
Finally we fix two consecutive indices $i\in I_1$ and $i+1\in I_2$. Then $\{i,i+1\}$ is an allowed subset, then by lemma \ref{Annpartperm},
$$P(X_i=X_{i+1}=0)\in \frak P_{N-2},$$
and by \ref{Premprop} 1), $\deg P(X_i=X_{i+1}=0)=N-2$. Then
$$\deg P_1(X_i=X_{i+1}=0)=\deg P_1(X_i=0)\leq \Card I_1 -2$$
$$\deg P_2(X_i=X_{i+1}=0)=\deg P_2(X_{i+1}=0)\leq \Card I_2 -2$$
which yields
$$N-2=\deg P_1(X_i=X_{i+1}=0) + \deg P_2(X_i=X_{i+1}=0)\leq N-4$$
a contradiction.\\
{\bf Second case}: $A\neq \emptyset$. We prove the result by induction on $N\geq 3$. The result for $N=3$ is true by example \ref{Exemplesdepol}. Let $N\geq 4$ and suppose, in order to simplify the notations, that $N-1\in A$. We have
$$P_A(X_0,\ldots,X_{N-1}) = X_{N-1}\bigl( P_A(X_{N-1}=1)-P_A(X_{N-1}=0)\bigr) + P_A(X_{N-1}=0),$$
with
$$\begin{array}{lcl}
Q(X_0,\ldots,X_{N-2})&:=&P_A(X_{N-1}=1)-P_A(X_{N-1}=0)\\
&&\\
&=&\prod_{i\neq N-1}X_i \quad {\rm mod}\ (X_0,\ldots,X_{N-2})^{N-2}
\end{array}$$
$$R(X_0,\ldots,X_{N-2}):=P_A(X_{N-1}=0)=\prod_{i\neq N-1}X_i \quad {\rm mod}\ (X_0,\ldots,X_{N-2})^{N-2}.$$
Since $\{N-1\}$ is an allowed subset for $A$, $R\in\frak P_{N-1}$ by \ref{Annpartperm}.
Now, by induction hypothesis, $R\in \bb Z[X_0,\ldots,X_{N-2}]$ is irreducible. By Eisenstein criterion, it is sufficient to prove that $R$ does not divide $Q$. But $\deg R=\deg Q=N-1$ and both polynomials have the same dominant monomial. Therefore we have to check that $R\neq Q$.
\begin{itemize}
\item If $A\neq \bb Z/N\bb Z$, we may suppose that $N-2\not\in A$ and $N-1\in A$, then $\{N-2,N-1\}\in\cal P_A$ and the monomial $M_{N-3}=\prod_{0\leq i\leq N-3}X_i $ is  in $P_A$, hence  in $R$, however $M_{N-3}X_{N-1}$ is not in $P_A$ hence not in $Q$.
\item If $A = \bb Z/N\bb Z$, $P_A$ contains $X_{N-1}$, therefore $Q(0,\ldots,0)=1$ though \\$R(0,\ldots,0)=0$.\hfill$\Box$
\end{itemize}

\begin{Rem}{\rm If $N=2$ second assertion of the preceeding proposition
is wrong as it can be
seen in example \ref{Exemplesdepol}.}
\end{Rem}

\begin{Prop}\label{caracterisationpol}
 Let $\displaystyle\frak P ' = \bigcup_{N\geq 0} \frak P_{N}' $ be a
family of polynomials where 
$$\frak P_{N}' \subset {\bf Q}\lbrack X_{0}, \ldots
, X_{N-1}\rbrack$$
satisfy the following conditions:
\begin{description}
\item{i)} For every $0\leq N\leq 2$, $ \frak P_{N}' \ =\ \frak P_{N}$,
\item{ii)} For every $N\geq 0$, $\Card \frak P_{N}'  =\Card \frak P_{N}$,
\item{iii)} If $P \in \frak P_{N}' $, then $\deg P=N$, and its homogeneous part of degree
$N$ is $\displaystyle
\prod_{0\leq i\leq N-1}X_i$,
\item{iv)} If $N\geq 3$ and $P \in \frak P_{N}' $, there exists $A=A_{P}\subset  {\bb Z}/N{\bb Z}$ such
that for every generating allowed subset $B\in {\cal P}_A$ we have
$$P(X_{i}=0,\ i\in B)\ \in \ \frak P_{N-\Card B}'.$$
Moreover, for every monomial $\lambda \prod_{i\notin C}X_{i}$ of $P$, where
$C\not=
\emptyset$ and $\lambda\in {\bb Q}$, there exists a generating allowed
subset $B$ such that
$B\subset C$.
\end{description}
Then, for every $N\geq 0$, $\frak P_{N} ' = \frak P_{N}$.
\end{Prop}
{\bf Proof}: We show by induction on $N\ge 2$ that $\frak P'_N=\frak P_N$. By i) let $N\ge 3$.\\
 Let $P\in \frak P'_N$. By condition iv), there exists $A=A_P\subset \bb Z/N\bb Z$ such that for every $B\in\cal P_A$
$$P(X_i=0, i\in B)\in \frak P_{N-\Card B}.$$
We are going to show that $P=P_A$. Both polynomials have the same dominant monomials $\prod_{0\leq i\leq N-1}X_i$. Let $B\in \cal P_A$ and $\prod_{i\not\in B}X_i$ one of the monmials of $P_A$. By iv), induction hypothesis and Proposition \ref{Premprop}, 3), 
$$P(X_i=0,i\in B)=\prod_{i\not\in B}X_i \quad {\rm mod}\ (X_i,i\not\in B)^{N-\Card B-1}$$
hence this monomials belongs to $P$ and by iii) each monomial of $P_A$ belongs to $P$. Conversely let $\l\prod_{i\not\in C}X_i$ be a monomial of $P$, let $B\in\cal P_A$ such that $B\subset C$. Denoting by $B'$ the complement of $B$ in $\bb Z/N\bb Z$ and $N'=\Card B'$, there exists $A'\subset \bb Z/N'\bb Z$ for which
$$P(X_i=0,i\in B)=P_{A'}$$
and by lemma \ref{Annulationdevariables}, $C\in \cal P_A$ and $\l=1$.\\
We have now, $\frak P'_N\subset \frak P_N$. We conclude by ii).\hfill$\Box$\\

To end this section we give a property of these polynomials which will allow to compute the discriminant of singularities whose exceptional set is associated to concatenation of sequences $\s=\s'\s''$.
\begin{Prop}\label{concatenation} Let $A'=\{a'_1,\ldots,a'_{p'}\}\subset \bb Z/N'\bb Z$ and $A''=\{a''_1,\ldots,a''_{p''}\}\subset \bb Z/N''\bb Z$. We identify $A'$ (resp. $A''$) with the subset of $\bb Z/N\bb Z$ (denoted in the same way) 
$$A'=\{a'_1,\ldots,a'_{p'}\}\subset \bb Z/N\bb Z \quad ({\rm resp.}\  A''=\{a''_1+N',\ldots,a''_{p''}+N'\}\subset \bb Z/N\bb Z )$$
 where
$$N=N'+N'', \quad and \quad 0\leq a'_1<\cdots<a'_{p'}<N'\leq a''_1+N'<\cdots < a''_{p''}+N'<N.$$
Setting $A=A'\cup A''\subset \bb Z/N\bb Z$ we have
$$\begin{array}{lcl}
P_A(X_0,\ldots,X_{N-1}) & = & P_{A'}(X_0,\ldots,X_{N'-1})P_{A''}(X_{N'},\ldots, X_{N'+N''-1})\\
&&\\
&&+ P_{A'}(X_0,\ldots,X_{N'-1}) + P_{A''}(X_{N'},\ldots, X_{N'+N''-1}).
\end{array}$$
\end{Prop}
{\bf Proof}: With the same identification as in the statement we have
$$\begin{array}{lcl}
\cal P_A &=& \bigl\{B'\cup B'' \mid B'\in\cal P_{A'}, B''\in\cal P_{A''}\bigr\}\\
&&\\
&& \bigcup \bigl\{B'\cup \{N',\ldots,N-1\} \mid B'\in \cal P_{A'}\bigr\} \bigcup \bigl\{ \{0,\ldots,N'-1\} \cup B''  \mid B''\in \cal P_{A''}\bigr\}
\end{array}$$
and this gives the three terms of the decomposition.\hfill$\Box$

\subsection{Main results}
\begin{Th}[Main theorem]\label{Maintheorem} Let $\s=\s_0\cdots\s_i\cdots\s_{N+\r}$ be a sequence of integers such that there are $N\geq 1$ singular sequences $\s_{i_j}=s_{k_j}$, $0\leq j\leq N-1$ and $0\leq \r\leq N$ regular sequences $r_{m_l}$, $0\leq l\leq \r-1$. Let $A \subset \bb Z/N\bb Z$ defined by
$$A=A(\s):=\left\{0\leq j\leq N-1  \mid \s_i\ {\rm is\ a\ regular\ sequence\ for}\ i=i_j +1 \quad {\rm mod}\ N+\r\right\}.$$
Then we have
$$\det M(\s)=P_A(k_0,\ldots,k_{N-1})^2.$$
\end{Th}
\begin{Cor} Let $S$ be a minimal surface containing a GSS with $n=b_2(S)$, rational curves $D_0,\ldots,D_{n-1}$ and intersection matrix $M(S)=(D_iD_j)=-M(\s)$. Then
\begin{description}
\item{i)} The index of the sublattice $\sum_{i=0}^{n-1}\bb Z D_i$ in $H_2(S,\bb Z)$ is
$$\left[H_2(S,\bb Z):\sum_{i=0}^{n-1}\bb Z D_i\right] = P_{A(\s)}(k_0,\ldots,k_{N-1});$$
\item{ii)} The curves $D_0,\ldots,D_{n-1}$ form a basis of $H_2(S,\bb Q)$ if and only if $\s\neq r_n$;
\item{iii)} The curves $D_0,\ldots,D_{n-1}$ form a basis of $H_2(S,\bb Z)$ if and only if $\s=s_1r_{n-1}=(3,2,\ldots,2)$ for $n\geq 1$ or $\s=s_1s_1=(3,3)$ if $n=2$.\\
 In this case we have the following matrices:
\begin{itemize}
\item $n=1$, $M(S)=-1$,
\item $n=2$, $M(S)=\left(\begin{array}{cc}-1&1\\1&-2\end{array}\right), \left(\begin{array}{cc}-1&0\\0&-1\end{array}\right)$,
\item $n\ge 3$,
$$\left(\begin{array}{ccccccc}
-3&0&1&0&\ldots&0&1\\
0&-2&1&0&\ldots&\ldots&0\\
1&1&-2&1&\ddots&&\vdots\\
0&0&1&\ddots&\ddots&\ddots&\vdots\\
\vdots&\vdots&\ddots&\ddots&\ddots&\ddots&0\\
0&0&&\ddots&\ddots&\ddots&1\\
1&0&\ldots&\ldots&0&1&-2
\end{array}\right)$$
\end{itemize}
\end{description}
\end{Cor}

\begin{Cor} Let $S$ be an even Inoue-Hirzebruch surface with  intersection matrices $M(S)=-M(\s)$ and $\s=s_{k_0}\cdots s_{k_{2q-1}}$.  Let  $\G$ and $\G'$ be the two cycles with intersection matrices $M(\G)$ and $M(\G')$, then
$$\begin{array}{lcl}
[H^2(\G,\bb Z):H_2(\G,\bb Z)]&=&|\det M(\G)|=P_{\bb Z/2q\bb Z}(k_0,\ldots,k_{2q-1})=|\det M(\G')|\\
&&\\
&=&[H^2(\G',\bb Z):H_2(\G',\bb Z)].
\end{array}$$
\end{Cor}
\subsection{A multiplicative topological invariant associated to singularities}

\begin{Def} A {\bf simple sequence} $\s$ is a sequence of the form 
$$\s=s_{k_0}\cdots s_{k_{N-1}}r_m$$
 with $N\ge 1$. A singularity is called {\bf simple} if it is obtained  by the contraction of a divisor whose dual graph is associated to a simple sequence. Of course any sequence $\s=\s_0\cdots\s_{N+\r-1}$ where $\s_i$ is singular or regular splits into $\r$ simple sequences.\\
The polynomial associated to any  singularity $(X,x)$ of type $\s$ is defined by
$$\D_\s(X_0,\ldots,X_{N-1}):=P_{A(\s)}(X_0,\ldots,X_{N-1})+1.$$
The integer $k=\D_\s(k_0,\ldots,k_{N-1})$ will be called the {\bf twisting coefficient} of the singularity.
\end{Def}

\begin{Lem} Let $\s'=s_{k'_0}\cdots s_{k'_{N'-1}}r_{m'}$ and $\s''=s_{k''_0}\cdots s_{k''_{N''-1}}r_{m''}$ be two  simple sequences. Then, denoting by $\s=\s'\s''$ the sequence obtained by concatenation of $\s'$ and $\s''$,  we have with $N=N'+N''$
$$\D_{\s'\s''}(X_0,\ldots,X_{N-1})=\D_\s'(X_0,\ldots,X_{N'-1})\D_{\s''}(X_{N'},\ldots,X_{N'+N''-1}).$$
\end{Lem}
{\bf Proof}: For $A'=A(\s')\subset \bb Z/N'\bb Z$, $A''=A(\s'')\subset \bb Z/N''\bb Z$, $N=N'+N''$, $A=A'\coprod A''\subset \bb Z/N\bb Z$, and $A=A(\s)$, we have by \ref{concatenation},
    $$\begin{array}{lcl}
    \D_{\s'\s''}(X_0,\ldots,X_{N'+N''-1})&=&P_A(X_0,\ldots,X_{N-1})+1\\
    &&\\
    &=& P_{A'}(X_0,\ldots,X_{N'-1})P_{A''}(X_{N'},\ldots,X_{N-1})\\
  &&\\  
    &&+P_{A'}(X_0,\ldots,X_{N'-1})+P_{A''}(X_{N'},\ldots,X_{N-1}) +1\\
  &&\\
    & =& (P_{A'}(X_0,\ldots,X_{N'-1})+1)(P_{A''}(X_{N'},\ldots,X_{N-1}) +1)\\
  &&\\
    &=&\D_{\s'}(X_0,\ldots,X_{N'-1})\D_{\s''}(X_{N'},\ldots,X_{N-1}).
    \end{array}$$
    \hfill$\Box$

Now we shall express the invariant $\D_\s$ for $\s$ simple, thanks to the determinant of the branch:

\begin{Lem}\label{invariantdebranche} Let $\s$ be a simple sequence with branch $B$ defined by 
$$Selfint(B)= \left\{
\begin{array}{l}
(\underbrace{2,\ldots,2}_{k_0-1},\  k_1+2,\  \underbrace{ 2,\ldots,2}_{k_2-1}, \ldots,\  k_{p-2}+2,\   \underbrace{2,\ldots,2}_{k_{p-1}-1},\  2) \\
\hfill  if\  p\equiv 1 ({\rm mod}\ 2)\\
\ \\
(\underbrace{2,\ldots,2}_{k_0-1},\  k_1+2,\  \underbrace{ 2,\ldots,2}_{k_2-1}, \ldots,\  k_{p-3}+2,\   \underbrace{2,\ldots,2}_{k_{p-2}-1},\  k_{p-1}+2) \\
\hfill  if\  p\equiv 0 ({\rm mod}\ 2)
\end{array}\right.$$
 then 
$$\D_\s(k_0,\ldots,k_{p-1})=\det B,$$
where $\det B$ is the determinant of the intersection matrix of the curves in $B$.
\end{Lem}
{\bf Proof}: For $p=1$, $Selfint(B)=(\underbrace{2,\ldots,2}_{k_0})$ and $\det B=k_0+1=P_\s(k_0)+1$.\\
For $p=2$, $Selfint(B)=(\underbrace{2,\ldots,2}_{k_0-1},k_1+2)$ and 
$\det B=k_0k_1+k_0+1=P_\s(k_0,k_1)+1$ (see Example \ref{Exemplesdepol}). By induction: we suppose that $p$ is odd, i.e. $p=2q+1$; the even case is  left to the reader. We shall use notations in (\ref{DefPol}). Since there is only one branch we have $\s=(s_{k_0}\cdots s_{k_{2q}}r_m)$, $N=2q+1$. For $A=\{2q\}\subset \bb Z/(2q+1)\bb Z$, we denote the allowed subsets by $\cal P_{2q+1}$. For the sequel we need the following observation: Let $C\in \cal P_{2q+1}$, then:
\begin{itemize}
\item if $2q\not\in C$ and $2q-1\not\in C$, $C\in\cal P_{2q-1}$ and $\sharp(C)$ is even;
\item if $2q\in C$ and $2q-1\not\in C$, $C=\{2q\}\cup C'$, $C'\in\cal P_{2q-1}$ and $\sharp(C')$ is even;
\item if $2q\not\in C$ and $2q-1\in C$, $C=\{2q-1,2q-2\}\cup C'$, $C'\in\cal P_{2q-1}$ and $\sharp(C')$ is even;
\item if $2q\in C$ and $2q-1\in C$, $C=\{2q,2q-1\}\cup C'$, $C'\in\cal P_{2q-1}$ and $\sharp(C')$ is odd or even.
\end{itemize}
Denote by $\D_{2q+1}$ the determinant of the branch when $\s$ contains  $2q+1$ singular sequences.
Applying lemma (\ref{DevN}) below, we have 
$$\begin{array}{lcl}
\D_{2q+1}(k_0,\ldots,k_{2q}) & = &(k_{2q}+1)\Bigl\{k_{2q-1}\D_{2q-1}(k_0,\ldots,k_{2q-2}-1) \\
&&\\
&&+  \D_{2q-1}(k_0,\ldots,k_{2q-2})\Bigr\} - k_{2q}\D_{2q-1}(k_0,\ldots,k_{2q-2}-1) \\
&&\\
&=&k_{2q}k_{2q-1} \D_{2q-1}(k_0,\ldots,k_{2q-2}-1)\\
&& + k_{2q}\Bigl\{ \D_{2q-1}(k_0,\ldots,k_{2q-2})-\D_{2q-1}(k_0,\ldots,k_{2q-2}-1)\Bigr\} \\
&&+ k_{2q-1}\D_{2q-1}(k_0,\ldots,k_{2q-2}-1)\\
&& + \D_{2q-1}(k_0,\ldots,k_{2q-2})
\end{array}$$
 
 In the sequel $\sum_{C'\in\cal P_{2q-1}}\prod_{i\not\in C'}k_i$ is shortened to $\sum_{C'\in\cal P_{2q-1}}$. Recall that $C'\in \cal P_{2q-1}$, i.e. $C'\subset \{\buildrel .\over{0},\ldots, \buildrel .\over{2q-2}\}=\bb Z/(2q-1)\bb Z$. By induction hypothesis,

$$\begin{array}{r}
\D_{2q+1}(k_0,\ldots,k_{2q}) 
 = \dps k_{2q}k_{2q-1}\left\{ \sum_{B'\in\cal P_{2q-1}\atop{2q-2\in B'}} \prod_{i\not\in B'}k_i + \sum_{B'\in\cal P_{2q-1}\atop{2q-2\not\in B'}} \left(\prod_{i\not\in B'\atop{i<2q-2}}k_i\right)(k_{2q-2}-1)+1\right\} \\

\dps+k_{2q}\left\{ \sum_{B'\in\cal P_{2q-1}\atop{2q-2\in B'}} + \sum_{B'\in\cal P_{2q-1}\atop{2q-2\not\in B'}} - \sum_{B'\in\cal P_{2q-1}\atop{2q-2\in B'}}  - \sum_{B'\in\cal P_{2q-1}\atop{2q-2\not\in B'}} \left(\prod_{i\not\in B'\atop{i<2q-2}}k_i\right)(k_{2q-2}-1)\right\}\hspace{10mm}\\
\dps+ k_{2q-1}\left\{\sum_{B'\in\cal P_{2q-1}\atop{2q-2\in B'}}+ \sum_{B'\in\cal P_{2q-1}\atop{2q-2\not\in B'}}\left( \prod_{i\not\in B'\atop{i<2q-2}}k_i\right)(k_{2q-2}-1)+1\right\}+ \sum_{B'\in\cal P_{2q-1}} +1 \hspace{12mm}\\
\\
=\dps k_{2q}k_{2q-1}\left\{ \sum_{2q-2\in B'}  + \sum_{2q-2\not\in B'} -  \sum_{2q-2\not\in B'}\left( \prod_{i\not\in\{2q-2\}\cup B'}k_i \right)+1\right\} 
+k_{2q}\left\{  \sum_{2q-2\in B'\atop{\sharp(B')\ odd}} \right\} \hspace{5mm} \\
\dps+ k_{2q-1}\left\{\sum_{2q-2\in B'}+ \sum_{2q-2\not\in B'} -  \sum_{2q-2\not\in B'} \left(\prod_{i\not\in \{2q-2\}\cup B'}k_i\right) +1\right\} 
+ \sum_{B'\in\cal P_{2q-1}} +1\hspace{8mm}\\
\end{array}$$
$$\begin{array}{lcl}
&\hspace{10mm}= &\dps k_{2q}k_{2q-1}\left\{ \sum_{2q-2\in B'\atop{\sharp(B')\ even}}  + \sum_{2q-2\not\in B'} +1 \right\} 
+ k_{2q}\left\{  \sum_{2q-2\in B'\atop{\sharp(B')\ odd}} \right\}\hspace{41mm}\\
&&\hspace{10mm}\dps+ k_{2q-1}\left\{ \sum_{2q-2\in B'\atop{\sharp(B')\ even}}
 + \sum_{2q-2\not\in B'} +1\right\}
  + \sum_{B'\in\cal P_{2q-1}} +1\hspace{41mm} \\
\\
&\hspace{10mm}=& \dps\sum_{B\in\cal P_{2q+1}\atop{2q\not\in B, 2q-1\not\in B}}  + \sum_{B\in\cal P_{2q+1}\atop{2q\not\in B, 2q-1\in B}} + \sum_{B\in\cal P_{2q+1}\atop{2q\in B, 2q-1\not\in B}}
 + \sum_{B\in\cal P_{2q+1}\atop{2q\in B, 2q-1\in B}}+1 \hspace{37mm}\\
 \\
&\hspace{10mm}=&\dps \sum_{B\in\cal P_{2q+1}}+1=P_\s(k_0,\ldots,k_{2q})+1.\hspace{66mm}
\Box

\end{array}$$

\begin{Prop}\label{produitdarbres} Let $\s=\s_0\cdots\s_{\rho-1}$ be a decomposition of $\s$ into simple sequences and let $B_0,\ldots,B_{\rho-1}$ be the branches of the dual graph, then
\begin{description}
\item{i)} $\D_\s=\prod_{i=0}^{\rho-1} \D_{\s_i}= \prod_{i=0}^{\rho-1} \det B_i$,
\item{ii)} $P_{A(\s)}=\prod_{i=0}^{\rho-1}(P_{A(\s_i)}+1) -1=\prod_{i=0}^{\rho-1} \det B_i   \ -1.$
\end{description}
(notice that different polynomials depend on different indeterminates).
\end{Prop}

\subsection{Twisted holomorphic 1-forms near the isolated singularity}
Let $S$ be a surface containing a GSS such that $b_2(S)=n$, with maximal divisor $D=\sum_{i=0}^{n-1}D_i$.  We assume that the intersection matrix $M(S)=-M(\s)$ is negative definite. Therefore we have
$$D=\G+\sum_{i=0}^{\rho-1}B_i,$$
 where $B_0,\ldots,B_{\rho-1}$ denote the branches of the dual graph.

\begin{Th} If $2n<\s_n(S)<3n$, then there exists a non-vanishing closed twisted logarithmic 1-form 
$$\o\in H^0(S,\O^1(Log D)\otimes L^k)$$
where the integer $k=k(S)\ge 2$ satisfies
$$k(S)=\prod_{i=0}^{\rho -1} \det B_i.$$
In particular in the complement of the singular point there is a non-trivially twisted non-vanishing holomorphic 1-form.
\end{Th}
We recal that in the notation $L^\a$, $\a\in\bb C^\star$ is the defining parameter of the topologically trivial line bundle.\\
Proof: By \cite{DO} p1537, there exists a global twisted logarithmic 1-form on  $S$ which does not vanish. The positive integer $k=k(S)$ is the integer which appears in any contraction $F$ associated to $S$ (see lemma 2.7 and thm 2.8 in \cite{DO}). By \cite{F} p480, the germ $F$ is conjugate to a germ of class 4 (conjugate by $(z,w)\mapsto (w,z)$ !) 
$$F(z,w)=(\mu zw^s + P(w), w^k),$$
and by \cite{F2} p35, we have 
$$\det M(S)=(-1)^n(k-1)^2.$$
With the main theorem \ref{Maintheorem} and Proposition \ref{produitdarbres} we conclude that
  $$k=\prod_{i=0}^{\rho - 1} \det B_i.$$
We obtain, in the complement of the singularity, a non-vanishing section on $\O^1\otimes \Pi_\star L^k$ where $ \Pi_\star L^k$ is not trivial since $M(S)$ is negative definite and $L^k$ has no non-vanishing sections on a neighbourhood of $D$.\hfill $\Box$

\section{Proof of the main theorem}

The aim is to compute the discriminant of the quadratic form using the family of polynomials previously introduced.\\
{\bf Sketch of proof}: 1)  When we compute the determinant  of $M(\s)$, a singular sequence $s_k=(k+2,2,\ldots,2)$
 produces a monomial containing $k^2$ because $k$ appears two times: one time because of the entry $k+2$ and a second time, with lemma  \ref{easycase}, because of the sequence 
 $$(\underbrace{2,\ldots,2}_{k-1}).$$
 By the same lemma, a regular sequence $r_m$ produces the integer $m$  at most at degree one. Therefore the determinant is a polynomial in the variables $k_0,\ldots,k_{N-1}$, and $m_0,\ldots, m_{\rho-1}$.
The idea is to develop the determinant splitting it into pieces which have a geometrical meaning. For example, consider $M=M(s_{k_0}r_m s_{k_1})$. Its dual graph is
\begin{center}
\includegraphics[width=7cm]{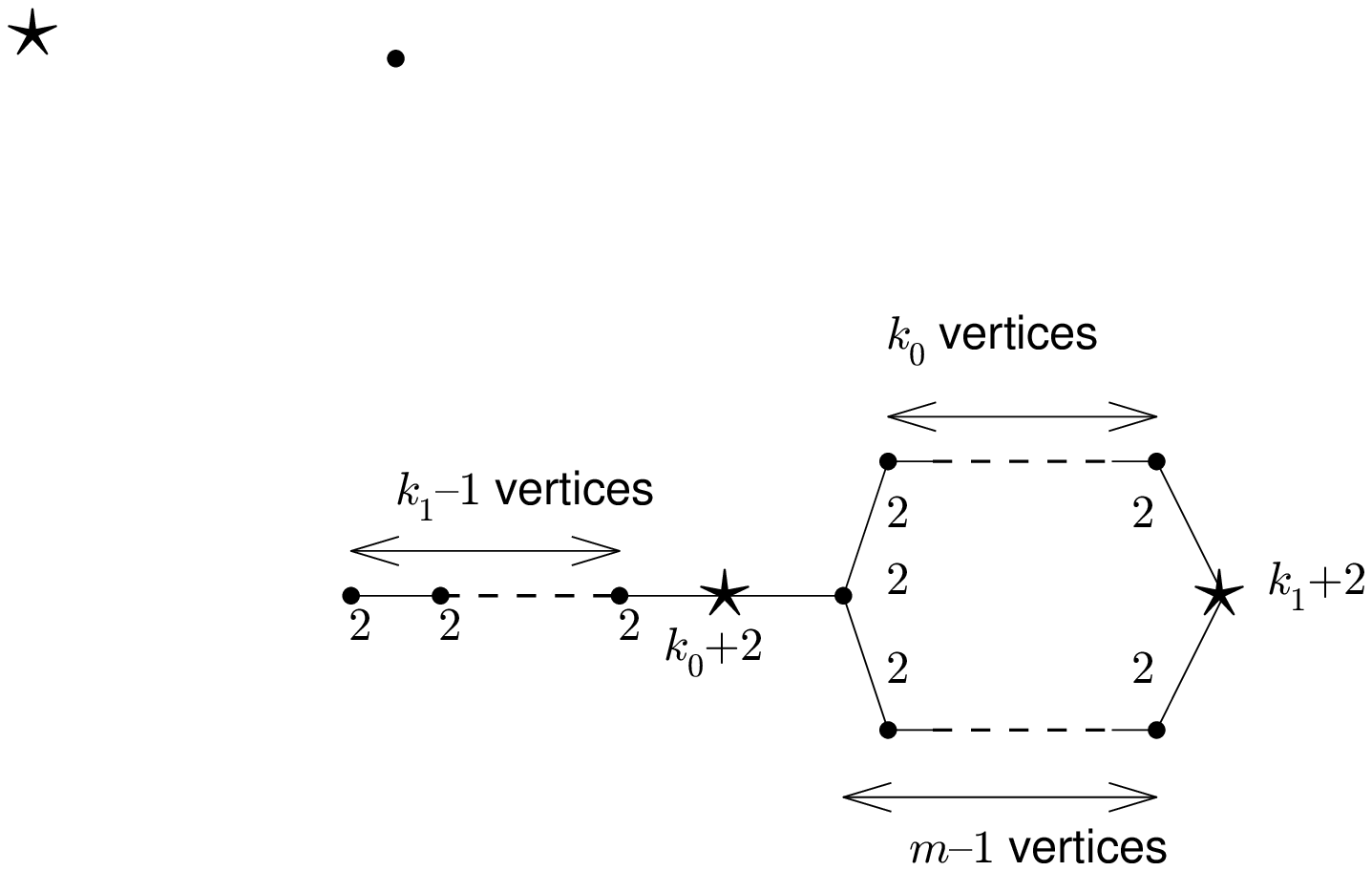}
\end{center}
 The vertices with weight $2$ are represented by a bullet, the vertices with weight $\ge 3$ are represented by a star. It splits into
\begin{center}
\includegraphics[width=12cm]{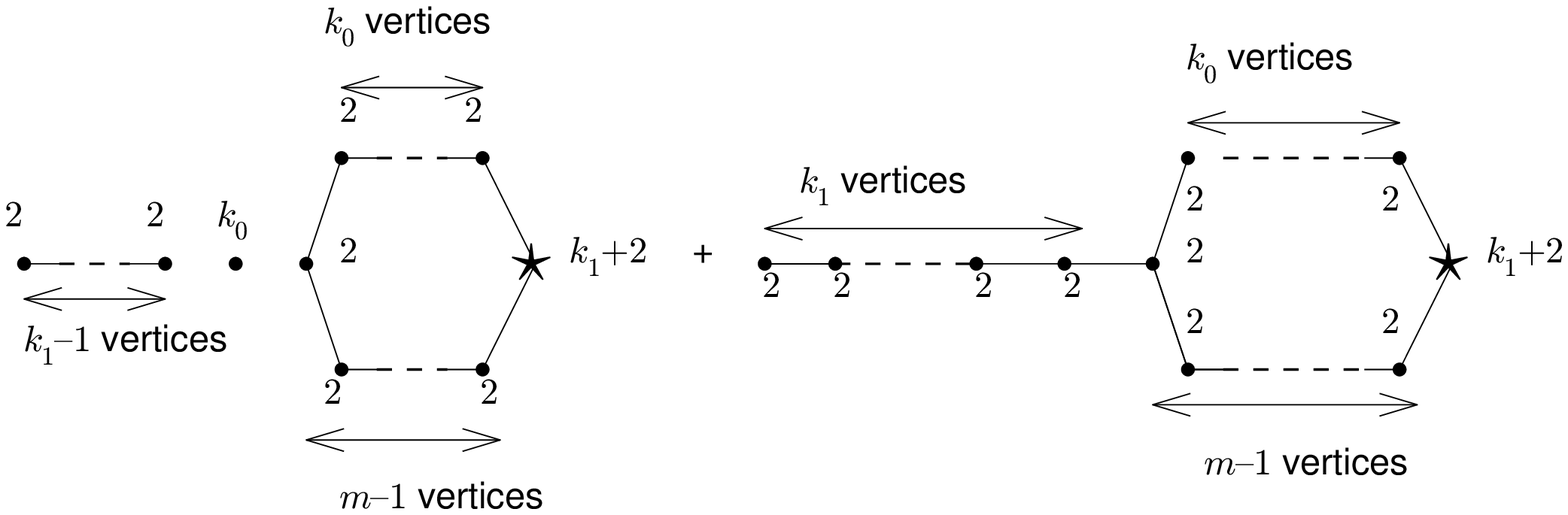}
\end{center}
which corresponds to the developement of the determinant along the $k_1$-th column by the splitting
$$\left(\begin{array}{c}\vdots\\-1\\k_0+2\\-1\\ \vdots\end{array}\right) = 
\left(\begin{array}{c}\vdots\\0\\k_0\\0\\ \vdots\end{array}\right) +
\left(\begin{array}{c}\vdots\\-1\\2\\-1\\ \vdots\end{array}\right),$$

$$\begin{array}{l}
\det M=\\
\\
\quad\left| \begin{array}{ccccccccccc}
2&-1&&&&&&&&&\\
-1&\ddots&\ddots&&&&&&&&\\
&\ddots&2&-1&&&&&&&\\
&&-1&k_0+2&-1&&&&&&\\
&&&-1&2&-1&0&\cdots&\cdots&0&-1\\
&&&&-1&\ddots&\ddots&&&&0\\
&&&&0&\ddots&2&\ddots&&&\vdots\\
&&&&\vdots&&\ddots&k_1+2&\ddots&&\vdots\\
&&&&\vdots&&&\ddots&2&\ddots&0\\
&&&&0&&&&\ddots&\ddots&-1\\
&&&&-1&0&\cdots&\cdots&0&-1&2
\end{array}\right|\scriptsize\begin{array}{l} 0\\ \\ \\ k_1-1\\ \\ k_1\\  \\  \\  \\ \\ \\ \\ \\ \\ \\ \\ \\ k_0+k_1+m-1\end{array}\\
\end{array}$$
$$=k_0 
\left|\begin{array}{cccc}2&-1&&\\
-1&\ddots&\ddots&\\
&\ddots&\ddots&-1\\
&&-1&2
\end{array}\right|{ \tiny\begin{array}{l} 0\\ \\ \\ \\  \\ \\ \\ k_1-2 \end{array}}
\left|\begin{array}{ccccccc}2&-1&&&&&-1\\-1&\ddots&\ddots&&&&\\
&\ddots&2&\ddots&&&\\
&&\ddots&k_1+2&\ddots&&\\
&&&\ddots&2&\ddots&\\
&&&&\ddots&\ddots&-1\\
-1&&&&&-1&2
\end{array}\right| { \tiny\begin{array}{l} k_1\\ \\ \\ \\ \\ \\ \\ \\ \\ \\ \\ \\  \\ \\ \\ k_0+k_1\\+m-1 \end{array}}$$
$$+ \left| \begin{array}{ccccccccccc}
2&-1&&&&&&&&&\\
-1&\ddots&\ddots&&&&&&&&\\
&\ddots&2&-1&&&&&&&\\
&&-1&2&-1&&&&&&\\
&&&-1&2&-1&0&\cdots&\cdots&0&-1\\
&&&&-1&\ddots&\ddots&&&&0\\
&&&&0&\ddots&2&\ddots&&&\vdots\\
&&&&\vdots&&\ddots&k_1+2&\ddots&&\vdots\\
&&&&\vdots&&&\ddots&2&\ddots&0\\
&&&&0&&&&\ddots&\ddots&-1\\
&&&&-1&0&\cdots&\cdots&0&-1&2
\end{array}\right|\scriptsize\begin{array}{l} 0\\ \\ \\ \\ k_1-1 \\ \\ k_1\\  \\  \\  \\ \\ \\ \\ \\ \\ \\ k_0+k_1\\+m-1\end{array}
$$
2) Noticing that a determinant is always a square of an integer, we prove that the determinant is in fact obtained as a square of a polynomial in $k_0,\ldots,k_{N-1}$ and hence {\it  the integers $m$ do not appear} in the development.\\
In the sequel we shall associate to $M$ a family of matrices obtained by this type of developement. The easy cases are those of a chain and  of a cycle with all diagonal entries equal to $2$:
\begin{Lem}\label{easycase}
Let $\d_m$ and $\D_m$ be the determinants of order $m\ge 1$ defined by
$$\d_1=2,\quad \d_2=\left|\begin{array}{cc}2&-1\\-1&2\end{array}\right|, \quad \D_1=0,\quad \D_2=\left|\begin{array}{cc}2&-2\\-2&2\end{array}\right|$$
$$\d_m= \left|\begin{array}{ccccc}2&-1&&&\\
-1&\ddots&\ddots&&\\
&\ddots&\ddots&\ddots&\\
&&\ddots&\ddots&-1\\
&&&-1&2
\end{array}\right|, \quad \D_m= \left|\begin{array}{ccccc}2&-1&&&-1\\
-1&\ddots&\ddots&&\\
&\ddots&\ddots&\ddots&\\
&&\ddots&\ddots&-1\\
-1&&&-1&2
\end{array}\right|, \quad m\ge 3,$$
then 
$$\d_m=m+1 \quad {\rm and}\quad \D_m=0.$$
\end{Lem}
{\bf Proof}: left to the reader.\hfill$\Box$

\begin{Lem}\label{DevN} Let $N=(n_{ij})_{0\leq i,j\leq p-1}$ be a matrix of order $p\ge 2$, of the form 
$$N=\left(\begin{array}{ccccccc}
2&-1&&&&&\\
-1&\ddots&\ddots&&&&\\
&\ddots&\ddots&\ddots&&&\\
&&\ddots&\ddots&-1&&\\
&&&-1&n_{mm}&\ldots&n_{m,p-1}\\
&&&   &\vdots   &         &\vdots\\
&&&   &n_{p-1,m}&\ldots&n_{p-1,p-1}
\end{array}\right)\scriptsize\begin{array}{c}0\\ \\ \\ \\ \\ \\ m-1\\ \\ m\\ \\ \\ p-1\end{array},$$
where $m\leq p-2$. For $J\subset \{0,\ldots,p-1\}$, we denote by $N_J=(n_{ij})_{i,j\in J}$ the submatrix whose entries depend on indices in $J$. Then
$$\det N=(m+1)\ \det N_{\{m,\ldots,p-1\}} -m\ \det N_{\{m+1,\ldots,p-1\}}.$$
\end{Lem}
{\bf Proof}: The result is trivial if $m=0$. If $m\ge 1$, development along the first column yields with induction hypothesis
$$\begin{array}{lcl}
\det N&=&2\det N_{\{1,\ldots,p-1\}} - \det N_{\{2,\ldots,p-1\}}\\
&&\\
&=&2\Bigl( m\  \det N_{\{m,\ldots,p-1\}}-(m-1)\ \det N_{\{m+1,\ldots,p-1\}}\Bigr)\\
&&-\Bigl( (m-1)\ \det N_{\{ m,\ldots,p-1\}} - (m-2)\ \det N_{\{m+1,\ldots,p-1\}}\Bigr)\\
&&\\
&=&(m+1)\ \det N_{\{m,\ldots,p-1\}} - m\ \det N_{\{m+1,\ldots,p-1\}}.
\end{array}$$
\hfill$\Box$
\subsection{Expression of determinants by polynomials}
\begin{Not}\label{Notationsmatrices} {\rm Let $N\geq 0$ and $\r\geq 0$ be integers such that $\r=1$ if $N=0$ and $\r\leq N$ if $N\geq 1$.\\
Let $M=M(\s)$ where $\s=\s_0\cdots\s_{N+\r-1}=(a_0,\ldots,a_{n-1})$ contains $N$ singular sequences $s_{k_i}$, $i=0,\ldots,N-1$ and $\r$ regular sequences $r_{m_j}$, $j=0,\ldots,\r-1$. Let
$$n=\sum_{i=0}^{N-1}k_i + \sum_{j=0}^{\r-1}m_j$$
be the order of $M=(m_{ij})_{0\leq i,j \leq n-1}$ or  the number of vertices of the associated dual graph.\\
We denote by $\cal C$ the set of subsets $J\subset\{0,\ldots,n-1\}$ which satisfy the following condition 
$$ \left\{\begin{array}{l} {\rm  let}\ 0\le l\le N+\r-1,\ {\rm and}\ \s_l=(a_r,\ldots,a_s).\\\\
{\rm If}\  \a \ {\rm satisfies}\ r+1\leq\a\leq s \ {\rm and}\  \a\in J \\
\\
{\rm then\  for\  all}\ \b \ {\rm such\  that}\  r+1\leq\b\leq s,
 \ {\rm we\  have}\  \b\in J.
 \end{array}\right.\leqno{(C)}$$
}
\end{Not}

 Splitting the graph into some pieces or changing the weights of some vertices, we associate to $M$ a family $\cal M$ of matrices in the following way:\\
 For $J\in\cal C$, let  $K_J$ defined by
 $$K_J=\left\{ j\in J \mid m_{jj}>2\right\}.$$
For $K\subset K_J$, denote by $M_J^K$ the matrix
$$M_J^K:=(m'_{ij})_{i,j\in J}$$
where 
$$\left\{\begin{array}{lcll} m'_{kk}&=&2&{\rm if}\ k\in K\\m'_{ij}&=&m_{ij}&{\rm in\ other\ cases}\end{array}\right.$$
The family $\cal M$ is
$$\cal M=\left\{M_J^K \mid J\in \cal C, K\subset K_J\right\}.$$

Now, for a fixed matrix $M_J^K$, we consider
\begin{itemize}
\item a partition $J=J'\cup J''$ of $J$, where $J'$ (resp. $J''$) is the subset of indices of vertices of the cycle (resp. of the branches), and
\item another partition of $J'$ and of $J''$ depending on $K$, composed of subsets of the following two types:
\begin{description}
\item{(1)} singletons $\{i\}$ such that $m_{ii}>2$,
\item{(2)} when elements of type $(1)$ are removed, connected components of vertices $j$ with weight  $m_{jj}=2$
\end{description}
\end{itemize}
To end, denote by $\nu_1(M_J^K)$ (resp. $\nu_2(M_J^K)$) the total number of subsets of type $(1)$ (resp. type $(2)$) in the partitions of $J'$ and $J''$ and we set
$$\nu(M_J^K)=\nu_1(M_J^K)+\nu_2(M_J^K).$$

\begin{Ex} Let $M=M(r_1s_1s_2)=M(2,3,42)$. Its dual graph is
\begin{center}
\includegraphics[width=4cm]{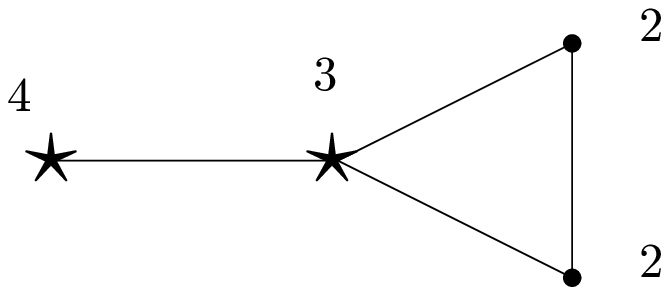}
\end{center}
\begin{itemize}
\item If $J=\{0,1,2,3\}$ and $K=\emptyset$ then $J'=\{0,1,3\}$, $J''=\{2\}$ and $\nu(M_J^K)=3$ and the dual graph of $M_J^K$ is
\begin{center}
\includegraphics[width=4cm]{2342}
\end{center}
\item If $J=\{0,1,2,3\}$ and $K=\{1,2\}$ then $\nu(M_J^K)=2$ and the dual graph of $M_J^K$ is
\begin{center}
\includegraphics[width=4cm]{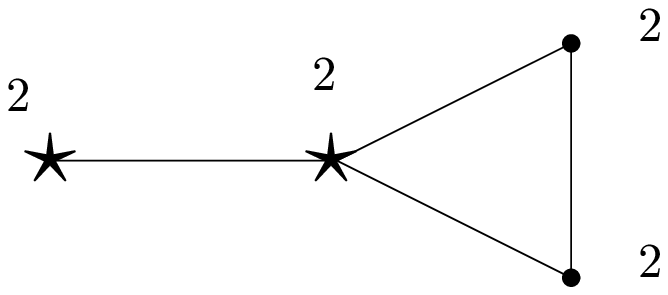}
\end{center}
\item If $J=\{1,2,3\}$ and $K=\{2\}$ then $J'=\{1,3\}$, $J''=\{2\}$, $\nu(M_J^K)=3$,
\begin{center}
\includegraphics[width=4cm]{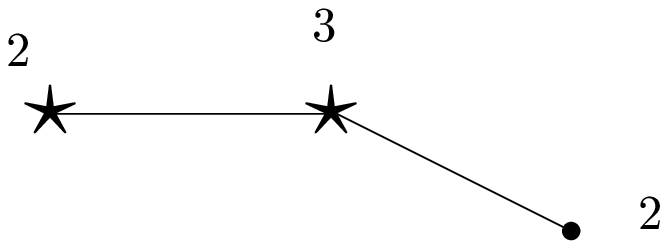}
\end{center}
\end{itemize}
\end{Ex}

\begin{Lem} \label{Existpoly} Let $M_J^K\in\cal M$, $\nu_1=\nu_1(M_J^K)$ and $\nu_2=\nu_2(M_J^K)$. Then there exists a polynomial
$$Q\in \bb Z[X_0,\ldots,X_{\nu_1-1},Y_0,\ldots,Y_{\nu_2-1}]$$
of degree $1$ respectively each variable, such that
$$\det M_J^K=Q(k_{i_0},\ldots,k_{i_{\nu_1-1}},m_0,\ldots,m_{\nu_2-1})$$
where subsets $\{i_j\}$ are of type $(1)$ and $m_j$ are the cardinals of the subsets of type $(2)$ which compose the partition of $J$.
\end{Lem}
{\bf Proof}:  By induction on $\nu=\nu_1+\nu_2\geq 1$. We have $\nu_1\leq N$ and $\nu_2\leq N+\r$ by condition $(C)$.\\
If $\nu=1$, either $\nu_1=1$, i.e. the determinant is of order $1$ and the result is clear, either $\nu_2=1$ and the results derives from lemma \ref{easycase}.\\
If $\nu\ge 2$, we may suppose that $M_J^K=(m'_{ij})$ is irreducible because reducible case is an immediate consequence of the induction hypothesis. Several cases may happen:\\
1) $M_J^K$ is a matrix of a cycle: Since $\nu\ge 2$ there exists an index $j\in J$ such that $m'_{jj}=k_{i_j}+2$. The decomposition of the $j$-th column
$$\left(\begin{array}{c} \vdots\\0\\-1\\k_{i_j}+2\\-1\\0\\ \vdots\end{array}\right) = 
 \left(\begin{array}{c} \vdots\\0\\0\\k_{i_j}\\0\\0\\ \vdots\end{array}\right) + \left(\begin{array}{c} \vdots\\0\\-1\\2\\-1\\0\\ \vdots\end{array}\right)$$
yields the relation

$$\det M_J^K = k_{i_j}\ \det M_{J\setminus\{j\}}^K + \det M_J^{K\cup\{j\}}\leqno{(\dag)}$$

where $M_{J\setminus\{j\}}^K$ (resp. $M_J^{K\cup\{j\}}$) is a matrix of a chain (resp. of a cycle). Setting
$$\nu'_i=\nu_i(M_{J\setminus\{j\}}^K), \quad \nu''_i=\nu_i(M_J^{K\cup\{j\}}), \qquad i=1,2,$$
we have
$$\left\{\begin{array}{ll} \nu'_1=\nu_1-1 & \nu'_2=\nu_2\\ &\\ \nu''_1=\nu_1-1&\nu''_2\leq \nu_2 \end{array}\right.$$
(there is one exception : when all entries of the cycle are $\ge 3$. We have $\nu''_1=\nu_1-1$ but $\nu''_2= \nu_2+1$ but then we repete the procedure ).\\
By induction hypothesis there exists polynomials
$$Q\in \bb Z[X_0,\ldots,\widehat{X_j},\ldots,X_{\nu_1-1},Y_0,\ldots,Y_{\nu_2-1}]$$
$$R\in \bb Z[X_0,\ldots,\widehat{X_j},\ldots,X_{\nu_1-1},Y_0,\ldots,Y_l,\widehat{Y_{l+1}},\ldots,Y_{\nu_2-1}]$$
such that, by a suitable numbering of the indices
$$
\det M_{J\setminus\{j\}}^K=
 Q(k_{i_0},\ldots,\widehat{k_{i_j}}, \ldots, k_{\nu_1-1}, m_0,\ldots,m_{\nu_2-1})
$$
$$\begin{array}{l}
\det M_J^{K\cup\{j\}} =\\
\hspace{1cm} R(k_{i_0},\ldots,\widehat{k_{i_j}}, \ldots, k_{\nu_1-1}, m_0,\ldots,\ldots,m_{l-1},m_l+m_{l+1}+1,m_{l+2},\ldots,m_{\nu_2-1}).
 \end{array}$$
We conclude replacing in $(\dag)$.\\
2) $M_J^K$ is not the matrix of a cycle: then the dual graph is a part of a cycle or contains bits of branches of $M$. In any cases, the dual graph contains a terminal vertex
\begin{center}
\includegraphics[width=5cm]{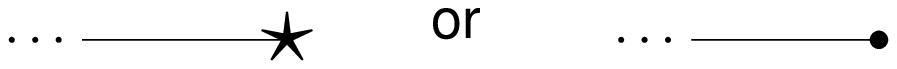}
\end{center}
\begin{itemize}
\item If in this chain there is a vertex with weight $>2$, we develop as before,
\item If not, all vertices have a weight equal to $2$, but since $\nu\ge 2$, this chain leads to a bifurcation
\begin{center}
\includegraphics[width=4cm]{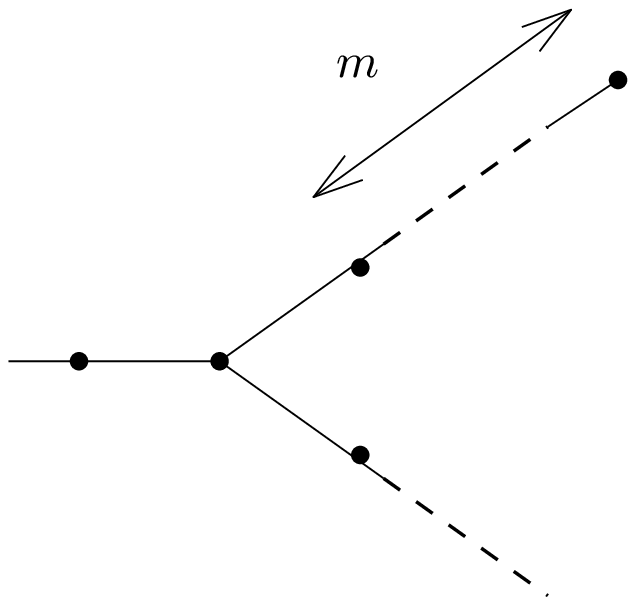}
\end{center}
Either the vertex of bifurcation has a weight $>2$ and we develop as before, either we apply lemma \ref{DevN} with appropriate numbering of entries of $M_J^K$:
$$\det M_J^K = (m+1)\ \det M_{J\setminus\{0,\ldots,m-1\}}^K - m\ \det (M_J^K)_{J\setminus\{0,\ldots,m\}}.\leqno{(\ddag)}$$
\end{itemize}
The matrix $(M_J^K)_{J\setminus\{0,\ldots,m\}}$ obtained by deletion of the branch with its root may not be in $\cal M$, however applying once again lemma \ref{DevN}, we obtain a matrix in $\cal M$ thanks to the explicit description of $M$ given by theorem 2.39 in \cite{D1}. We apply then induction hypothesis and $(\ddag)$.\hfill $\Box$

\begin{Lem}\label{lemme38} Let $M=M(\s_1\cdots\s_{N+\r})$ be a matrix satisfying notations \ref{Notationsmatrices}. Then, there exists a polynomial 
$$Q\in\bb Z[X_0,\ldots,X_{N-1},Y_0,\ldots,Y_{\r-1}]$$
 of degree at most $2$ (resp. $1$) relatively $X_i$, $i=0,\ldots,N-1$ (resp. $Y_j$, $j=0,\ldots,\r-1$) which satisfies
$$\det M=Q(k_0,\ldots,k_{N-1},m_0,\ldots,m_{\r-1}).$$
\end{Lem}
{\bf Proof}: We have $M=M_{\{0,\ldots,n-1\}}^\emptyset \in\cal M$ and by theorem \ref{Arbresetcycle}, $\nu_1=N$, $\nu_2\leq N+\r$ (with $\r=\r(S)$). More precisely (with notations of \ref{Arbresetcycle}), if there exists an integer $s$ such that $p_s=0 \ {\rm mod}\ 2$, we have for $t=s+1\ {\rm mod}\ \r$, the chain
\begin{center}
\includegraphics[width=8cm]{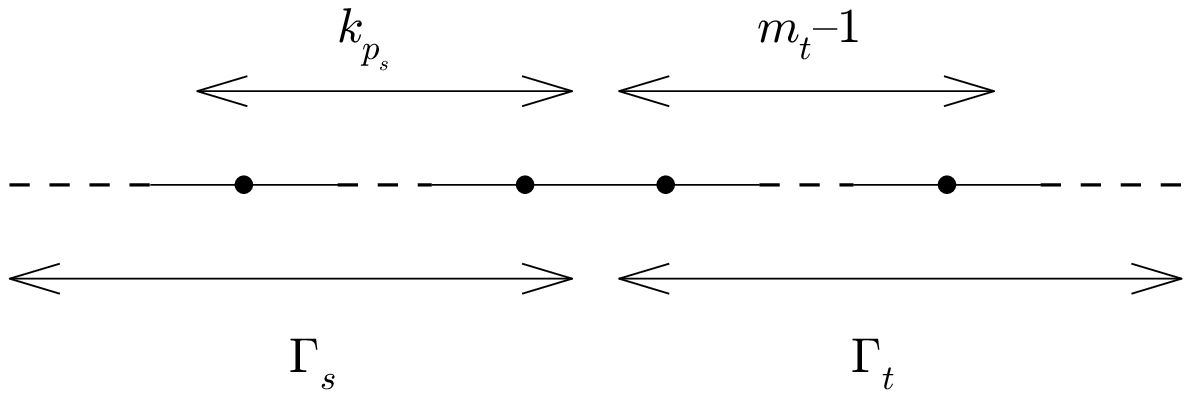}
\end{center}
Let $\frak S=\{t \mid p_s=0\ {\rm mod}\ 2,\ {\rm for}\ s=t-1\}$. Then
$$\nu_2=N+\r-\Card \frak S.$$
Lemma \ref{Existpoly} gives a polynomial in $2N+\r-\Card \frak S$ indeterminates
$$Q\in\bb Z[X_0,\ldots,X_{N-1},Y_0,\ldots,Y_{N-1},Y_N,\ldots,\widehat{Y_{N+t}},\ldots, Y_{N+\r-1} \mid t\in \frak S]$$
such that for suitable indices $i(t)\leq N-1$,
$$\begin{array}{l}
\det M(\s_0\cdots\s_{N+\r-1})=\\
\hspace{2cm} Q(k_0,\ldots,k_{N-1}, k_0,\ldots, k_{i(t)}+m_t,\ldots,k_{N-1},m_0,\ldots,\widehat{m_t},\ldots,m_{\r-1})
\end{array}$$
Setting $Y_i=X_i$ for $i\leq N-1$, $i\neq i(t)$, $t\in\frak S$, and substituing $X_{i(t)}+m_t$ in $Y_{i(t)}$ for $t\in\frak S$, the wished polynomial is obtained.\hfill $\Box$\\

Here is the key lemma for the reduction lemma of the following section:
\begin{Lem}\label{Polcarre} 1) Let $P,Q$ be two polynomials in $ \bb Q[X_0,\ldots,X_{n-1}]$. Suppose there exists an integer $N$ such that for $k_0\ge N, \ldots, k_{n-1}\ge N$ the following equality
$$P(k_0,\ldots,k_{n-1})=\pm Q(k_0,\ldots,k_{n-1})$$
holds. Then $P=Q$ or $P=-Q$.\\
2) Let $P\in\bb Q[X_0,\ldots,X_{n-1}]$ of degree at most $2$ relatively to each indeterminate. Suppose that there exists an integer $N$ such that for $k_0\ge N, \ldots, k_{n-1}\ge N$, $P(k_0,\ldots,k_{n-1})$ is the square of a rational. Then there exists $Q\in\bb Q[X_0,\ldots,X_{n-1}]$ satisfying
$$P=Q^2.$$
In particular, if $\deg_{\! X_i}P\leq 1$, $P$ does not depend on $X_i$.
\end{Lem}
{\bf Proof}: 1) By induction on $n\ge 1$. \\
2) The statement is true for $n=1$ without condition on the power by \cite{R}. Then by induction: suppose $n\ge 2$ and fix $k_0,\ldots,k_{n-2}\ge N$. Set
$$\begin{array}{lcl}
A(X_{n-1})&=&P(k_0,\ldots,k_{n-2},X_{n-1})\\
&&\\
&=&X_{n-1}^2 P_2(k_0,\ldots,k_{n-2})+X_{n-1}P_1(k_0,\ldots,k_{n-2})+P_0(k_0,\ldots,k_{n-2}).
\end{array}$$
For each $k_{n-1}\ge N$, $A(k_{n-1})$ is the square of a rational, hence by the  one variable case, $P_0(k_0,\ldots,k_{n-2})$ and $P_2(k_0,\ldots,k_{n-2})$ are squares of rationals. Induction hypothesis shows that there exists polynomials $Q_0, Q_1\in \bb Q[X_0,\ldots,X_{n-2}]$ which satisfy
$$P_0=Q_0^2, \quad {\rm and}\quad P_2=Q_1^2.$$
Replacing, one obtains
$$P_1(k_0,\ldots,k_{n-2})=\pm 2 Q_0(k_0,\ldots,k_{n-2})Q_1(k_0,\ldots,k_{n-2}).$$
By 1), one conclude that
$$P=(X_{n-1}Q_1 \pm Q_0)^2.$$
\hfill$\Box$

\subsection{The reduction lemma}
In this section we shall prove that the polynomial which computes a determinant depends on the positions of regular sequences in $\s$, however not  on their length.
\begin{Lem}[Reduction lemma]\label{reductionlemma} Let $M=M(\s_0\cdots\s_{N+\r-1})$ be a matrix which fulfils conditions \ref{Notationsmatrices}. Then, there exists a polynomial $P_\s\in \bb Q[X_0,\ldots,X_{N-1}]$ of degree at most $1$ relatively each indeterminate $X_i$, $i=0,\ldots,N-1$ such that
$$\det M(\s)=P_\s(k_0,\ldots,k_{N-1})^2.$$
In particular the determinant of $M$ does not depend on the length of the regular sequences.
\end{Lem}
{\bf Proof}: By lemma \ref{lemme38} there exists $Q\in\bb Q[X_0,\ldots,X_{N-1},Y_0,\ldots,Y_{\r-1}]$ such that
$$\det M= Q(k_0,\ldots,k_{N-1},m_0,\ldots,m_{\r-1}),$$
of degree at most $2$ in $X_i$ and at most $1$ in $Y_j$. The matrix $-M$ is an intersection matrix hence $\det M$ is the square of an integer by \ref{matriceintdetcarre} when $k_i\ge 1$, $i=0,\ldots,N-1$ and $m_j\ge 1$, $j=0,\ldots,\r$. Then lemma \ref{Polcarre} implies the existence of
$$P_\s\in\bb Q[X_0,\ldots,X_{N-1},Y_0,\ldots,Y_\r]$$
which satisfies $Q=P^2$. But $\deg_{Y_j}Q\leq 1$, therefore $P$ and $Q$ do not depend on $Y_j$.\hfill$\Box$
\subsection{Relation between determinants and polynomials of $\frak P$}
The next step is to prove that the polynomials $P_\s$ (see \ref{reductionlemma}) belong in fact in the family $\frak P$ previously defined. We shall apply the caracteristic properties of $\frak P$ given in \ref{caracterisationpol}. We start with examples.

\begin{Ex}\label{ExemplesPsigma} 1) Case $N=0$: then $M=M(\s)=M(r_m)$ and $\det M=0$. Therefore $P_\s=0$.\\
2) Case $N=1$: If $M=M(s_k)=\left(\begin{array}{ccccc}
k+2&-1&&&-1\\
-1&2&-1\\
&\ddots&\ddots&\ddots\\
&&\ddots&\ddots&-1\\
-1&&&-1&2
\end{array}\right)$, then 
$$\det M=k\d_{k-1}+\D_k=k^2, \quad{\rm and}\quad P_\s(X)=X.$$
If $M=M(s_kr_m)$ we have by the reduction lemma \ref{reductionlemma}
$$\det M=\det M(s_kr_1)=\left|\begin{array}{ccccc}
k&&&&-1\\&2&-1\\&-1&\ddots&\ddots\\&&\ddots&\ddots&-1\\-1&&&-1&2
\end{array}\right| = k\d_k - \d_{k-1}=k(k+1)-k=k^2,$$ and $P_\s(X)=X$.\\

\begin{center}
\includegraphics[width=5cm]{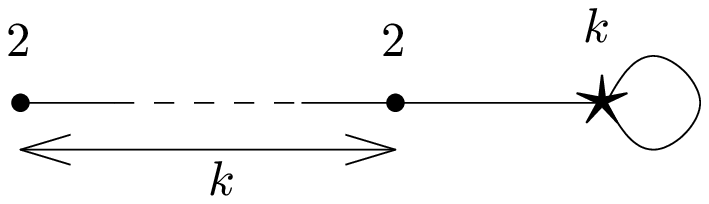}
\end{center}

3) Case $N=2$: If $M=M(s_{k_0}s_{k_1})$ the matrix is reducible and $\det M=(k_0k_1)^2$.

\begin{center}
\includegraphics[width=6cm]{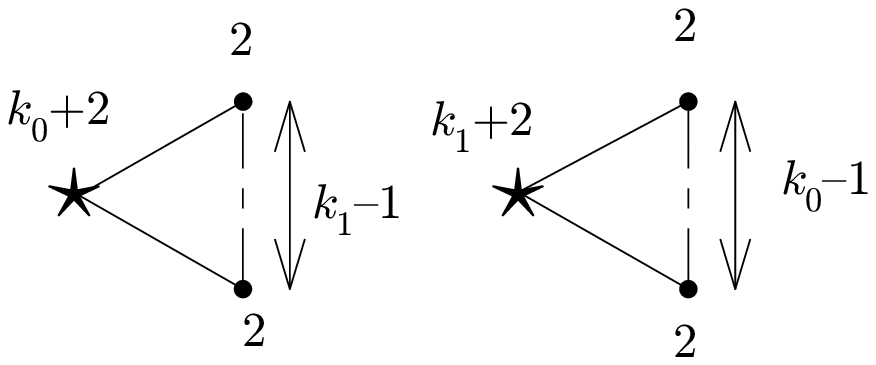}
\end{center}

If $M=M(s_{k_0}r_ms_{k_1})$ we have by \ref{reductionlemma},
$$\det M=\det M(s_{k_0}r_1s_{k_1})=(k_0k_1+k_1)^2$$

\begin{center}
\includegraphics[width=6cm]{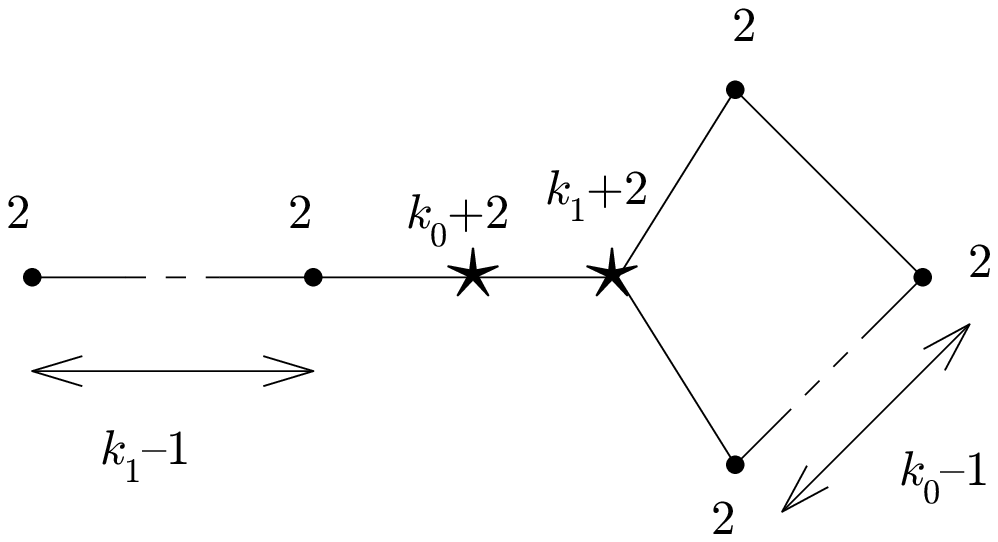}
\end{center}

If $M=M(s_{k_0}r_{m_0}s_{k_1}r_{m_1})$ we have
$$\det M = \det M(s_{k_0}r_{1}s_{k_1}r_{1})=(k_0k_1+k_0+k_1)^2.$$

\begin{center}
\includegraphics[width=8cm]{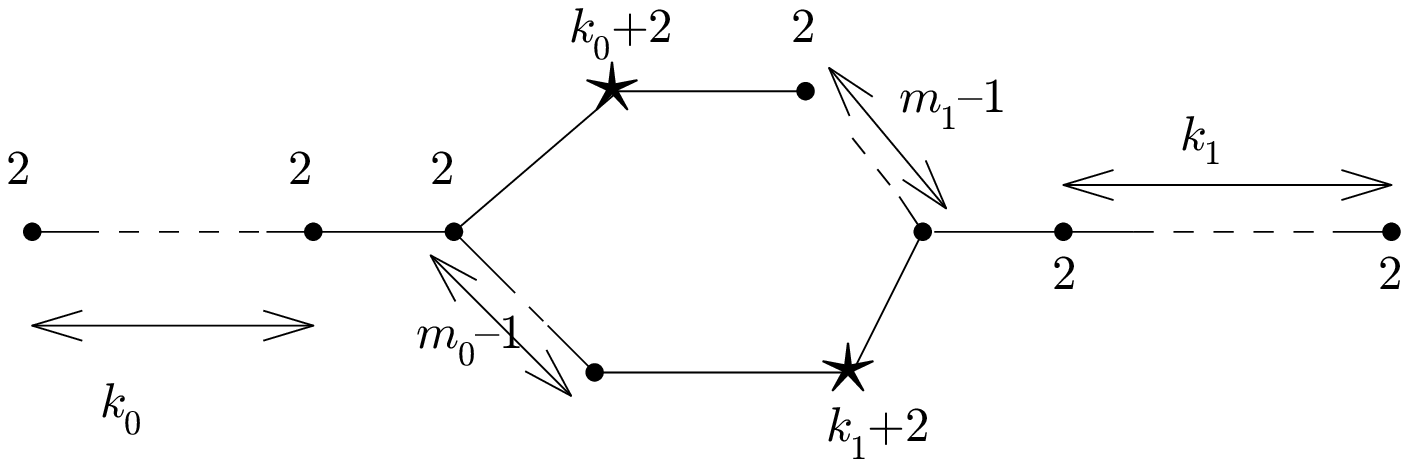}
\end{center}
\end{Ex}

\begin{Prop} Let $\frak P'_N$ be the family of polynomials $P_\s\in\bb Q[X_0,\ldots,X_{N-1}]$ such that
$$\det M(\s)=\det M(\s_0\cdots\s_{N+\r-1})=P_\s(k_0,\ldots,k_{N-1})^2,$$
with notations of (\ref{Notationsmatrices}). 
\begin{description}
\item{1)} For any $N\ge 0$, $\frak P_N=\frak P'_N$,
\item{2)} Let $\s_{i_j}=s_{k_j}$, $0\le i_j\le N+\r-1$, $0\le j\le N-1$ be the singular sequences in $\s$ and let $A\subset \bb Z/N\bb Z$ be the subset  of indices $j$ such that $\s_{i_j +1}$ is a regular sequence, then
$$P_\s=P_A.$$
\end{description}
\end{Prop}
{\bf Proof}: To prove 1) it is sufficient to check conditions i) to iv) of proposition \ref{caracterisationpol}.\\
a) Condition i) has been checked in examples \ref{Exemplesdepol} and \ref{ExemplesPsigma}. It is not possible to have two adjacent regular sequences, hence there are $2^N$ ways to insert regular sequences among $N$ singular sequences, therefore we have ii).\\
b) We suppose now that $N\ge 3$. Let $A$ be the subset (perhaps empty) of indices $j$ in $\bb Z/N\bb Z$ such that the singular sequence $\s_{i_j}=s_{k_j}$ is followed by a regular sequence. Let $\l \prod_{j\in J}X_j$, $\l\in\bb Q$, $J\subset\{0,\ldots,N-1\}$ be a monomial of $P_\s$. We are going to prove that:
\begin{center}
If $i\not\in J$ but $i-1\in J$ and $i+1\in J$, then $i\in A$.\\
\end{center}
On that purpose, suppose that $i\not\in A$, $i-1\in J$ and $i+1\in J$. Since $\det M=P_\s(k_0,\ldots,k_{N-1})^2$, it is sufficient to show that in the development of $\det M$, any term which contains the factor $(k_{i-1}k_{i+1})^2$ must also contain the factor $k_i^2$, or more simply the factor $k_i$. In view of the reduction lemma \ref{reductionlemma}, there are two possible cases:
\begin{itemize}
\item $\s$ contains the sequence $s_{k_{i-1}}r_1s_{k_i}s_{k_{i+1}}$: By theorem \ref{Arbresetcycle}, the dual graph of $M$ contains one of the two subgraphs
\begin{center}
\includegraphics[width=11cm]{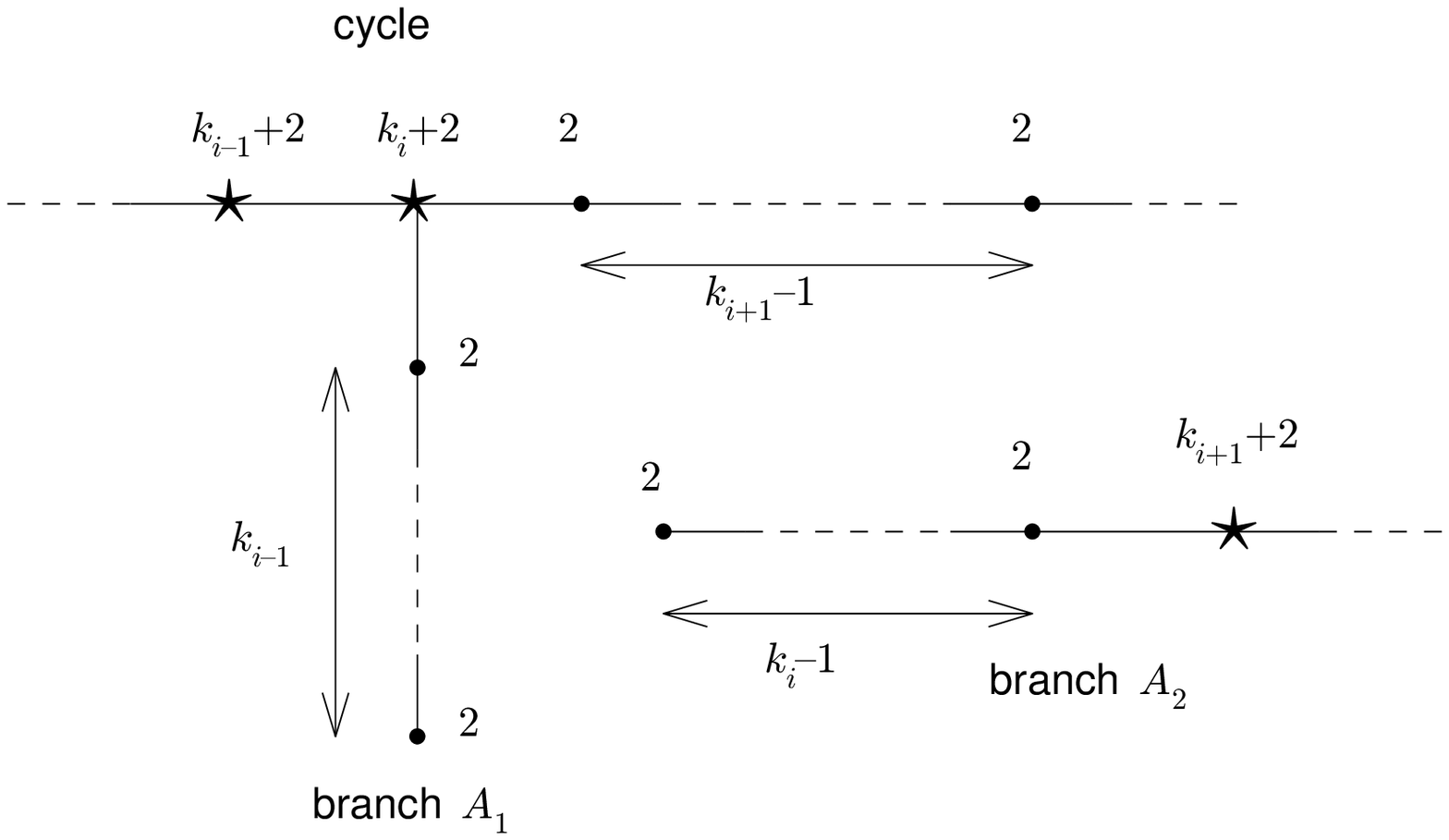}
\end{center}
\begin{center}
\includegraphics[width=11cm]{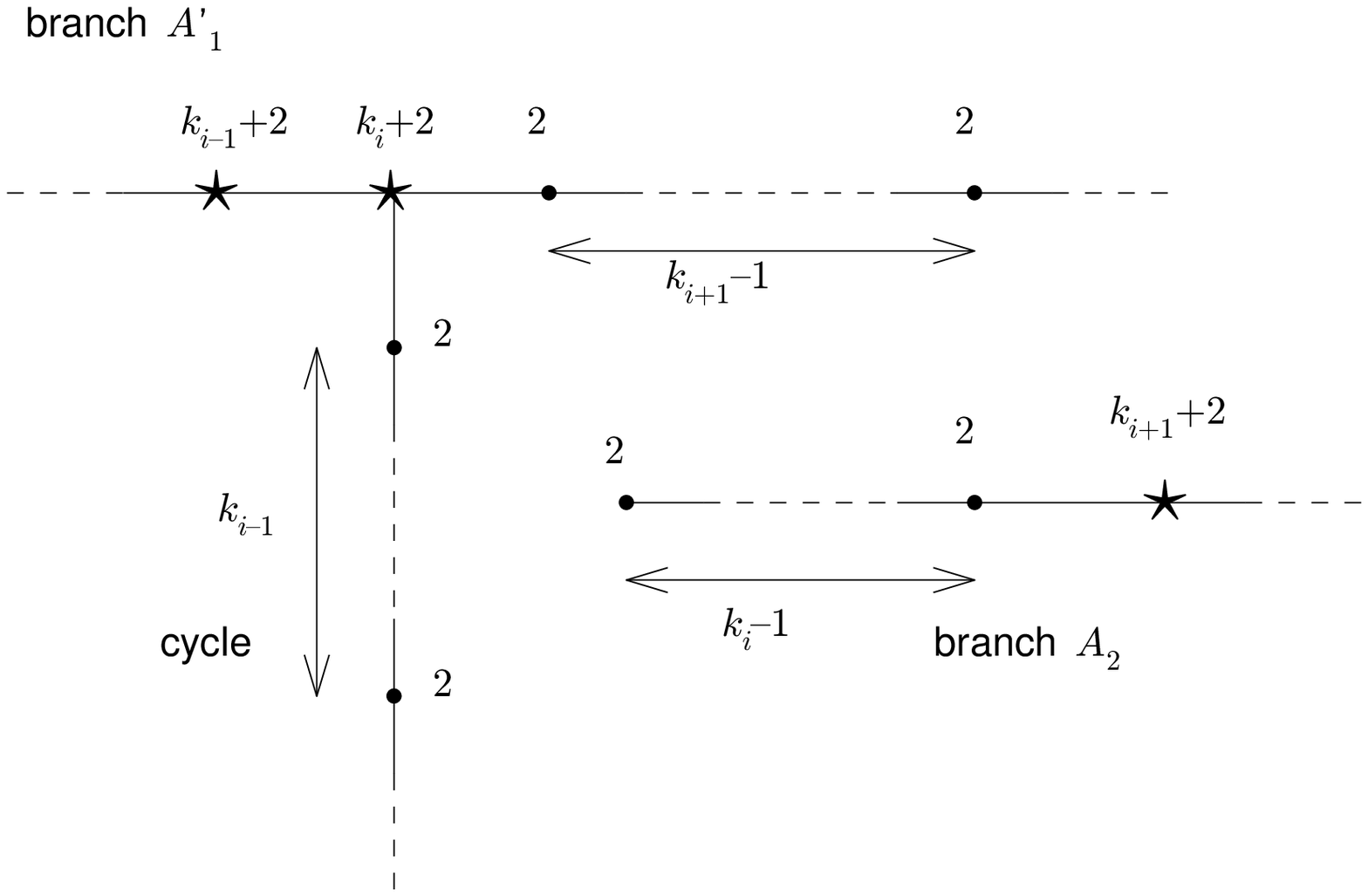}
\end{center}
Notice  that to obtain the factor $(k_{i-1}k_{i+1})^2$ one has to develop the determinant relatively to the branch $A_2$, then each term containing the factor $k_{i+1}^2$ has to contain $k_ik_{i+1}^2$. 
\item $\s$ contains the sequence $s_{k_{i-1}}s_{k_i}s_{k_{i+1}}$: a similar argument gives the result.
\end{itemize}
c)  Any monomial of $P_\s$ may be written as $\l \prod_{j\not\in C}X_j$, $\l\in\bb Q$ ($C$ is the complement of $J$ !). Suppose that $C\neq \emptyset$. Then, either $C$ contains an element of $A$, either $C$ doesn't, however by b), $C$ contains a pair $\{j,j+1\}$. We have proved that in all cases $C$ contains a generating allowed subset, hence we have the second part of iv).\\
d) In order to see that for each allowed subset $B\in \cal P_A$ we have
$$P_\s(X_i=0,\ i\in B)\in \frak P'_N$$
it is sufficient to check this property for generating allowed subsets $B$. By theorem \ref{Arbresetcycle}, 
\begin{itemize}
\item If $i\in A$, then the dual graph of $M$ contains the subgraph
\begin{center}
\includegraphics[width=5cm]{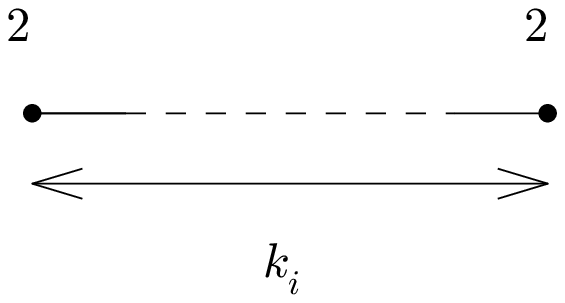}
\end{center}
with $k_i$ vertices (and not $k_i-1$). Vanishing of $k_i$ yields a configuration of a branch $A_s$ and part of cycle $\G_s$ whose parity is changed (see \ref{Arbresetcycle}).
\item If $\{i,i+1\}$ is generating allowed pair, the dual graph contains the subgraphs
\begin{center}
\includegraphics[width=11cm]{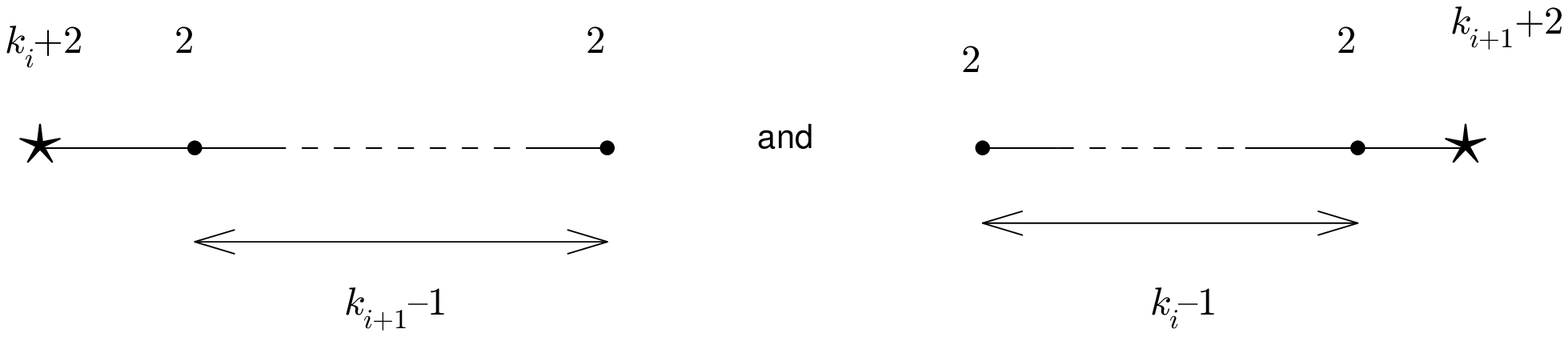}
\end{center}
Vanishing of $k_i$ and $k_{i+1}$ yields the graph of $M(\s')$, where $\s'$ is obtained from $\s$ deleting the sequences $s_{k_i}$ and $s_{k_{i+1}}$.
\end{itemize}
e) To end we have, on one hand, to compute the homogeneous parts of degrees $N$ and, on second hand, to compute the homogeneous part of degree $N-1$ of $P_\s$ to check that $P\s=P_A$ thanks to proposition \ref{Premprop}.
By reduction lemma, $\deg_{X_i} P_\s\leq 1$, hence if we show that $P_\s$ contains the monomial $\prod_{i=0}^{N-1}X_i$, it is necessarily its homogeneous part of highest degree. \\
If $A=\emptyset$, the dual graph contains one or two cycles without branches. To obtain in the developement of $\det M$ the term $(\prod_{i=0}^{N-1}k_i)^2$, it is sufficient to develop successively relatively each vertex of weight $>2$. By b), $P_\s$ contains no monomials of degree $N-1$, which gives the result in this case.\\
If $A\neq \emptyset$, we may suppose by reduction lemma, that all regular sequences are equal to $r_1$. By theorem \ref{Arbresetcycle}, all roots of the branches have weight $>2$. If we develop successively relatively to each column corresponding to a vertex of weight $>2$, we obtain:
$$\det M(\s)=\prod_{i=0}^{N-1} k_i \ \det B  + \sum_{i=0}^{N-1} \prod_{j\neq i} k_j \ \det B_i \quad {\rm mod}\quad (k_0,\ldots,k_{N-1})^{2N-2},$$
where the dual graph of $B$ is obtained from the one of $M(\s)$ by deletion of all the vertices of weights $>2$, and the dual graph of $B_i$ by deletion all the vertices of weights $>2$ but the one of  $k_i+2$ and setting the weight of the latter equal to $2$.
Now the graph of $B$ is composed of connected components which are chains of the form
\begin{center}
\hspace{3cm}\includegraphics[width=7cm]{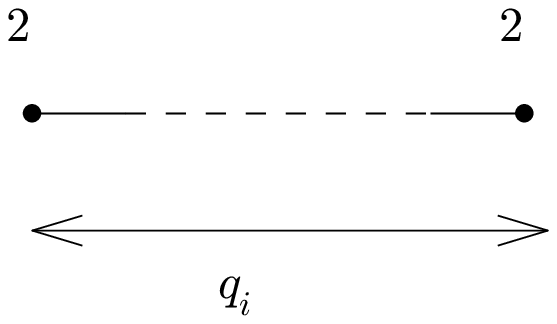}
\end{center}
where $q_i= k_i-1$ (resp. $q_i= k_i$) if   the  sequence  which follows $ s_{k_i}$   is  singular (resp. regular). Therefore the contribution of this term is
$$\prod_{i\not\in A}k_i\prod_{i\in A}(k_i+1)=\prod_{0\le i\le N-1}k_i + \sum_{i\in A} \prod_{j\neq i} k_j \quad {\rm mod}\quad (k_0,\ldots,k_{N-1})^{N-2}.$$
It remains to compute $\det B_i$: By lemma \ref{Existpoly}, $\det B_i$ is a polynomial of degree at most $N$ and we have to determine when this degree is precisely $N$. On that purpose, suppose that the index $i$ corresponds to a vertex between two chains of vertices of weight $2$, that is to say we have a subgraph
 \begin{center}
\includegraphics[width=8cm]{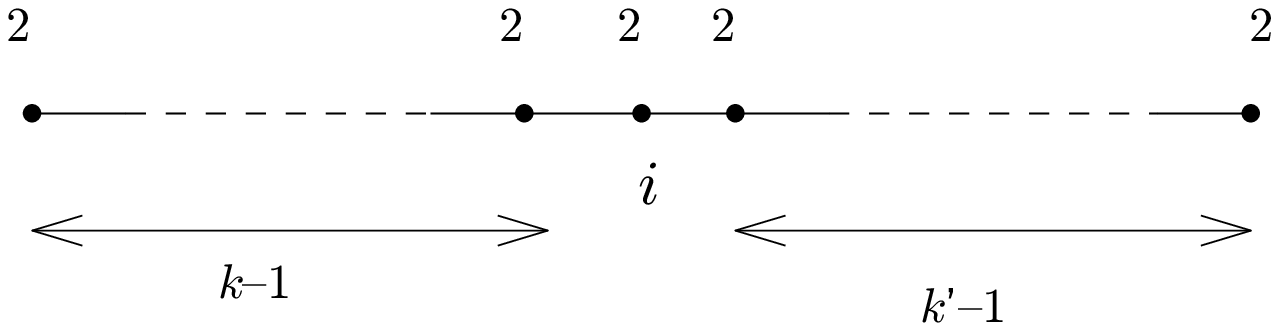}
\end{center}
By lemma \ref{easycase} the determinant of this connected component is $k+k'$, 
hence of degree $1$ and $\det B_i$ will be of degree at most $N-1$. Therefore we are only interested in vertices which are the root of a branch or linked to a root. By theorem  \ref{Arbresetcycle}, for $\r(S)\ge 1$ and $t=s+1$ mod $\r(S)$ there are four possible situations:\\
 \begin{center}
\includegraphics[width=12cm]{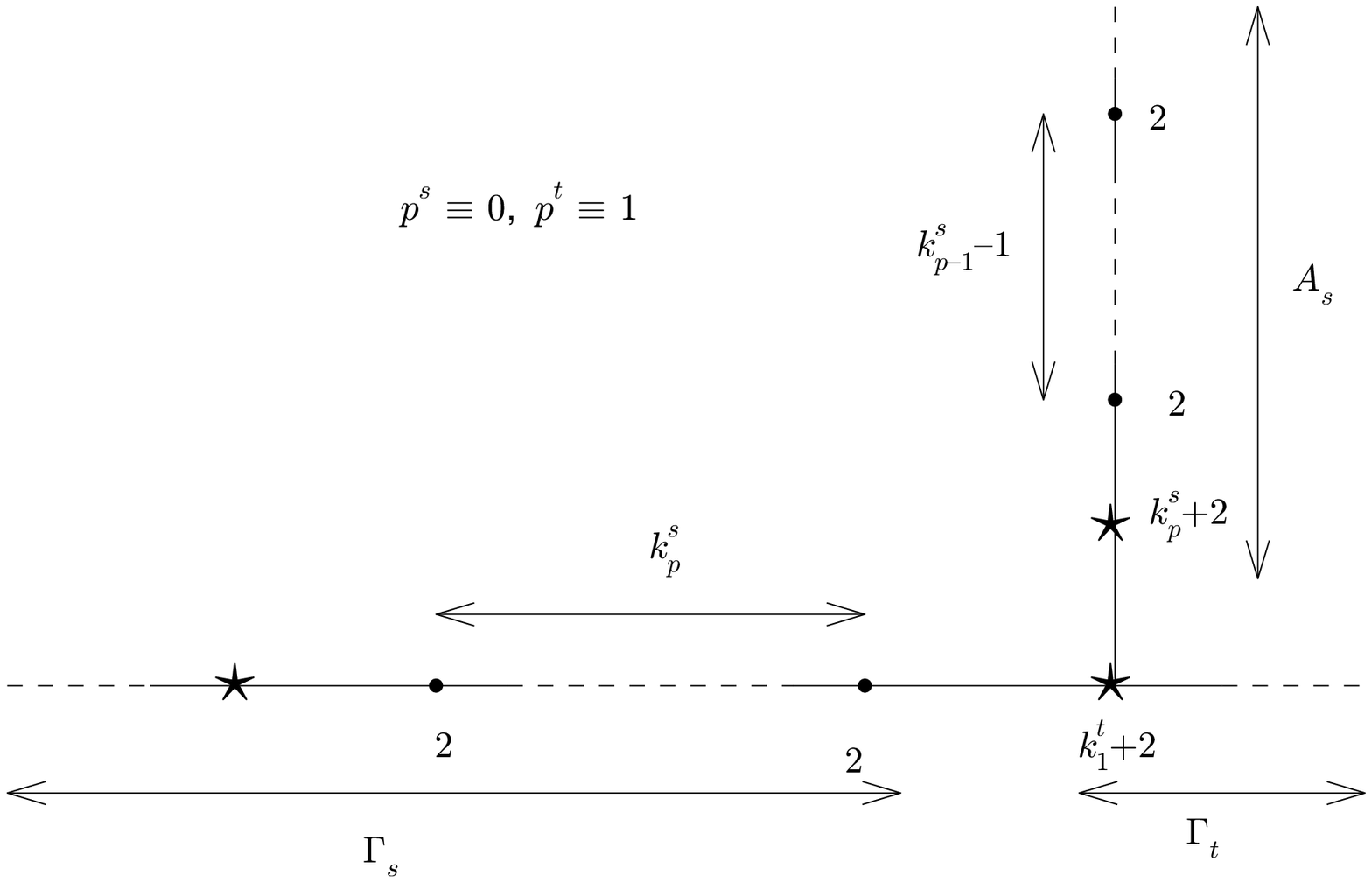}
\end{center}
\begin{center}
\includegraphics[width=6cm]{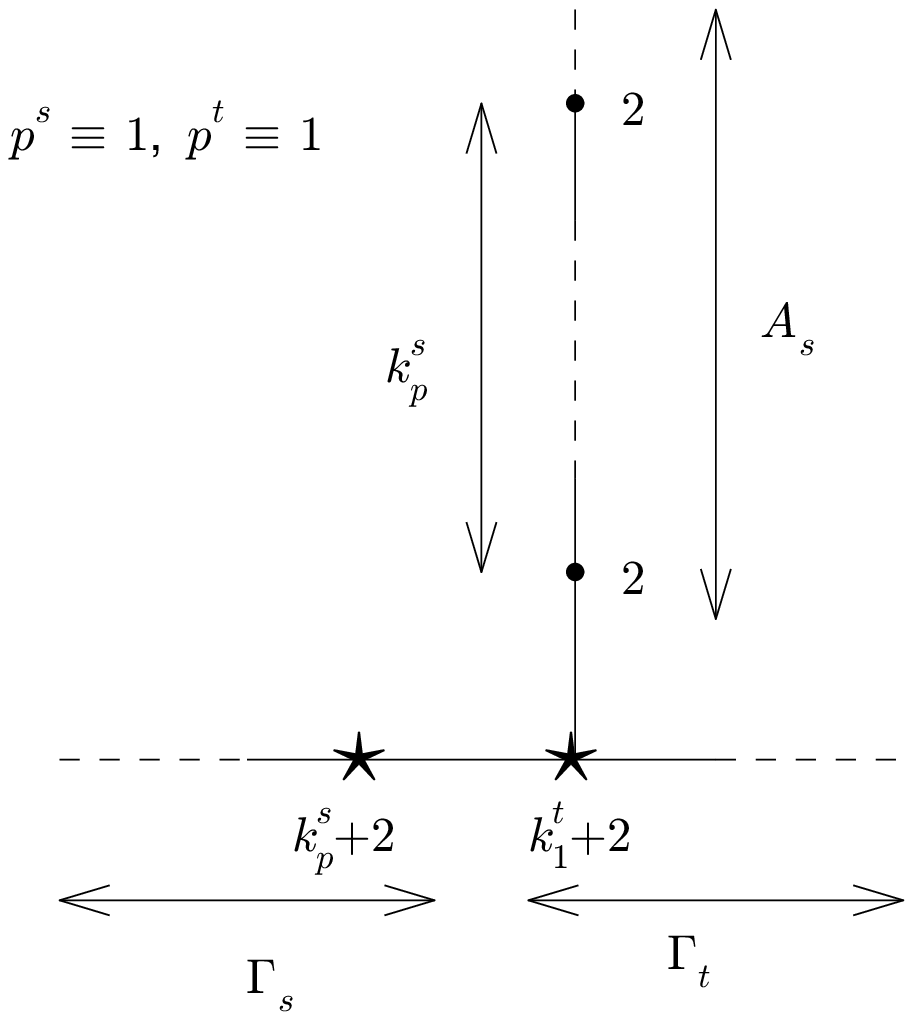}
\end{center}
 \begin{center}
\includegraphics[width=10cm]{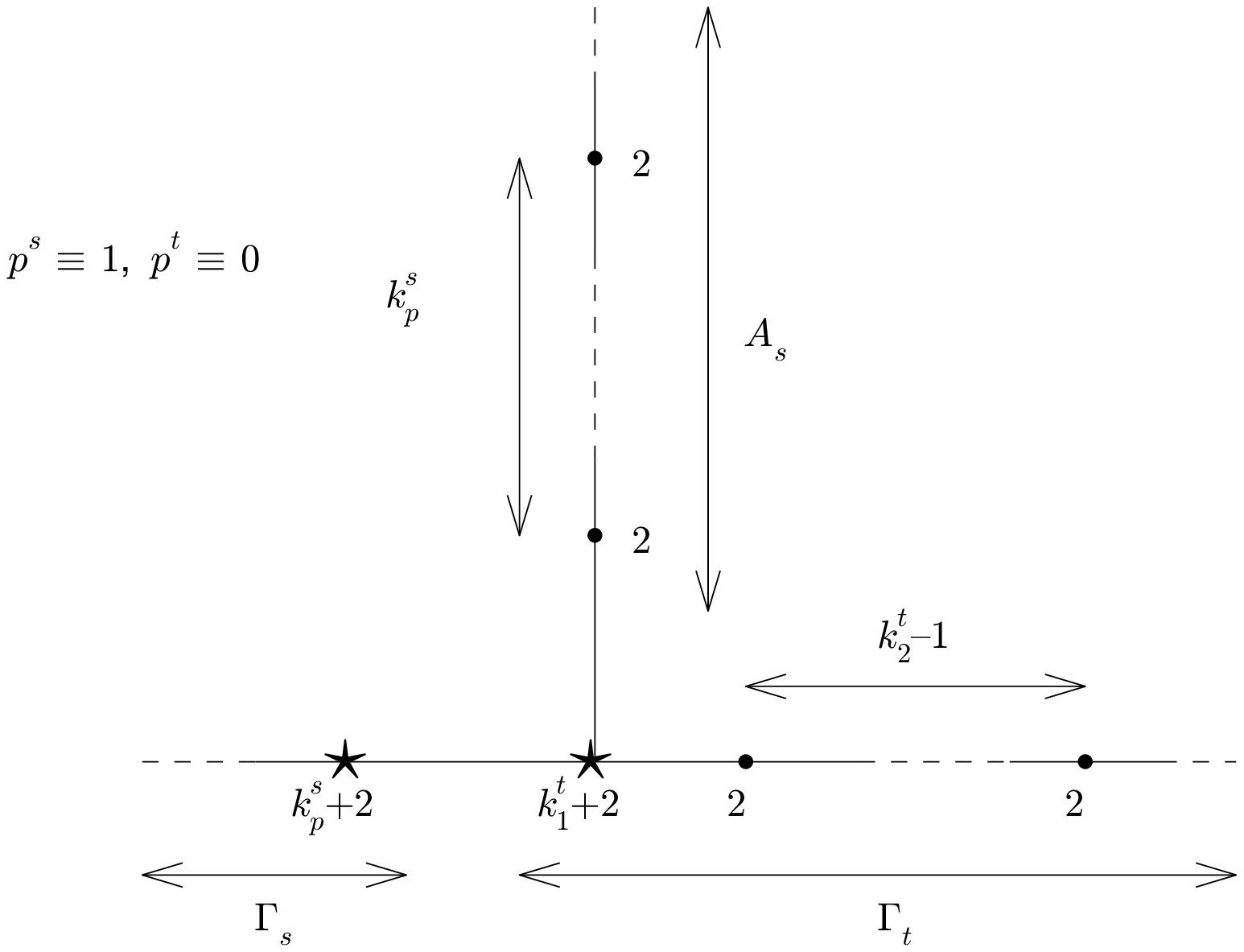}
\end{center}
 \begin{center}
\includegraphics[width=9cm]{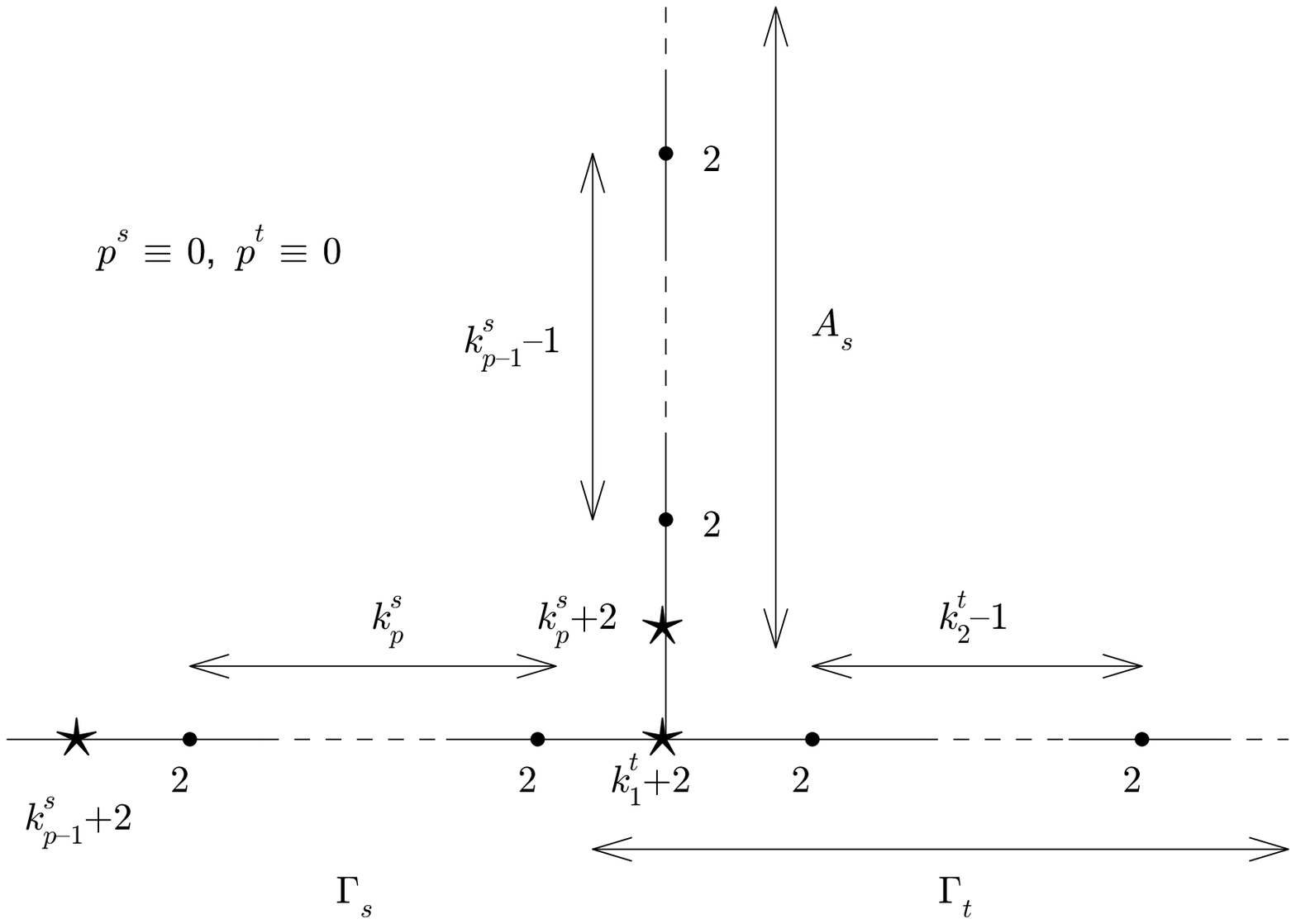}
\end{center}
We see that the two only involved vertices are thoose of weight $k_p^s+2$ and $k_1^t+2$.
\begin{itemize} 
\item If $s_{k_i}$ is followed by a regular sequence, i.e. $s_{k_i}=s_{k_p^s}$ or ($s_{k_i}=s_{k_1^t}$ and $p^t=1$): 
\begin{description}
\item In the first case the graph of $B_i$ contains one of the subgraphs
   \begin{center}
     \includegraphics[width=11cm]{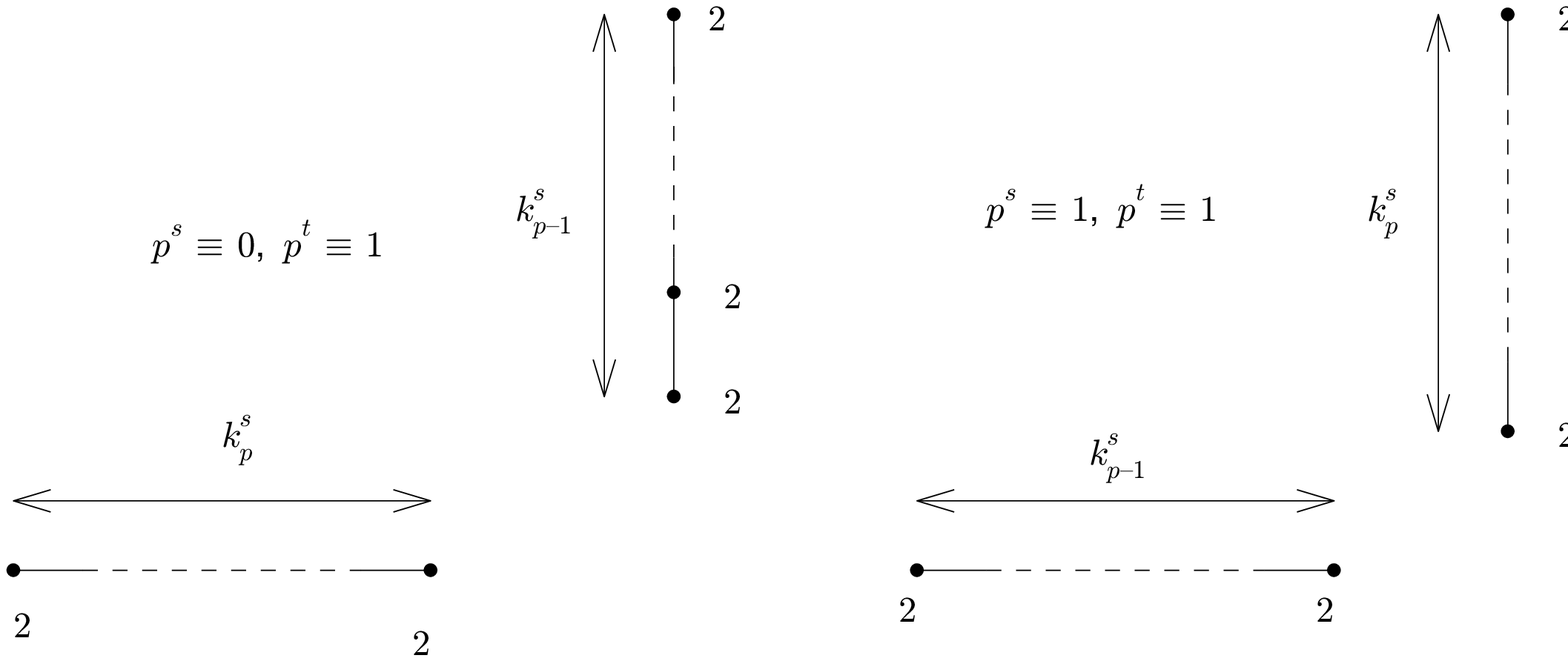} 
    \end{center} 
      \begin{center}
    \includegraphics[width=8cm]{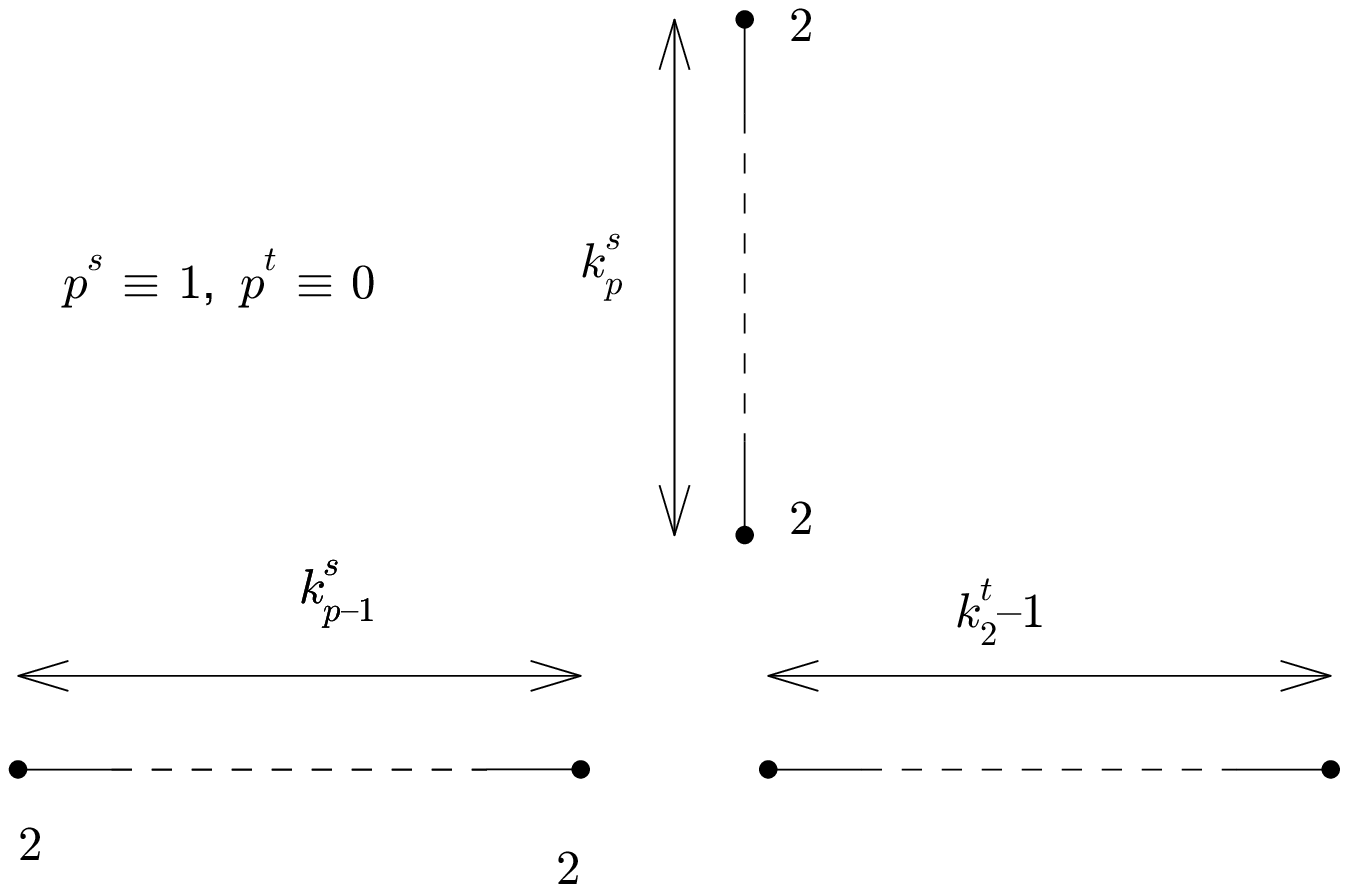}
    \end{center} 
   \begin{center}
    \includegraphics[width=10cm]{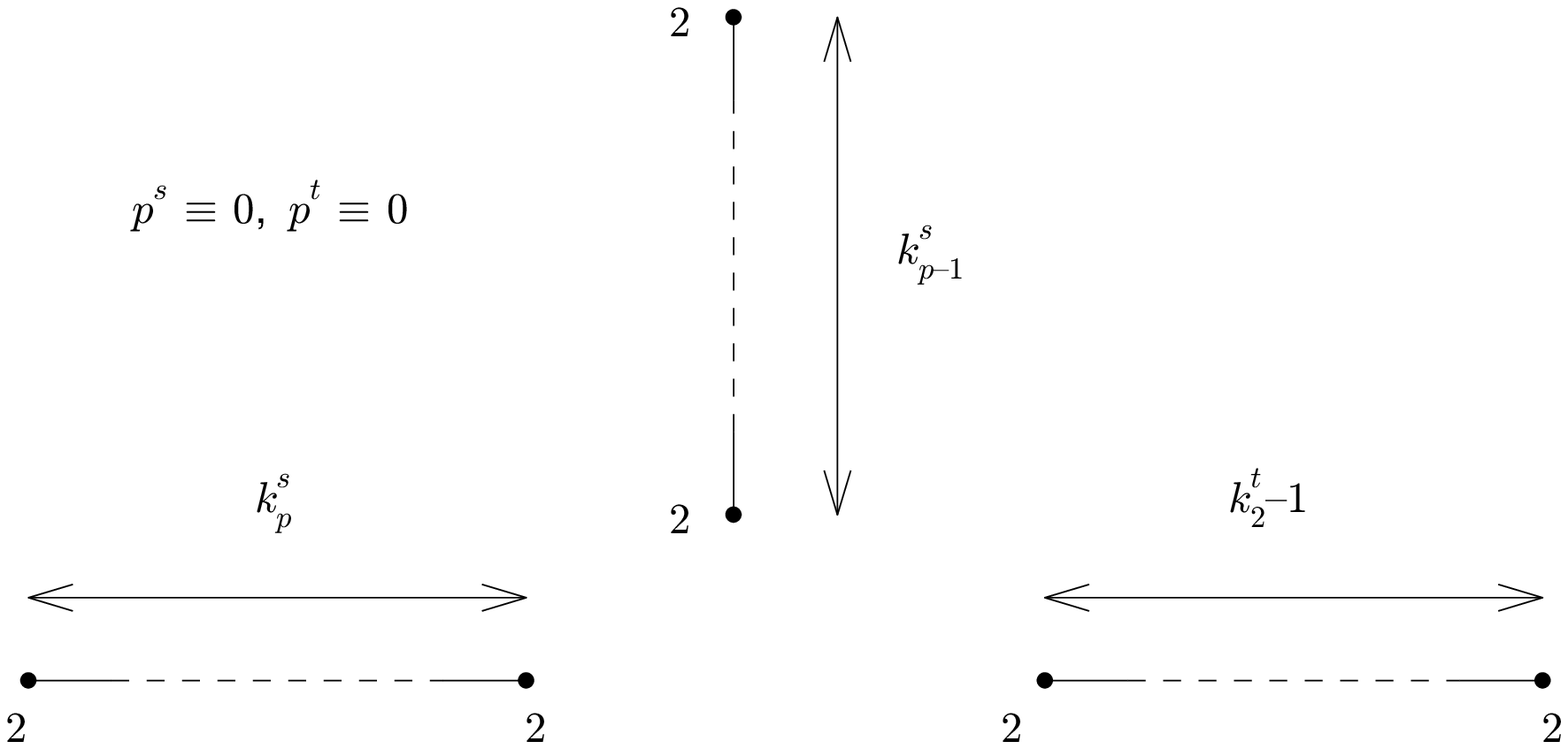}
    \end{center} 
    \item If $p^t=1$ and  $p^s$ is any integer, we have the following  connected component 
    \begin{center}
    \includegraphics[width=4cm]{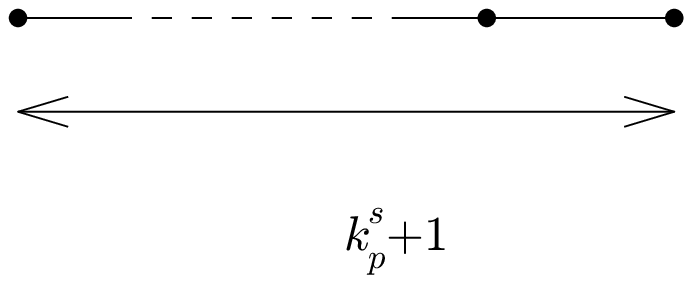}
    \end{center} 
\end{description}
in all these cases $\deg B_i=N$ with contribution $\prod_{i=0}^{N-1}k_i$.

\item If  $s_{k_i}$ is followed by a singular sequence, i.e. $s_{k_i}=s_{k_1^t}$ and (if $p^t\equiv 1$ then $p^t\ge 3$):  the dual graph contains the following connected components
 \begin{center}
    \includegraphics[width=11cm]{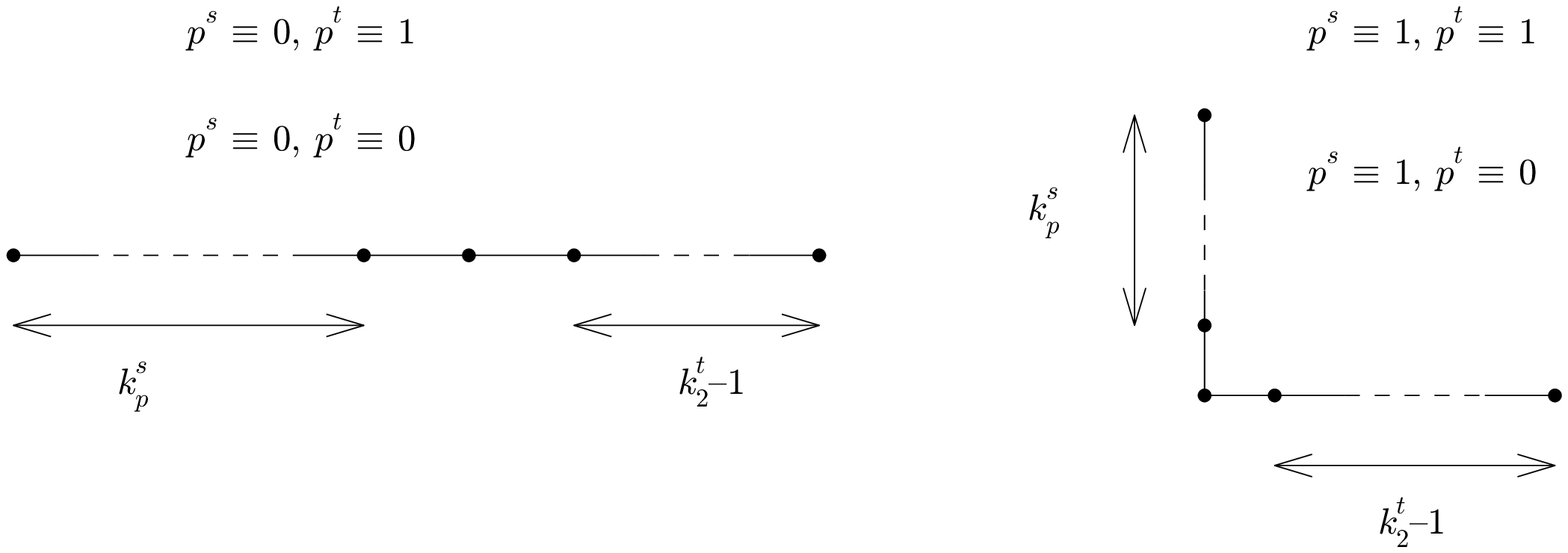}
    \end{center} 
In all these cases, $\deg \det B_i    =N-1$.
\end{itemize} 
 Finally, we have
 $$\begin{array}{lcl}
 P_\s(k_0,\ldots,k_{N-1})^2 &=& \det M(\s)\\
 & = &\dps\prod_{i=0}^{N-1}k_i \left(\prod_{i=0}^{N-1}k_i + \sum_{i\in A} \prod_{j\neq i} k_j\right) + \sum_{i\in A} \prod_{j\neq i}k_j \prod_{i=0}^{N-1}k_i \\
 &&\\
 & = &\dps \left(\prod_{i=0}^{N-1}k_i\right)^2 + 2  \sum_{i\in A} \prod_{j\neq i}k_j \prod_{i=0}^{N-1}k_i \\
&&\\
 &&\quad {\rm mod}\quad (k_0,\ldots,k_{N-1})^{2N-2}
 \end{array}$$
 and $P_\s=P_A$ by proposition \ref{Premprop}, 3) as wanted.\hfill $\Box$

Georges Dloussky, UMR6632 CNRS\\
Centre de Mathématiques et d'Informatique, 
Universit\'e d'Aix-Marseille 1,\\
 39, rue F.Joliot-Curie, 13453 Marseille Cedex 13, France\\
dloussky@cmi.univ-mrs.fr


\begin{thebibliography}{AUT}

\bibitem{D1} {\sc Dloussky G.}: Structure des surfaces de
Kato. {\em M\'emoire de  la S.M.F 112.$\rm n^{\circ}14$ (1984).} 
\bibitem{D2germe} {\sc Dloussky G.}: Sur la classification des germes d'applications holomorphes contractantes. {\em Math. Ann. 280 (1988), 649-661.}
\bibitem{D2} {\sc Dloussky G.}: Une construction \'el\'ementaire des
surfaces d'Inoue-Hirzebruch. {\em Math. Ann. 280, (1988), 663-682.} 
\bibitem{DO} {\sc Dloussky G., Oeljeklaus K.}: Vector fields and 
foliations on compact surfaces of class VII$_0$. {\em Ann. Inst. Fourier 49 (1999), 1503-1545.}
\bibitem{DO2} {\sc Dloussky G., Oeljeklaus K.}: Surfaces de la classe VII$_0$ et automorphismes de Hénon. {\em C.R.A.S. 328, série I, p.609-612, 1999.}

\bibitem{E} {\sc Enoki I.}: Surfaces of class VII$_0$ with curves. 
{\em T\^ohoku Math. J. 33,
(1981), 453-492.} 
\bibitem{F} {\sc Favre, Ch.}: Classification of $2$-dimensional
contracting rigid germs,   
{\em Jour. Math. Pures Appl. 79, (2000), 475-514.}
\bibitem{F2} {\sc Favre, Ch.}: Dynamique des applications rationnelles. 
{\em Thèse pour le grade de Docteur en Sciences. Université de Paris XI Orsay (2000). http://tel.archives-ouvertes.fr/tel-00003577/fr/}
\bibitem{H}{\sc Hirzebruch F.}:  Hilbert modular surfaces. {\em L'enseignement Math. 19 (1973), 183-281.}
\bibitem{I} {\sc Inoue M.}: New surfaces with no meromorphic 
functions II.
 {\em Complex Analysis and Alg. Geom. 91-106. Iwanami Shoten Pb. 
1977.} 
\bibitem{KA}{\sc Karras U.} Deformations of cusps singularities. {\em Proc. of Symp. in pure Math. 30, 37-44, AMS, Providence.}
\bibitem{K} {\sc Kato Ma.} Compact complex manifolds containing ``global spherical shells'' I {\em Proc. of the Int. Symp. Alg. Geometry, Kyoto
(1977) Iwanami Shoten Publ.}
\bibitem{KO}{\sc Kodaira K.} On the structure of compact complex 
analytic 
surfaces I, II.  {\em Am. J. of Math.
vol.86, 751-798 (1964); vol.88, 682-721 (1966)}.

\bibitem{L2}{\sc  Laufer H.} On minimally elliptic singularities. 
{\em Amer. J. of Math. 99, 
p1257-1295, (1977).}
\bibitem{LW}{\sc Looijenga E. \& Wahl J.} Quadratic functions and smoothing surface singularities. {\em Topology 25 (1986), 261-291.}
\bibitem{MER} {\sc M\'erindol J.Y.} Surfaces normales dont le 
faisceau 
dualisant est trivial. {\em C.R.A.S. 293,
417-420 (1981)}.
\bibitem{N1} {\sc Nakamura I.}: On surfaces of class $\rm VII_0$ with
curves. {\em Invent. Math. 78,(1984), 393-443.} 
\bibitem{N2} {\sc Nakamura I.}: On surfaces of class $\rm VII_0$ with Global Spherical Shells. {\em Proc. of the Japan Acad. 59, Ser. A, No 2 (1983), 29-32}
\bibitem{N3}{\sc Nakamura I.}  On the equations 
$x^{p}+y^{q}+z^{r}-xyz=0$.{\em Advanced Studies in pure
Math. 8, Complex An. Singularities, 281-313 (1986). }
\bibitem{N4}{\sc Nakamura I.} Inoue-Hirzebruch surfaces and a duality of hyperbolic unimodular singularities I. {\em Math. Ann. 252, 221-235 (1980).}
\bibitem{PIN} {\sc Pinkham H.}: Singularit\'es rationnelles de 
surfaces.{\em  Appendice. S\'eminaire sur les
singularit\'es des surfaces. Lecture Notes  777 . Springer-Verlag 
1980}.
\bibitem{R}{\sc Ribenboim R.} Polynomials
whose values are powers. {\em J. f\H ur die reine und ang. Math. 268/269, 34-40
(1974)}.
\bibitem{SA}{\sc Sakai F. }Enriques classification of normal 
Gorenstein 
surfaces.{\em Am. J. of Math. 104, 1233-1241 (1981)}.
\bibitem{SE}{\sc Serre, J.P} Cours d'arithmétique
{\em Presses Universitaires de France (1970)}




\end{thebibliography}
\end{document}